\documentclass[letterpaper]{article}
\usepackage{authblk}
\usepackage{url}
\usepackage[margin=1in,footskip=0.25in]{geometry}
\usepackage{amssymb}
\usepackage{graphicx}
\usepackage{amsmath}
\usepackage{algorithmic} 
\usepackage{algorithm}
\usepackage{float}
\usepackage{lipsum}
\usepackage{subfigure}
\usepackage{wrapfig}
\usepackage{color}
\usepackage{easylist}
\usepackage{lineno}
\usepackage{stmaryrd}
\usepackage{booktabs}
\usepackage{cleveref}
\usepackage{multirow}
\usepackage{hhline}

\def \grad{\nabla}
\def \half{\frac{1}{2}}
\def \p   {\partial}

\newcommand   \D [2]{\frac{\partial     #1}{\partial   #2   }}

\usepackage{blindtext}

\def \II{\mathbb{I}}

\def \F{\vec{F}}

\def \N{\vec{N}}

\def \b{\vec{b}}
\def \t{\vec{t}}
\def \U{\vec{U}}
\def \V{\vec{V}}
\def \X{\vec{\chi}}
\def \Y{\vec{\xi}}
\def \b{\vec{b}}
\def \f{\vec{f}}

\def \n{\vec{n}}
\def \s{\vec{X}}
\def \f{\vec{f}}
\def \u{\vec{u}}

\def \x{\vec{x}}
\def \grad{\nabla}
\def \half{\frac{1}{2}}
\def \p{\partial}

\def \Dx{{\mathrm d} \x}

\def \DA{{\mathrm d} A}

\newcommand{\volume}{{\ooalign{\hfil$V$\hfil\cr\kern0.08em--\hfil\cr}}}

\def \dt{\Delta t}

\def \cI{\mathcal I}
\def \cJ{\vec{\mathcal J}}

\def \cS{\vec{\mathcal S}}
\def \cT{\mathcal T}

\def \mfac{M_\text{fac}}
\def \CL{C_\text{L}}
\def \CD{C_\text{D}}

\def \Re{\textrm{Re}}
\def \St{\textrm{St}}
\def \Lwake{L_\text{wake}}
\def  \thetas{\theta_\text{s}}
\def \L2{L^2}
\def \Linf{L^{\infty}}

\def \grad{\nabla}
\def \half{\frac{1}{2}}
\def \p   {\partial}
\def \Mfac {M_{\text{fac}}}

\renewcommand{\vec}[1]{\ensuremath\boldsymbol{#1}}

\newcommand*{\rttensor}[1]{\overline{\overline{#1}}}

\newcommand{\sspan}{\mathop{\mathrm{span}}}

\title{An Immersed Interface Method for Discrete Surfaces}

\author[1,2]{Ebrahim M.~Kolahdouz} 
\author[3]{Amneet Pal Singh Bhalla}
\author[2]{Brent A.~Craven}
\author[4,5,6]{Boyce E.~Griffith} 
\affil[1]{Department of Mathematics, University of North Carolina, Chapel Hill, NC, USA}
\affil[2]{Division of Applied Mechanics, Office of Science and Engineering Laboratories, Center for Devices and Radiological Health, United States Food and Drug Administration, Silver Spring, MD, USA}
\affil[3]{Department of Mechanical Engineering, San Diego State University, San Diego, CA, USA}
\affil[4]{Departments of Mathematics, Applied Physical Sciences, and Biomedical Engineering, University of North Carolina, Chapel Hill, NC, USA}
\affil[5]{Carolina Center for Interdisciplinary Applied Mathematics, University of North Carolina, Chapel Hill, NC, USA}
\affil[6]{McAllister Heart Institute, University of North Carolina, Chapel Hill, NC, USA}
\affil[ ]{\texttt{ebrahimk@email.unc.edu} and \texttt{boyceg@email.unc.edu}}

\begin{document}

\maketitle

\begin{abstract}
Fluid-structure systems occur in a range of scientific and engineering applications.
The immersed boundary (IB) method is a widely recognized and effective modeling paradigm for simulating fluid-structure interaction (FSI) in such systems, but a difficulty of the IB formulation of these problems is that the pressure and viscous stress are generally discontinuous at fluid-solid interfaces.
The conventional IB method regularizes these discontinuities, which typically yields low-order accuracy at these interfaces.
The immersed interface method (IIM) is an IB-like approach to FSI that sharply imposes stress jump conditions, enabling higher-order accuracy, but prior applications of the IIM have been largely restricted to numerical methods that rely on smooth representations of the interface geometry.
This paper introduces an immersed interface formulation that uses only a $C^{0}$ representation of the immersed interface, such as those provided by standard nodal Lagrangian finite element methods.
Verification examples for models with prescribed interface motion demonstrate that the method sharply resolves stress discontinuities along immersed boundaries while avoiding the need for analytic information about the interface geometry.
Our results also demonstrate that only the lowest-order jump conditions for the pressure and velocity gradient are required to realize global second-order accuracy.
Specifically, we demonstrate second-order global convergence rates along with nearly second-order local convergence in the Eulerian velocity field, and between first-~and second-order global convergence rates along with approximately first-order local convergence for the Eulerian pressure field.
We also demonstrate approximately second-order local convergence in the interfacial displacement and velocity along with first-order local convergence in the fluid traction along the interface.
As a demonstration of the method's ability to tackle more complex geometries, the present approach is also used to simulate flow in a patient-averaged anatomical model of the inferior vena cava, which is the large vein that carries deoxygenated blood from the lower extremities back to the heart.
Comparisons of the general hemodynamics and wall shear stress obtained by the present IIM and a body-fitted discretization approach show that the present method yields results that are in good agreement with those obtained by the body-fitted approach.

\end{abstract}

\noindent \textbf{Keywords:}	immersed boundary method, immersed interface method, finite element, fluid-structure interaction, jump conditions, complex geometries

\section{Introduction and overview}
\label{sec:intro}

Stable and accurate yet simple and computationally tractable methods for treating interfaces are of great importance in fluid dynamics problems involving immersed boundaries.
Body-fitted discretization approaches, which approximate the fluid dynamics on a computational domain that conforms to the interface geometry, can provide high accuracy for such problems but are limited by the continued difficulties of mesh generation, especially for problems involving substantial motion of the immersed interface. 
A widely used alternative approach to these types of interactions is the immersed boundary (IB) method introduced by Peskin \cite{peskin1972flow,peskin1977}.
Unlike body-fitted methods, the IB approach avoids dynamic mesh regeneration and allows for the use of fast structured-grid fluid solvers.
In the classical IB method, the Navier-Stokes equations are solved on a fixed background Eulerian grid, and the immersed structure is represented as a collection of Lagrangian markers.
The interaction of Eulerian and Lagrangian frameworks takes place in two steps: 1) spreading structural forces from the Lagrangian markers to the Eulerian grid using a regularized Dirac delta function and 2) interpolating velocities from the Eulerian grid back to the Lagrangian markers using the same smoothed delta function.

A limitation of the conventional IB method is that it yields low-order accuracy at fluid-structure interfaces.
The continuous form of the IB equations uses integral equations with singular kernels to connect the Lagrangian and Eulerian frames, and this formulation is equivalent to a formulation involving jump conditions \cite{lai2001remark}.  
However, the use of regularized delta functions in the conventional numerical realization of the IB method has the effect of regularizing stress jumps at those interfaces, which implies that stresses do not converge pointwise at the interface.
Much work to improve the IB method has focused on improving its accuracy while retaining the original method's use of nonconforming grids.
The first formally second-order version of the method was introduced by Lai and Peskin \cite{lai2000}.
This method was formulated to be second-order accurate if applied to a hypothetical problem in which the regularized delta functions are replaced by fixed smooth functions, independent of the mesh \cite{lai2000}.
This method was refined by Griffith and Peskin and applied to a specific FSI problem for which the formally second-order accurate method was able to achieve its designed order of accuracy \cite{griffith2005order}, although for general FSI problems, the formally second-order accurate method still only realizes first-order convergence rates.
This approach was further extended by Griffith et al.~\cite{griffith2007adaptive} to use adaptive mesh refinement. 
This methodology has enabled modeling in a number of application areas, including cardiac dynamics \cite{
	griffith2009simulating,
	griffith2012immersed,
	XYLuo12,
	HGao14-iblv_mi,
	VFlamini16-aortic_root,
	WWChen16-coupled_lv,
	HGao17-coupled_mv_lv,
	hasan2017image,
	LYFeng18-MV},
platelet adhesion \cite{TSkorczewski14-platelet_adhesion},
esophageal transport \cite{
	kou2015fully,
	WKou17-continuum_eso,
	WKou18-abnormalities},
heart development \cite{
	LWaldrop16-peristalsis,
	NABattistaXX-heart_development},
insect flight \cite{
	SKJones15-lift_vs_drag,
	santhanakrishnan2018flow}, and
undulatory swimming \cite{
	SAlben13,
	HNguyen14,
	APSBhalla14-EHD,
	EDTytell14,
	RBale15-evolution,
	APHoover17-active_jellyfish,
	NNangia17-optimal}.

Meanwhile, motivated by the same objective, a number of sharp interface methods have also been developed.
Among them are
the immersed interface method (IIM) \cite{leveque1994immersed,li2001immersed,lee2003immersed},
the ghost-fluid method \cite{fedkiw1999non,fedkiw2002coupling},
the cut-cell method \cite{berger1989adaptive, berger1998aspects,clarke1986euler,ye1999accurate,barad2009adaptive},
the hybrid Cartesian-immersed boundary method \cite{mohd1997simulations,fadlun2000combined},
and the curvilinear immersed boundary method \cite{borazjani2008curvilinear,gilmanov2005hybrid}.
Many of these methods achieve higher-order accuracy by adopting approaches that are similar to body-fitted discretization methods, such as local modifications to the finite difference stencils, to allow for the accurate reconstruction of boundary conditions in the vicinity of the immersed interface.

This paper introduces a new immersed interface scheme for the incompressible Navier-Stokes equations that uses only a $C^{0}$-continuous representation of the interface geometry.
Immersed interface methods resemble the classical IB method in that the motion is derived from a singular surface force.
In the classical IB method and other continuum surface force approaches \cite{chang1996level,cortez2000blob}, these forces are smoothed out over a small region near the interface using a smooth kernel function.
In the IIM, these forces are used to determine discrete jump conditions that are imposed in the finite difference discretization of the fluid equations.
In fact, in most modern IIMs for FSI, generalized Taylor series expansions are used to extend the jump values from the interface to the finite difference discretization \cite{xu2006systematic}.
In this work, we view this procedure as constructing a discrete force-coupling operator that is tailored to a particular fluid discretization method.

The IIM was introduced by Leveque and Li \cite{leveque1994immersed} for the solution of elliptic PDEs with discontinuous coefficients, or in the presence of singular forces.
This initial method was extended to the two-dimensional solution of the incompressible Stokes \cite{leveque1997immersed,li2007augmented} and Navier-Stokes \cite{li2001immersed, lee2003immersed, le2006immersed} equations, typically in combination with geometrical representations of the interface based on level set methods \cite{hou1997hybrid, sethian2000structural, xu2006level}.
The systematic derivation of the jump conditions for the velocity and its first and second normal derivatives, as well as the pressure and its first normal derivative at the interface, was detailed by Lai and Li \cite{lai2001remark}.
They then coupled their method with a second-order accurate projection algorithm to solve the full Navier-Stokes equations in two spatial dimensions and empirically demonstrated that the selected jump conditions were adequate to achieve full second-order accuracy for the velocity and nearly second-order accuracy for the pressure in the maximum norms \cite{li2001immersed}.
Later, Xu and Wang used a generalized surface parameterization for the interface representation and derived jump conditions for the first-, second-, and third-order spatial derivatives of velocity and pressure along with jump conditions in the first-~and second-order time derivatives of the velocity \cite{xu2006systematic}.
They used their method in both two-dimensional \cite{xu2006} and three-dimensional \cite{xu20083d} applications.
In the infinity norm, second-order accuracy of both velocity and pressure was demonstrated in two spatial dimensions, although some deviations from the designed order of accuracy were observed in empirical tests.
In the three-dimensional cases, nearly second-order accuracy was shown in the infinity norm of the velocity, and between first-~and second-order accuracy was observed for the pressure.  
Local truncation error analysis suggests that achieving pointwise second-order accuracy in solving an elliptic PDE with a solution that includes a discontinuity along an internal interface requires correction terms up to the third normal derivative \cite{li2006immersed, xu2006}.
However, several authors have reported empirical results demonstrating local second-order accuracy when including jumps only up to the second normal derivative, provided that all the spatial terms in the original equation are approximated to second-order accuracy \cite{li2001immersed, le2006immersed,beale2007accuracy}.  
Because the interfacial discontinuities are one dimension lower than the solution domain, prior work has also studied how these conditions are sufficient despite the truncation error being reduced to first-order at the interface \cite{li2006immersed, beale2007accuracy}.
The textbook of Li and Ito \cite{li2006immersed} provides additional details on the IIM, which is now routinely used to simulate various physical phenomena \cite{hou1997hybrid, lombard2004numerical, bergou2007passive, jayathilake2010deformation, el2012computational,kolahdouz2015numerical}.


Almost all previous work using the IIM represents the immersed boundary as a smooth interface, and in many cases, it is assumed that an analytic description of the interface geometry is available.
An exception is the work of Xu and Wang \cite{xu20083d}, in which triangular patches are used to find the intersections between the interface and the finite difference stencils.
However, in that work, the boundary parameterization still relies on an analytic description of the interface for computing geometrical quantities \cite{xu20083d}, which limits the applicability of that method to situations with more general boundary geometries.

The present numerical method does not use analytical information about the interface geometry.
Instead, all geometrical information is provided by the discrete finite element representation of the interface.
The ability to use such representations facilitates simulations with complex interfacial geometries.
However, because even the lowest-order jump conditions require geometrical quantities that are discontinuous on a $C^0$ representation (e.g.~the normal vector to the surface), the use of these representations requires a substantial extension of prior IIM formulations.
In the present work, an $L^2$ projection is used to construct jump condition values that are continuous along the $C^0$ interface.
An $L^2$ projection is also used to construct a velocity interpolation scheme that accounts for the known velocity discontinuities along the interface.


Empirical results from verification tests demonstrate that the method attains second-order global convergence rates along with nearly second-order local convergence in the Eulerian velocity field, and between first-~and second-order global convergence rates along with approximately first-order local convergence for the Eulerian pressure field.
These tests also show that along the interface, the method yields approximately second-order local accuracy in the velocity along with first-order local accuracy in the fluid traction (pressure and wall shear stress).
In this work, we impose only the lowest order jump conditions in the pressure and velocity gradient.
To our knowledge, all prior work using the IIM that has achieved similar levels of accuracy has imposed additional higher-order jump conditions.
We also show that this approach yields a substantial improvement in accuracy as compared to a conventional formally second-order accurate IB method, and that it is necessary to use both pressure and velocity jump conditions along with corrected velocity interpolation for the present immersed interface scheme to achieve its full accuracy.
As a demonstration of the method's versatility in treating complex geometries, this work also presents initial results of IIM simulations of flow in a patient-averaged anatomical model of the human inferior vena cava, which is the major vein that is responsible for returning deoxygenated blood from the lower extremities  back to the heart.

\section{Mathematical formulation}
\label{sec:formulation}

This section introduces the governing equations and physical interface conditions in the continuous formulation of the IB method.  
The alternative treatment of singular forces in the form of jump conditions for the pressure and velocity derivatives is then described.
The resulting jump conditions will be used in the discretization technique described in Sec.~\ref{sec:methodology}.
Throughout the paper, Eulerian quantities are indicated using lower case variables, and Lagrangian quantities are indicated by upper case variables.
Bold face variables are used for vectors, and bold double-bar symbols are used for tensors.
Superscripts `$+$' (`$-$') indicate values obtained from the `exterior' (`interior') side of an interface.
 
\subsection{Governing equations}
\label{subsec:formulation-GE}

\begin{figure}[t!!]
	\centering
	\includegraphics[width=0.45\textwidth]{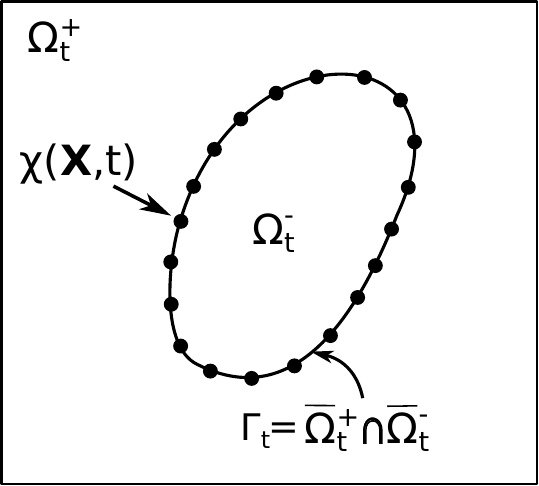}
	\caption[]{
		Lagrangian and Eulerian coordinate systems for a thin massless interface undergoing infinitesimal displacements. 
		At time $t$, the physical position of a material point $\s \in \Gamma_0$ attached to the interface and parameterized 
		by $\s=\s(\iota,\nu)$ is $\X(\s, t) \in \Gamma_t$, in which $(\iota,\nu)$ are curvilinear coordinates. 
		 In practice, the $C^0$ interface discretization is the only geometrical information used in the numerical method.}
\label{fig:Lag_Eul_schematic} 
\end{figure}

This study concerns the flow of a viscous, incompressible Newtonian fluid interacting with an infinitesimally thin immersed boundary.
The computational domain is $\Omega$, and $\x \in \Omega$ indicates fixed physical coordinates.
We consider the case in which $\Omega$ is divided into an exterior region $\Omega^+ = \Omega^+_t$ and an interior region $\Omega^- = \Omega^-_t$, each parameterized by time $t$.
We require that $\Omega = \overline{\Omega^+_t} \cup \overline{\Omega^-_t}$, so that the interface is $\Gamma_t = \overline{\Omega^+_t} \cap \overline{\Omega^-_t}$.
Here, we consider only cases in which the interface is stationary or moves in a prescribed manner.
We describe the motion of the interface through reference coordinates $\s \in \Gamma_0$ attached to the interface at time $t=0$.
To simplify the numerical treatment of this problem, we use a penalty method that allows for small deviations between the prescribed position of the interface and its actual physical configuration.
The prescribed physical position of material point $\s$ at time $t$ is $\Y(\s,t)$, which has velocity $\vec{W}=\partial\Y(\s,t)/\partial t$, and the actual physical position of material point $\s$ at time $t$ is $\X(\s,t)$, which has velocity $\U(\s,t) = \partial\X(\s,t)/\partial t$.
Because there will generally be a discrepancy between the prescribed and actual configurations of the boundary, we use the actual configuration, which is determined by the computed motion map $\X : \Gamma_0 \times t \mapsto \Omega$, to determine the location where jump conditions are imposed and where the interface velocity is evaluated,
see Fig.~\ref{fig:Lag_Eul_schematic}.
The interface can also be described using curvilinear coordinates, say $(\iota,\nu)$, so that $\s=\s(\iota,\nu)$ and $\X = \X(\s(\iota,\nu),t)$.
The present numerical method does not use a global curvilinear coordinate system.
The finite element discretization of the interface geometry implicitly defines a local curvilinear coordinate system (i.e., in terms of the coordinates of the reference elements), but because this is a mesh-dependent coordinate system, we describe the interface geometry using reference coordinates attached to the interface.
This approach leads to mesh-independent parameters and can be readily implemented within an isoparametric finite element formulation.

We consider the case in which both the mass density $\rho$ and dynamic viscosity $\mu$ of the fluid are uniform across $\Omega$.
The effect of the immersed boundary appears as a singular force distribution in the equation of momentum balance. 
In this case, the governing equations are
   \begin{align}
       \label{eq:momentum} \rho \left( \D{\u}{t}(\x,t) + \u(\x,t) \cdot \grad \u(\x,t) \right) &= -\grad p(\x,t) + \mu \grad^2 \u(\x,t) + \f(\x,t) \\ 
       \label{eq:continuity}  \grad \cdot \u(\x,t) &= 0, \\
       \label{eq:fsiconstraint} \f(\x,t)             &= \int_{\Gamma_0} \F(\s,t) \, \delta(\x - \X(\s,t)) \, {\mathrm d} A, \\
        \label{eq:fsiVel} \D{\X}{t}(\s,t)      &= \int_\Omega \u(\x,t) \, \delta(\x - \X(\s,t)) \, \Dx,
    \end{align}
in which $\u(\x,t)$ is the velocity, $p(\x,t)$ is the pressure, $\delta(\x)=\prod_{i=1}^{d}\delta(x_{i})$ is the $d$-dimensional delta function, $\F(\s,t)$ is the interfacial force along the immersed boundary, and $\f(\x,t)$ is the corresponding singular Eulerian force density.
Notice that interaction between the interface and the fluid occurs along $\Gamma_t = \X(\Gamma_0,t)$.
Additionally in Eq.~(\ref{eq:fsiVel}), $\U(\s,t) = \partial\X(\s,t)/\partial t = \u(\X(\s,t),t)$ because the velocity is continuous across the boundary as a result of the no-slip and no-penetration conditions along the interface.

For a rigid interface or an interface with prescribed kinematics, the interfacial force density is a Lagrange multiplier for the imposed motion.
In this study, we consider a penalty formulation similar to that proposed by Goldstein et al.~\cite{goldstein1993modeling}. 
In this formulation, the rigidity constraint is inexactly imposed through an approximate Lagrange multiplier force,
\begin{equation}
\F(\s,t) = \kappa\left(\Y(\s,t)-\X(\s,t)\right) + \eta \left(\vec{W}(\s,t) - \vec{U}(\s,t)\right).
 \label{eq:penalty_force}
\end{equation} 
The first term on the right-hand side is a linear spring force, with $\kappa$ a spring stiffness penalty parameter, and the second term is a damping force, with $\eta$ a damping penalty parameter.
Both forces act to provide energetic penalization if the structure deviates from its prescribed position.
Note that as $\kappa \rightarrow \infty$, the formulation imposes a hard constraint on the deformation.
Numerical experiments using both the standard IB method and the immersed interface approach also demonstrate that $\eta > 0$ can help reduce numerical oscillations, particularly in problems with high pressure forces.
In the special case of a stationary interface, $\Y(\s,t) = \s$ and $\vec{W}(\s,t) \equiv \vec{0}$, and the penalty force becomes simply
\begin{equation}
\F(\s,t) = \kappa\left(\s-\X(\s,t)\right) - \eta \vec{U}(\s,t).
 \label{eq:penalty_force_stationary}
\end{equation} 
A separate line of research is focused on how to enforce these types of constraints exactly \cite{kallemov2016immersed, balboa2017hydrodynamics}, but doing so requires solving a coupled system of equations involving an exact Lagrange multiplier force along with the Eulerian velocity and pressure fields.
The application of these types of approaches to the current sharp interface scheme is not addressed here.

\subsection{Physical jump conditions}
\label{subsec:formulation-PJC}

Throughout this subsection, $\s$ and $\x$ are taken to be corresponding positions in the reference and current configuration at time $t$, so that $\x = \X(\s,t)$.
A jump in a scalar field $\psi$ at position $\x = \X(\s,t)$ along the interface is
\begin{equation}
	\llbracket \psi(\vec{x},t) \rrbracket = \lim_{\epsilon \to 0^+} \psi(\vec{x} + \epsilon \vec{n}(\x,t),t) - \lim_{\epsilon \to 0^-} \psi(\vec{x} - \epsilon \vec{n}(\x,t),t)
	= \psi^{+}(\vec{x},t) - \psi^{-}(\vec{x},t),
	\label{eq:jump_definition}
\end{equation}
in which $\llbracket \cdot \rrbracket$ indicates the jump value, $\vec{n}(\x,t)$ is the outward unit normal vector along the interface $\Gamma_t$ in the current configuration, and $\psi^{+}(\vec{x},t)$ and $\psi^{-}(\vec{x},t)$ are the limiting values as approaching the interface position $\vec{x}$ from the exterior region $\Omega^{+}_t$ and interior region $\Omega^{-}_t$, respectively.
This definition can be extended for jumps in vectors in a component-wise manner.

For an incompressible Newtonian fluid, the fluid stress is
\begin{equation}
    \label{eq:fluidstress} \rttensor{\vec{\sigma}}_\text{f}(\x,t) = -p(\x,t) \, \II + \mu \left( \grad\u(\x,t) + \grad\u^T(\x,t) \right).
\end{equation} 
Because the boundary is infinitesimally thin, the interfacial force is balanced by a discontinuity in the fluid traction $\vec{\tau}_\text{f} = \rttensor{\vec{\sigma}}_\text{f}(\x,t) \cdot \vec{n}(\x,t)$ along the interface,
\begin{equation}
    \label{eq:basicjump}
	-\jmath^{-1}(\s,t) \, \F(\s,t) = \llbracket \vec{\tau}_\text{f}(\x,t) \rrbracket = \vec{\tau}_\text{f}^{+}(\x,t)-\vec{\tau}_\text{f}^{-}(\x,t),
\end{equation}
in which $\jmath$ is the Jacobian determinant that converts the surface force density from force per unit area in the current
 configuration to force per unit area in the reference configuration.
(The calculation of $\jmath$ is detailed in Sec.~\ref{subssubec:3.3.2}.)
 
\begin{figure}[t!!]
		\centering
			\includegraphics[width=0.45\textwidth]{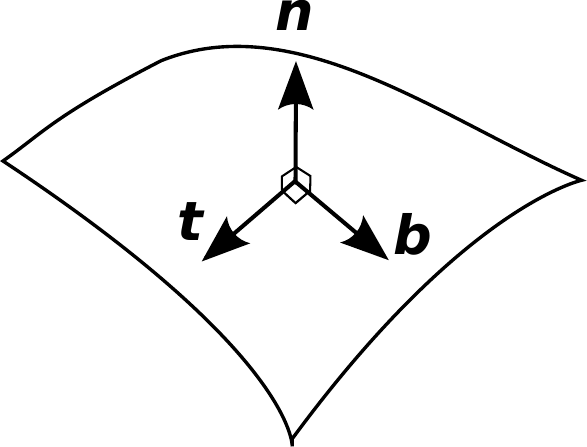}
		
		\caption{
		Unit normal ($\n$) and tangent ($\t$ and $\b$) vectors on the interface.
		}	
		\label{fig:schematic-normals} 
\end{figure}

Both velocity and pressure distributions are affected by the forces imposed along the immersed boundary. 
Away from the interface, these quantities can be assumed to be smooth, with the only possible discontinuities occurring at the immersed boundary.
The no-slip and no-penetration conditions require the velocity at the immersed interface to be continuous,
\begin{equation}
\llbracket \vec{u}(\x,t) \rrbracket=0.
\end{equation} 
Taking the dot-product of both sides in Eq.~(\ref{eq:fluidstress}) by $\vec{n}(x,t)$, and using Eq.~(\ref{eq:basicjump}), it is straightforward to show \cite{leveque1997immersed, lee2003immersed}
\begin{equation}
\label{eq:pj_codim1} \llbracket p(\vec{x},t) \rrbracket =  \jmath^{-1}(\s,t) \, \F(\s,t) \cdot \vec{n}(\vec{x},t). 
\end{equation}
Jump conditions for the shear stress can be derived by taking the dot-product of Eq.~(\ref{eq:fluidstress}) with tangential directions $\vec{b}(\vec{x},t)$ and $\vec{t}(\vec{x},t)$, see Fig.~\ref{fig:schematic-normals}.
If $\vec{n}=(n^x,n^y,n^z)$, $\vec{b}=(b^x,b^y,b^z)$, and $\vec{t}=(t^x,t^y,t^z)$ are the components of the local vectors in a local Cartesian coordinate system aligned with the interface, then
\begin{equation}
		\left \llbracket \begin{array}{c}
			\mu\frac{\partial\vec{u}}{\partial x}(\x,t)\\
			\mu\frac{\partial\vec{u}}{\partial y}(\x,t)\\
			\mu\frac{\partial\vec{u}}{\partial z}(\x,t)
		\end{array} \right \rrbracket =
				\left(\begin{array}{ccc}
			{n}^x 	& {b}^x 	& {t}^x \\
			{n}^y   & {b}^y		& {t}^y \\
			{n}^z 	& {b}^z 	& {t}^z	 
		\end{array} \right)
		\left \llbracket \begin{array}{c}
			\mu\frac{\partial\vec{u}}{\partial \vec{n}}(\x,t)\\
			\mu\frac{\partial\vec{u}}{\partial \vec{b}}(\x,t)\\
			\mu\frac{\partial\vec{u}}{\partial \vec{t}}(\x,t)
		\end{array} \right \rrbracket.
		\label{eq:CartTransformSystem}
\end{equation}
After some manipulation, the jump conditions become
\begin{align}
	\label{eq:du_dx_jump_codim1}     \left\llbracket \mu\frac{\partial\vec{u}}{\partial x}(\x,t) \right\rrbracket &= -(\rttensor{\vec{I}} - \vec{n}(\x,t)\otimes\vec{n}(\x,t)) \, \jmath^{-1}(\s,t) \, \F(\s,t) \, {n}^x, \\
	\label{eq:du_dy_jump_codim1}    \left\llbracket \mu\frac{\partial\vec{u}}{\partial y}(\x,t) \right\rrbracket &= -(\rttensor{\vec{I}} - \vec{n}(\x,t)\otimes\vec{n}(\x,t)) \, \jmath^{-1}(\s,t) \, \F(\s,t) \, {n}^y, \\
	\label{eq:du_dz_jump_codim1}   \left\llbracket \mu\frac{\partial\vec{u}(\x,t)}{\partial z} \right\rrbracket &= -(\rttensor{\vec{I}} - \vec{n}(\x,t)\otimes\vec{n}(\x,t)) \, \jmath^{-1}(\s,t) \, \F(\s,t) \, {n}^z,
\end{align}
which are convenient to use in the numerical implementation.
Together, these jump conditions for the pressure and velocity gradient state that the discontinuity in the fluid stress across the interface balances the force concentrated along the interface.

We remark that although there are generally discontinuities in the velocity gradient at the interface, $\grad \cdot \u(\x,t) \equiv 0$ holds throughout $\Omega$.
This can be seen by noting that the incompressibility condition holds for both the interior and exterior fluids, i.e., $\grad \cdot \u^{\pm}(\x,t) = 0.$
Applying the jump operator to both sides of this relation yields $\llbracket \grad \cdot \u(\x,t) \rrbracket = 0$.
Similarly, because the trajectories of material points cannot cross the interface, there also are no jump conditions associated with the convective derivative $\mathrm{D}\u/\mathrm{D}t(\x,t)$.

\section{Numerical Methods}
\label{sec:methodology}

This section describes the numerical approach, including the discretizations of the Eulerian and Lagrangian fields.
For simplicity, the numerical scheme is explained in two spatial dimensions.
The extension of the method to three spatial dimensions is straightforward.

\subsection{Eulerian finite difference approximation}
\label{subsec:methodology-FD}
A staggered-grid (MAC) discretization is used for the incompressible Navier-Stokes equations, which approximates the pressure at cell centers and the velocity and forcing terms at the edges (in two dimensions) or faces (in three dimensions) of the grid cells.
Standard second-order accurate centered approximation schemes are used for the divergence, gradient, and Laplace operators. 
The discrete divergence of the velocity $\vec{D} \cdot \u$ is evaluated at the cell centers, whereas the discrete pressure gradient $\vec{G} p$ and the components of the discrete Laplacian of the velocity $L \vec{u}$ are evaluated at the cell edges (or, in three dimensions, faces).
A staggered-grid version \cite{griffith2012volume,griffith2009} of the xsPPM7 variant \cite{rider2007accurate} of the piecewise parabolic method (PPM) \cite{colella1984piecewise} is used to approximate the nonlinear advection terms.
Previously described methods for physical boundary conditions \cite{griffith2009,griffith2012immersed} are used along the boundaries of the computational domain $\Omega$. 
In some tests, we use a locally refined Eulerian discretization approach described by Griffith \cite{griffith2012immersed} that employs Cartesian grid adaptive mesh refinement (AMR).

To determine the effects of the known jump conditions, we classify all Eulerian grid locations as \textit{regular} or \textit{irregular}.
Regular points are those for which none of the associated finite difference stencils cross the immersed boundary.
Because the solution is continuous away from the interface, the discretization at these points does not need to be modified.
For irregular points, however, at least one of the finite difference stencils associated with that point will cross the immersed boundary, and because the solution is discontinuous across the interface, standard discretizations that do not account for these jumps lead to large errors that do not decrease under grid refinement.
This motivates the application of the jump conditions through correction terms in the discretization.

Using a Taylor series expansion \cite{li2001immersed,xu2006systematic}, it can be shown that if the interface cuts between two Cartesian grid points at location $\x_{\circ}$, such that $x_{i,j,k}\leq x_{\circ} < x_{i+1,j,k}$, with $\x_{i,j,k}\in \Omega^{-}$ and $\x_{i+1,j,k}\in \Omega^{+}$, then for a piecewise differentiable quantity $\psi$ we have, for instance,
\begin{align}
\frac{\partial\psi}{\partial x}(\x_{i+\half,j,k}) &=\frac{\psi_{i+1,j,k}-\psi_{i,j,k}}{\Delta x}+\frac{\mathrm{sgn}\{n^x\}}{\Delta x}\displaystyle \sum_{m=0}^{2}\frac{(d^{+})^{m}}{m!}{\left\llbracket \frac{\partial^{m}\psi}{\partial x^{m}}\right\rrbracket}_{x_{\circ}} + O({\Delta x}^{2}), \\
\frac{\partial^{2}\psi}{\partial x^{2}}(\x_{i,j,k}) &=\frac{\psi_{i+1,j,k}-2\psi_{i,j,k}+\psi_{i-1,j,k}}{\Delta x^{2}}+\frac{\mathrm{sgn}\{n^x\}}{{\Delta x}^{2}}{\displaystyle \sum_{m=0}^{3}\frac{(d^{+})^{m}}{m!}{\left\llbracket \frac{\partial^{m}\psi}{\partial x^{m}}\right\rrbracket}_{x_{\circ}} + O({\Delta x}^{2})},
\end{align}
in which $\Delta x$ is the grid spacing in the $x$ direction, $\psi_{i,j,k}=\psi(\x_{i,j,k})$, $d^{+}= x_{i+1,j,k}-x_{\circ} > 0$, and $n^x$ is the $x$-component of the normal vector $\vec{n}=(n^x,n^y,n^z)$ at the intersection point $\x_\circ$.

\begin{figure}[t]
		\centering
			\includegraphics[width=0.9\textwidth]{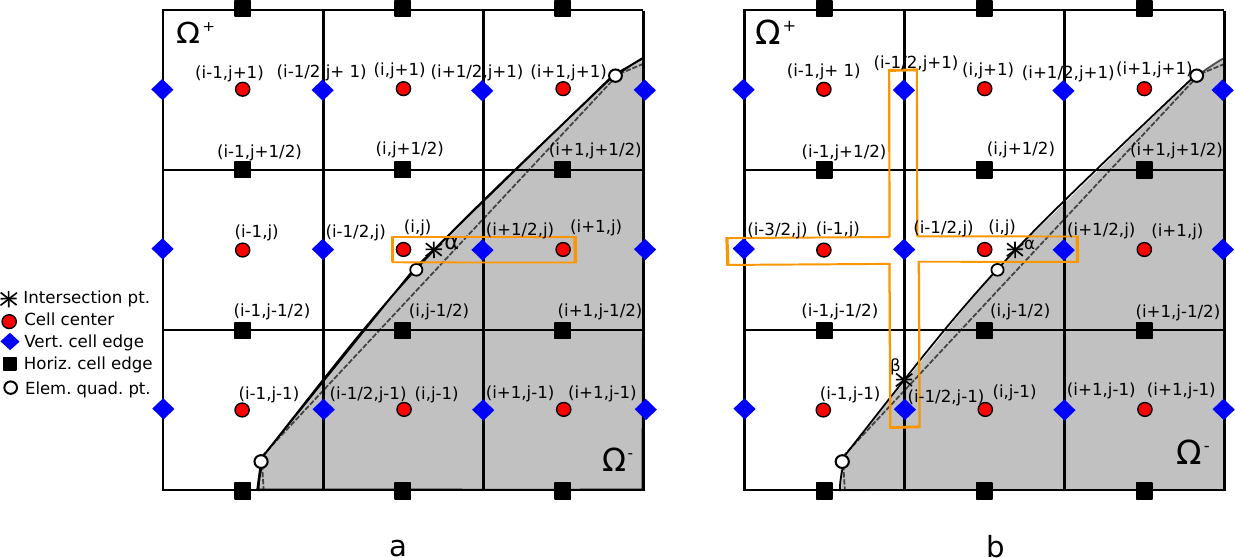}
		
		\caption{Two-dimensional illustration of irregular stencils in the Navier-Stokes equations for the $x$-component of (a) the pressure gradient $(G^{x}p)_{i+\half,j}$ and  
		 (b) Laplacian of the $x$-component of the velocity $(L^{x} u)_{i-\frac{1}{2},j}$ }	
		\label{fig:Jump-grid-px-d2u}		
\end{figure}

We now demonstrate how these corrections can be applied to irregular stencils of the pressure and the viscous terms in the $x$-component of the Navier-Stokes equations for the two-dimensional staggered mesh arrangements shown in Fig.~\ref{fig:Jump-grid-px-d2u}.
The extension to three spatial dimensions is straightforward.
We consider only the lowest-order jump conditions, accounting for discontinuities in the pressure and the derivative of the velocity.
Let $\mathbf{G}p=(G^{x}p,G^{y}p)$ be the discrete pressure gradient.
The modified discretization including the correction term for the $x$-component of this vector, e.g.~as shown in Fig.~\ref{fig:Jump-grid-px-d2u}(a), is
\begin{equation}
(G^{x}p)_{i+\half,j} =\frac{p_{i+1,j}-p_{i,j}}{\Delta x} -\frac{\llbracket p(\vec{x},t) \rrbracket}{\Delta x}. 
\label{eq:pressure-stencil-modifed}
\end{equation}
Similarly we represent the finite difference approximation of the vector Laplacian of  $\vec{u}=(u,v)$ by $\vec{Lu}=(L^{x}u,L^{y}v)$, in which $L^{x}$ and $L^{y}$ are in fact the same differencing operator applied to different components of the velocity vector.
The modified discretization of the $L^{x}u$, e.g.~as shown in Fig.~\ref{fig:Jump-grid-px-d2u}(b), is
\begin{equation}
\begin{split}
(L^{x} u)_{i-\frac{1}{2},j} =& \frac{u_{i+\frac{1}{2},j}-2u_{i-\frac{1}{2},j}+ u_{i-\frac{3}{2},j}}{\Delta x^{2}} + \frac{u_{i-\frac{1}{2},j+1}-2u_{i-\frac{1}{2},j}+ u_{i-\frac{1}{2},j-1}}{\Delta y^{2}}\\
& \mbox{} + \frac{h_{\alpha}}{\Delta x^{2}}\left\llbracket \frac{\partial\vec{u}}{\partial x} \right\rrbracket_{\alpha} - \frac{h_{\beta}}{\Delta y^{2}}\left\llbracket \frac{\partial\vec{u}}{\partial y} \right\rrbracket_{\beta}, 
\end{split}
\end{equation}
in which the positive values of $h_{\alpha}$ and $h_{\beta}$ are defined as $h_{\alpha} = x_{i+\half,j} - x_{\alpha}$ and $h_{\beta}=y_{\beta}-y_{i-\half,j-1}$.

For the remainder of the paper, we assume an isotropic Cartesian grid, so that $\Delta x = \Delta y = \Delta z = h$.

\subsection{Lagrangian discretization}
\label{subsec:methodology-FE}

\begin{figure}[t!!]
		\centering
			\includegraphics[width=0.45\textwidth]{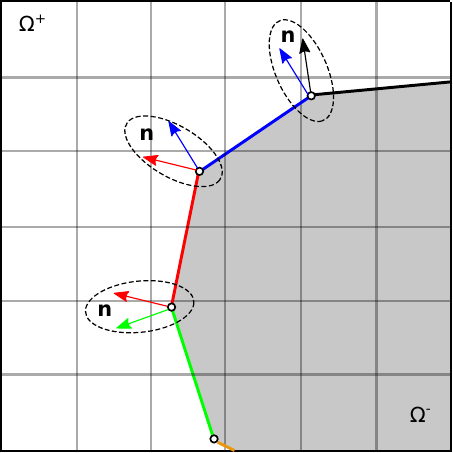}
		
		\caption{Two-dimensional schematic of a $C^{0}$-continuous immersed boundary. Normal vectors are discontinuous at the junction between neighboring elements.}	
		\label{fig:normals-elements} 
\end{figure}

To define a finite element approximation, consider $\cT_h$ as a triangulation of $\Gamma_0$, the reference configuration of the interface, with elements $U^e$ such that $\cT_h = \cup_{e} U^e$, in which $e$ indexes the mesh elements.
$\left\{\s_l\right\}_{l=1}^{M}$ are the nodes of the mesh and $\left\{\phi_l(\s)\right\}_{l=1}^{M}$ 
  are the corresponding nodal (Lagrangian) basis functions.
Using the finite element basis functions and the time-dependent physical positions of the nodes of the Lagrangian mesh $\{\X_l(t)\}_{l=1}^{M}$, the approximation to the interfacial deformation $\X_{h}(\s,t)$ is defined by
\begin{equation}
\X_{h}(\s,t)={\displaystyle \sum_{l=1}^{M}\X_{l}(t)\phi_{l}(\s)}.  \label{eq:discrete_deformation}
\end{equation}
From this representation of the deformation, it is straightforward to evaluate approximations to geometrical quantities such as the surface normal or surface Jacobian determinant by differentiating Eq.~\eqref{eq:discrete_deformation}.
Because we use $C^0$ basis functions in the present work, quantities that are obtained in terms of $\D{\X_h}{\s}$ are discontinuous in both the reference and current configurations, as shown in Fig.~\ref{fig:normals-elements}.

\subsection{Projected Jump conditions}
\label{subssubec:3.3.2}  
 
Because the normal direction to the interface is generally discontinuous on a $C^0$ mesh at the mesh nodes, 
pointwise jump conditions determined from the mesh geometry and the surface Jacobian $\jmath$ are generally discontinuous, see Fig.~\ref{fig:normals-elements}.
Instead of using these pointwise values directly, we project the jump conditions onto the subspace $S_h = \mathrm{span}\{\phi_l(\s)\}_{l=1}^{M}$.
Given a function $\psi \in L^{2}(\Gamma_0)$, its $L^{2}$ projection $P_{h} \psi$ onto the subspace $S_h$ is defined by requiring $P_{h} \psi$ to satisfy
\begin{equation}
   \int_{\Gamma_0} \big(\psi(\s) - P_{h}\psi(\s)\big) \, \phi_l(\s) \,  {\mathrm d} A = 0, \quad \forall l=1,\ldots,M.
   \label{eq:L2_proj_def}
\end{equation}
The $L^{2}$ projection of a vector-valued quantity is determined component-wise.
Because the $L^{2}$ projection is defined via integration, the function $\psi$ does not need to be continuous or even to have well-defined nodal values.
By construction, however, $P_h \psi$ will inherit any smoothness provided by the subspace $S_h$.
In particular, for $C^0$ Lagrangian basis functions, $P_h \psi$ will be at least continuous.
In our numerical scheme, we set $C^\text{n}_h(\s,t)$ to be the $L^2$ projection of the normal component of the surface force per unit current area, $\jmath^{-1} \F_h(\s,t)$, onto $S_h$, and we set $\vec{C}^\text{t}_h(\s,t)$ to be the projection of the tangential component of the force.
Solving for the projected jump conditions require solving linear systems of equations involving the mass matrix $M$ with components $M_{kl} = \int \phi_k \phi_l \, \DA$.
In practice, Eq.~\eqref{eq:L2_proj_def} is approximated using a numerical quadrature rule. 
To simplify notation, the subscript ``h" is mostly dropped in the remainder of the paper when showing numerical
approximations to the Lagrangian variables

%

In implementing the jump conditions, it is convenient to evaluate them in the Cartesian directions.
To do so, it is necessary to determine the normal and tangents to the interface.
We let vectors $\mathbf{e}_{\iota}$ and $\mathbf{e}_{\nu}$ be tangents to the local element coordinates $\iota$ and $\nu$, so that
\begin{equation}
\mathbf{e}_{\iota}=\left[\begin{array}{c}
\frac{\partial x}{\partial\iota}\\
\frac{\partial y}{\partial\iota}\\
\frac{\partial z}{\partial\iota}
\end{array}\right]
\text{ and }
\mathbf{e}_{\nu}=\left[\begin{array}{c}
\frac{\partial x}{\partial\nu}\\
\frac{\partial y}{\partial\nu}\\
\frac{\partial z}{\partial\nu}
\end{array}\right].
\end{equation}
The derivatives of the global coordinates with respect to local coordinates can be determined using the basis functions via
\begin{equation}
\frac{\partial \X}{\partial \iota}(\s,t)={\displaystyle \sum_{l=1}^{M}\X_{l}(t)\frac{\partial \phi_{l}(\s)}{\partial \iota}} \text{ and } 
\frac{\partial \X}{\partial \nu}(\s,t)={\displaystyle \sum_{l=1}^{M}\X_{l}(t)\frac{\partial \phi_{l}(\s)}{\partial \nu}}.  
\end{equation}
The \textit{area-weighted} normal vector $\n$ in the global coordinate is then obtained by evaluating $\n = \mathbf{e}_{\iota} \times \mathbf{e}_{\nu}$.
Similarly for the normal in the local coordinate, $\N$, we use the Lagrangian material coordinates,
\begin{equation}
\mathbf{E}_{\iota}=\left[\begin{array}{c}
\frac{\partial X}{\partial\iota}\\
\frac{\partial Y}{\partial\iota}\\
\frac{\partial Z}{\partial\iota}
\end{array}\right]
\text{ and }
\mathbf{E}_{\nu}=\left[\begin{array}{c}
\frac{\partial X}{\partial\nu}\\
\frac{\partial Y}{\partial\nu}\\
\frac{\partial Z}{\partial\nu}
\end{array}\right],
\end{equation}
yielding $\N = \mathbf{E}_{\iota} \times \mathbf{E}_{\nu}$.
The surface Jacobian, which is the area ratio between the two coordinates systems, is obtained via
\begin{equation}
\jmath = \frac{\mathrm d a}{\mathrm d A} = \frac{\|\n\|}{\|\N\|}. 
\end{equation}
In the present study, $\jmath \approx 1$ because we allow for only infinitesimal deformations of the interface.

\subsection{Intersection algorithm}
\label{subsec:Intersection algorithm}

Intersections between the Cartesian finite difference stencils and the interface are determined by a ray-casting algorithm.
In two dimensions, this is done by finding the intersections of Cartesian-oriented lines with first-order parametric elements, which amounts to solving a single linear equation.
In three-dimensional cases, we use the M{\"o}ller-Trumbore ray-triangle intersection algorithm, which is a fast method to calculate intersections with surface triangles \cite{moller2005fast}, 
and a similar efficient approach for convex planar quadrilateral is used to find intersections of stencils with rectangular surface elements \cite{lagae2005efficient}.
Although there should be at most one intersection between a given finite difference stencil and the interface, finite precision effects may yield spurious multiple intersections, or may completely miss intersections between the finite difference stencil and the interface.
To avoid these situations, before determining these intersections, we perturb the positions of the control points of the interface mesh away from the cell centers, nodes, edges, and faces of the background Cartesian grid by a factor proportional to $\sqrt{\epsilon_\text{mach}}$, 
in which $\epsilon_\text{mach}$ is machine precision. This is to ensure that we do not have to use specialized methods for treating apparent multiple intersections.
This perturbation in the position of the interface is much smaller than the accuracy of the overall method, and it avoids the need to implement specialized discrete geometry algorithms that explicitly treat such spurious intersections.

\subsection{Velocity interpolation}
\label{subsec:vel-interpolation}

The goal here is to determine a velocity-restriction operator $\cJ = \cJ[\X,\F]$ to determine the Lagrangian mesh velocity $\U$ from the Eulerian grid velocity $\u$.
The velocity interpolation used herein consists of two steps.
Ultimately, we project the interpolated velocity onto the Lagrangian basis functions to determine the motion of the interface, and so the first step of the interpolation procedure is to evaluate the Lagrangian velocity $\V$ at the quadrature points along the interface that are used to compute the $L^2$ projection of the interpolated velocity.
The second step projects the velocity field $\V$ onto the space $S_h$ spanned by the nodal finite element basis functions to determine the mesh motion $\U$.

\begin{figure}[t]
		\centering
			\includegraphics[width=0.45\textwidth]{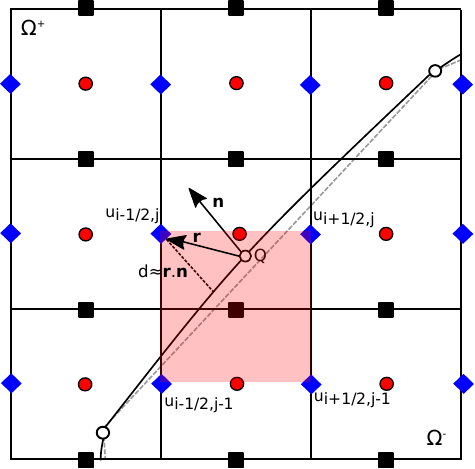}
		
		\caption{Two-dimensional schematic of the interpolation scheme for the $x$-component of the velocity at the interfacial location $\x = \X(\s,t)$.
			The cell used for calculating the modified bilinear interpolation stencil is highlighted in red.} 
		   \label{fig:interp_schematic}	 
\end{figure}

		

Because the projection step is standard, we focus here on evaluating $\V$ at a generic reference position $\s$ on the interface, with corresponding physical position $\x = \X(\s,t)$.
Our approach is to use modified bilinear (or, in three spatial dimensions, trilinear) interpolation that accounts for the known discontinuities in the velocity gradient at the interface.
Considering the two-dimensional schematic detailed in Fig.~\ref{fig:interp_schematic}, the general formula for the $x$-component of the velocity, $V^x$, at the interfacial location is
\begin{align}
  \label{eq:V_interp} V^x(\s,t) &=  (1-\zeta)(1-\lambda)u_{i-\frac{1}{2},j-1} + C_{i-\frac{1}{2},j-1} + \zeta(1-\lambda)u_{i+\frac{1}{2},j-1} + C_{i+\frac{1}{2},j-1}\\ \nonumber
          & \ \ \ \mbox{} + \zeta\lambda u_{i+\frac{1}{2},j} + C_{i+\frac{1}{2},j} + (1-\zeta)\lambda u_{i-\frac{1}{2},j} + C_{i-\frac{1}{2},j},
\end{align}
in which $\lambda=\frac{x-x_{i-\frac{1}{2},j-1}}{h}$ and $\zeta=\frac{y-y_{i-\frac{1}{2},j-1}}{h}$ and $\Delta x = \Delta y = h$ is the (isotropic) Cartesian grid spacing.
The corrections $C_{i-\frac{1}{2},j-1}$, $C_{i+\frac{1}{2},j-1}$, $C_{i+\frac{1}{2},j}$ and $C_{i-\frac{1}{2},j}$ have the following forms:
\begin{align}
   C_{i-\frac{1}{2},j-1} &=
      \begin{cases}
       h(1-\zeta)(1-\lambda) \big(\zeta n^x \llbracket\frac{\partial u}{\partial x}\rrbracket + \lambda n^y \llbracket\frac{\partial u}{\partial y}\rrbracket \big) & \x_{i-\frac{1}{2},j-1} \in \Omega^{+}, \\
        0 & \x_{i-\frac{1}{2},j-1} \in \Omega^{-},
     \end{cases} \\
   C_{i+\frac{1}{2},j-1} &=
      \begin{cases}
       -h\zeta(1-\lambda) \big((1-\zeta) n^x \llbracket\frac{\partial u}{\partial x}\rrbracket - \lambda n^y \llbracket\frac{\partial u}{\partial y}\rrbracket \big) & \x_{i+\frac{1}{2},j-1} \in \Omega^{+}, \\
        0 & \x_{i+\frac{1}{2},j-1} \in \Omega^{-},
     \end{cases} \\
   C_{i+\frac{1}{2},j} &=
      \begin{cases}
       -h\zeta\lambda \big((1-\zeta) n^x \llbracket\frac{\partial u}{\partial x}\rrbracket + (1-\lambda) n^y \llbracket\frac{\partial u}{\partial y}\rrbracket \big) & \x_{i+\frac{1}{2},j} \in \Omega^{+}, \\
        0 & \x_{i+\frac{1}{2},j} \in \Omega^{-},
     \end{cases} \\
   C_{i-\frac{1}{2},j} &=
      \begin{cases}
       h(1-\zeta)\lambda \big(\zeta n^x \llbracket\frac{\partial u}{\partial x}\rrbracket - (1-\lambda) n^y \llbracket\frac{\partial u}{\partial y}\rrbracket \big) & \x_{i-\frac{1}{2},j} \in \Omega^{+}, \\
        0 & \x_{i-\frac{1}{2},j} \in \Omega^{-}.
     \end{cases}
\end{align}
Notice that the jump conditions for the velocity gradient can be determined from the interfacial force $\F$.
This modified bilinear interpolation is similar to the approach of Tan et al.~\cite{tan2008immersed,tan2009immersed}.
An important difference is in the value used for the distance from the Cartesian grid locations to the interface in the Taylor series expansion of the correction terms.
The approach by Tan et al.~uses the distance between the grid points and the interfacial point ($|\vec{r}|$ in Fig.~\ref{fig:interp_schematic}).
Here, this is modified by projecting the location vector $\vec{r}$ to the normal vector $\n$ and using $|\vec{r}\cdot\n|$ instead.
Taking into account for the normal vector naturally provides a more accurate approximation to the actual distance between the grid point and the interface.

$\V(\s,t)$ can be evaluated at any $\s \in \Gamma_0$, but generally, $\V(\s,t) \not\in S_h = \sspan\left\{\phi_l(\s)\right\}$.
Therefore, the second step of the interpolation is to set $\U = P_h \V$, component-by-component.
This completes the construction of $\U = \cJ[\X,\F] \, \u$.

\subsection{Calculating the fluid traction}
\label{subsec:calc-traction}

Tests of the present method will evaluate components of the fluid traction along the interface.
For simplicity, we focus on evaluating exterior tractions, but the evaluation of interior tractions is similar.
To obtain an approximation to the exterior pressure at position $\x$ along the interface, we use
\begin{equation}
 p^{+}_h(\vec{x},t) =  \llbracket p(\vec{x},t) \rrbracket + \cI[p](\vec{x}^{-},t)
 \end{equation}
in which $p^{-}=\cI[p](\vec{x}^{-},t)$ is the interior pressure interpolated to a position $\x^{-}$ away from the interface in the opposite direction of the normal vector $\n$ and at a distance equal to $1.2$ times the diagonal size of one grid cell. 
Here, $\cI$ is the unmodified bilinear (or trilinear) interpolation operator.
By construction, $\x^-$ is sufficiently far from the interface to ensure that uncorrected interpolation may be used without degradation in accuracy.
To compute the exterior wall shear stress, the one-sided normal derivative of the velocity is calculated using the interfacial velocity reconstruction described in Sec.~\ref{subsec:vel-interpolation} along with the velocity value at a neighboring location in the direction of the normal vector $\vec{x}^{+}$.
As with the pressure, unmodified bilinear (or trilinear) interpolation is used to obtain the velocity away from the interface.
A one-sided finite difference formula is used to obtain an approximation to the normal derivative as
\begin{equation}
 \left(\frac{\partial \u}{\partial \vec{n}}\right)^{+}_h(\vec{x},t) = \frac{\cI[\u](\vec{x}^{+},t)-\u(\vec{x},t)}{\hat{h}}.
\end{equation}
where the distance $\hat{h}$ is chosen to be slightly larger than the diagonal size of the Cartesian mesh ($1.05$ times the diagonal size), so that regular bilinear (or trilinear) interpolation can be used to evaluate $\cI[\u](\vec{x}^{+},t)$.
A second-order formula using a three-point stencil can also used that requires interpolating an additional point in the normal direction. 
However, preliminary numerical experiments (data not shown) suggest the computation using only two points suppresses oscillations that appear to be associated with the higher-order stencil.
Moreover, as shown in our tests, this simple scheme is adequate to achieve the designed first-order accuracy of the wall shear stress.

Note that as with velocity interpolation, the pressure and wall shear stress can be evaluated at arbitrary locations along the interface.
As before, to obtain nodal values of these quantities along the interface, we use the $L^2$ projection of the interfacial values.

\subsection{Time-integration scheme}
\label{subsec:time-integration}

Starting from the values of $\X^{n}$ and $\vec{u}^{n}$ at time $t^n$ and $p^{n-\frac{1}{2}}$ at time $t^{n-\frac{1}{2}}$, we must compute $\X^{n+1}$, $\u^{n+1}$, and $p^{n+\frac{1}{2}}$. 
To do this, we first determine a preliminary approximation to the structure location at time $t^{n+\frac{1}{2}}$ by
\begin{equation}
\frac{\hat{\X}^{n+1}-\X^{n}}{\Delta t}=\U^n=\cJ^{n} \, \vec{u}^{n}
\end{equation}
with $\cJ^{n}=\cJ[\X^{n},\F^{n}]$ is the discrete velocity restriction operator described in Sec.~\ref{subsec:vel-interpolation}.
We also define an approximation to $\X$ at time $t^{n+\frac{1}{2}}$ by
\begin{equation}
\X^{n+\frac{1}{2}}=\frac{\hat{\X}^{n+1}+\X^{n}}{2}.
\end{equation}
Next, we solve for $\X^{n+1}$, $\u^{n+1}$, and $p^{n+\frac{1}{2}}$ via
\begin{align}
\label{eq:time-integ-momentum}\rho\left(\frac{\u^{n+1}-\u^{n}}{\Delta t}+\vec{A}^{n+\frac{1}{2}}\right) &=-\vec{G}p^{n+\frac{1}{2}}+\mu\vec{L}\left(\frac{\u^{n+1} + \u^{n}}{2}\right)+ \vec{f}^{n+\frac{1}{2}}, \\
\vec{D} \cdot \vec{u}^{n+1} &=0, \\
\frac{\X^{n+1}-\X^{n}}{\Delta t} &= \U^{n+\half} = \cJ^{n+\half} \left(\frac{\u^{n+1}+\u^{n}}{2}\right),
\end{align}
in which $\vec{G}$, $\vec{D} \cdot \mbox{}$, and $\vec{L}$ are, respectively, the discrete gradient, divergence, and Laplace operators, $\vec{A}^{n+\frac{1}{2}}=\frac{3}{2}\vec{A}^{n}-\frac{1}{2}\vec{A}^{n-1}$ is obtained from a high-order upwind spatial discretization of the nonlinear convective term $\u \cdot \grad \u$ \cite{griffith2009}, and $\vec{f}^{n+\frac{1}{2}}$ is discrete Eulerian body force on the Cartesian grid that corresponds to the sum of all the correction terms computed using the Lagrangian force $\vec{F} = \vec{F}[\X^{n+\frac{1}{2}}, \U^{n}, t^{n+\half}]$.
The semi-implicit Crank-Nicolson form in Eq.~(\ref{eq:time-integ-momentum}) requires only linear solvers for the time-dependent incompressible Stokes equations.
This system of equations is solved via the flexible GMRES (FGMRES) algorithm with an approximate factorization preconditioner based on the projection method that uses inexact subdomain solvers \cite{griffith2009}.
In the initial time step, a two-step predictor-corrector method is used to determine the velocity, deformation, and pressure, see Griffith and Luo \cite{BEGriffith17-ibfe} for further details.
Note that similar to the standard IB method, we can relate $\f$ and $\F$ by a force-spreading operator $\cS = \cS[\X]$ such that $\vec{f}=\cS[\X] \, \vec{F}$.
Unlike the conventional IB method, however, the present force-spreading and velocity-restriction operators are not adjoints.

\subsection{Standard IB formulation}
\label{subssubec:3.3.1}

For tests in this study that use the conventional IB method, we approximate the force using the nodal forces $\{\F_l(t)\}$ and the shape functions, so that
\begin{equation}
\mathbf{F}_{h}(\s,t)={\displaystyle \sum_{l=1}^{M}\mathbf{F}_{l}(t)\phi_{l}(\s)}.
\end{equation}
Using Eq.~(\ref{eq:penalty_force}), the discretized Lagrangian penalty force $\mathbf{F}_{h}(\s,t)$ can then be directly evaluated from $\X_h$ and the interface motion $\U_h$.
These forces are spread to the background grid using a discretized integral transform with a regularized delta function kernel.
The adjoint of this operator is used to determine the velocity of the mesh from the Cartesian grid velocity field.
In the present study, we use a piecewise linear regularized delta function.
Details are provided by Griffith and Luo \cite{BEGriffith17-ibfe}.

\subsection{Software implementation}
\label{sec:implementation}

The present approach is implemented in the open-source IBAMR software \cite{ibamr}, which is a C++ framework for FSI modeling using the IB formulation.
IBAMR provides support for large-scale simulations through the use of distributed-memory parallelism and adaptive mesh refinement (AMR).
IBAMR relies upon other open-source software libraries, including SAMRAI \cite{samrai, hornung2002managing}, PETSc \cite{petsc-web-page,petsc-user-ref, petsc-efficient}, \textit{hypre} \cite{hypre,falgout2002hypre}, and \texttt{libMesh} \cite{libmesh,Kirk:2006:LCL:1230680.1230693}, for various key functionalities.

\section{Numerical results}
\label{sec:results}

Verification examples in two and three spatial dimensions, including comparisons to benchmark computational and experimental studies, are used to investigate the accuracy of the proposed IIM.
We consider tests involving internal and external flow conditions, treating both stationary boundaries and boundaries with prescribed kinematics.
All the interfaces in two~and three-dimensions are covered by linear ($P^1$) Lagrangian elements.
Wherever there is an exact solution to compare against, the analytic interface geometry is used to determine the exact solution.
Those computations use only the discrete interface geometry, which in our tests are described using either one-dimensional line elements (in two spatial dimensions) or two-dimensional triangular or quadrilateral surface elements (in three spatial dimensions).
Analytic information about the interface geometry is only used to determine the discrete representation of the reference configuration of the interface.
When assessing the accuracy of the computed fluid traction, we always consider the traction on the side of the interface with the nontrivial flow.

In all cases, the penalty forces that approximately impose the interface kinematics act to support physical forces exerted by the flow on the interface, and the physical forces are a property of the model.
Thus, under grid refinement, we expect that the penalty forces will converge to the physical loading forces.
Notice also that the penalty forces are proportional to the displacement $\Y(\s,t)-\X(\s,t)$ between the prescribed and computed interface positions.
If we wish to achieve $\|\Y(\s,t)-\X(\s,t)\| = O(h^2)$, it is necessary that the penalty parameter $\kappa$ also satisfies $\kappa = O(1/h^2)$, so that an applied penalty force of the form $\F = \kappa (\Y(\s,t)-\X(\s,t))$ satisfies $\|\F\| = O(1)$ under grid refinement.
Because we maintain $\dt = O(h)$ in our convergence tests, which keeps the advective Courant-Friedrichs-Lewy (CFL) number fixed under grid refinement, it is convenient to choose $\kappa = \kappa_0 / \dt^2$.
We choose $\kappa_0$ by, first, choosing the finest grid spacing to be considered in the convergence study, say $h_\text{min}$, along with the corresponding time step size $\dt_\text{min}$.
We then empirically determine approximately the largest value of $\kappa_0$ that allows the scheme to remain stable when $h = h_\text{min}$ and $\dt = \dt_\text{min}$.
We then use the prescribed relationship between $\kappa$ and $\dt$ to determine $\kappa$ for all coarser cases.
This allows us to ensure that the numerical parameters are stable for all grid spacings considered in each test while using scalings that ensure that the method achieves its designed order of accuracy.
All computations use a tight relative convergence threshold of $\epsilon_\text{rel}=\text{1e-10}$ for all iterative linear solvers.

When assessing the accuracy of discontinuous quantities, such as the Eulerian pressure $p$, small perturbations in the interface position can result in $O(1)$ differences in the computed value at spatial locations $\x$ that are close to the interface, even if the interface position is computed very accurately.
At points within a mesh width of the interface, the pointwise errors can differ by the value of the jump in the solution at the interface.
These geometrical errors can result from approximation error as well as the finite precision effects.
To assess the impact of these unavoidable geometrical effects on the accuracy of the computed solutions, we report pointwise errors in the pressure both over the full domain $\Omega$, and over the subregion $\Omega^*$ that excludes a band with a width of two grid cells around the interface.  
Pointwise errors in the quantity denoted $p^*$ are computed over $\Omega^*$.
The numerical method is also shown to capture discontinuous quantities such as the pressure and wall shear stress on the interface itself, which indicates that errors in the location of the interface that effect the pointwise values of the Eulerian fields do not spoil the ability of the method to impose jump conditions accurately on those fields.

\subsection{Pressure-driven flow}
\label{subsec:pressure-driven}

The first example considers steady-state pressure-driven flow inside a channel.
Results are presented for both two-dimensional plane-Poiseuille channel flow as well as three-dimensional Hagen-Poiseuille flow inside a circular pipe.
In each case, both grid-aligned and skew structure configurations are used. 
For these tests, the interface comes into contact with the exterior boundary of the computational domain, $\partial \Omega$, and errors at these junctions can dominate some pointwise errors.
Because we wish to examine the accuracy of the method for fully developed flows, and not to assess the treatment of numerical boundary conditions for the flow solver, we perform our error analyses over a slightly smaller region that excludes positions near $\partial\Omega$ that are within 10\% of the length of the computational domain.

\subsubsection{Plane-Poiseuille flow}
\label{subsubsec:plane-poiseuille}

We first consider a two-dimensional domain $\Omega = [0,L]^2$.
A horizontal channel of width $H$ extends across the middle of the domain.
The two channel ends at $x=0$ and $x=L$ are subject to constant pressure $p_0$ and $-p_0$ boundary conditions, establishing a constant pressure gradient of $2p_0/L$ across the channel.  
The prescribed-pressure boundary condition is imposed by means of applying combined normal-traction and zero-tangential-slip boundary conditions.
At the remaining parts of the left ($x=0$) and right ($x=L$) sides of $\Omega$ along with the entire bottom ($y=0$) and top ($y=L$) boundaries, solid wall (zero-velocity) boundary conditions are imposed.
The channel walls are modeled using one-dimensional linear ($P^1$) Lagrangian elements with an element size twice that of the Eulerian grid spacing, which results in a mesh ratio of $\Mfac = 2$.
Note that in this context and throughout the rest of this paper, $\Mfac$ is the ratio of Lagrangian element size to the Eulerian grid spacing.
The steady-state solution of pressure-driven flow for such a channel is described by the plane Poiseuille equation,
  \begin{align}
	\label{eq:PPF_u_2D} u(y)&=\frac{p_{0}H}{\mu L}\left(y-y_{0}\right)\left(1-\frac{y-y_{0}}{H}\right), \\
	\label{eq:PPF_v_2D} v&=0, \\
	\label{eq:PPF_p_2D} p(x)&=p_{0}-2p_{0}x/L, 
\end{align}
in which $y_0$ is the height of the lower wall of the channel.
Remaining simulation parameters include $\rho=1$, $\mu=0.01$, $L=5$, $H=1$, $y_0=2$ and $p_0=0.2$, resulting in a maximum velocity of $U_{\text{max}}=1$, an average velocity of $\bar{U}=2/3$, and the Reynolds number $\Re=\frac{\rho H \bar{U}}{\mu}\approx 66.66$.
The time step size is $\Delta t=0.1h$, yielding a maximum advective CFL number of approximately $0.1$--$0.2$.
The penalty parameters are $\kappa=10^{-3}/\Delta t^2$ and $\eta = 0$.

Along with the horizontal plane-Poiseuille flow, an inclined version of the channel is also studied in which the immersed boundary is not aligned with the Cartesian directions.
In this slightly more challenging scenario, the channel is at an angle $\theta=\pi/12$ with the horizontal direction.
At the inlet and outlet of the channel, the rotated exact solution is used to determine normal traction and tangential velocity boundary conditions, so that the flow conditions are the same as in the aligned configuration.
At the remaining parts of $\partial\Omega$, solid-wall boundary conditions are imposed.
Other than the rotation in the geometry, all other simulation parameters are the same as the horizontal case.

\begin{figure}[t!!]
		\centering
			\includegraphics[width=0.9\textwidth]{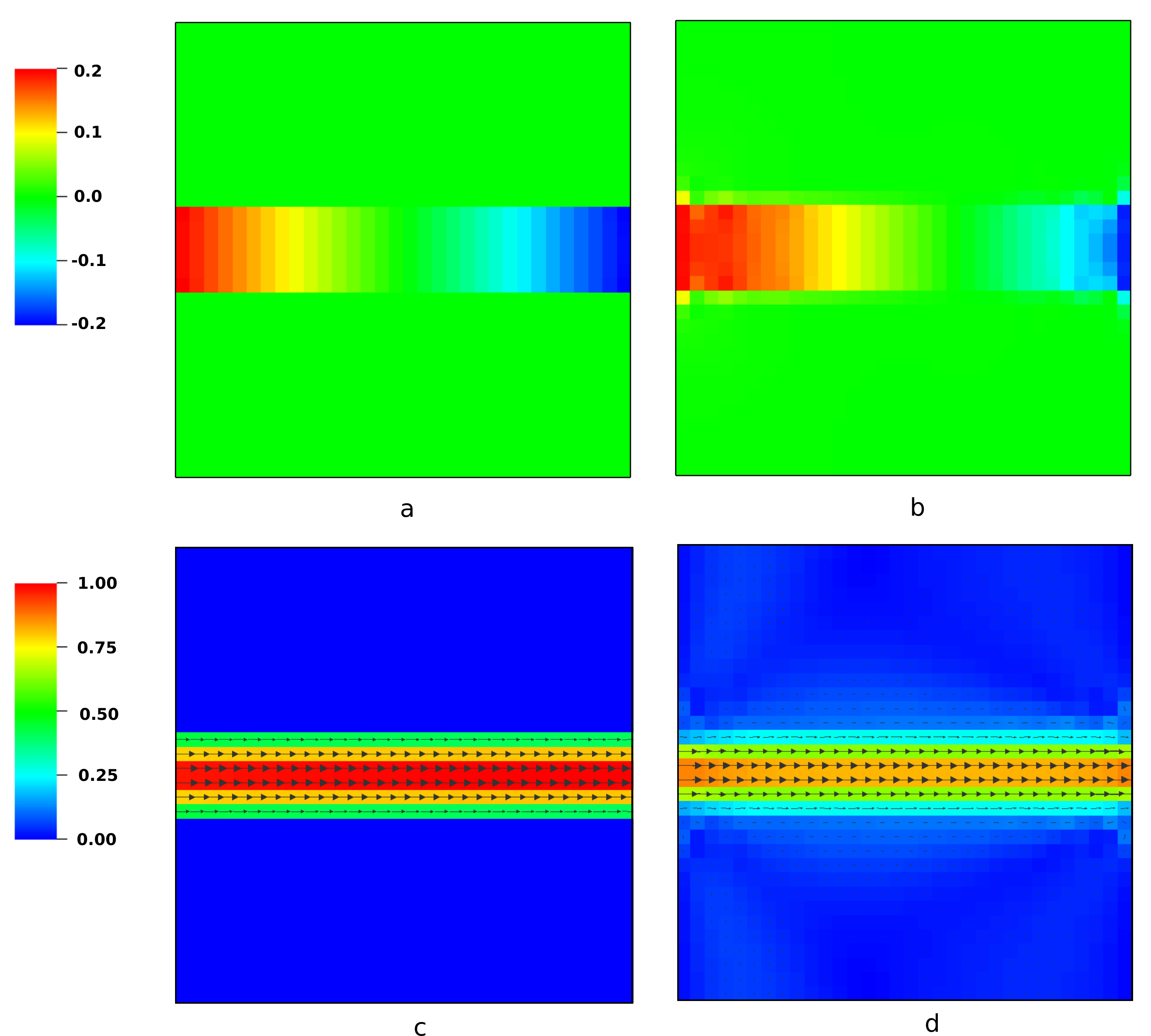}
		
		\caption{Comparison of the steady-state pressure (top panels) and velocity (bottom panels) distributions between the present IIM (left panels) and the standard IB method (right panels)
		     for the horizontal channel with grid spacing $h = 0.156$. The IIM sharply resolves both the pressure and flow, and yields higher flow rates at comparable grid spacings than the IB method.
		     See also Fig.~\ref{fig:compare_UP_profile_PPF}.}	
\label{fig:U_P_contours_straight_channel_2D_comparison} 
\end{figure}

\begin{figure}[t!]
		\centering
			\includegraphics[width=0.9\textwidth]{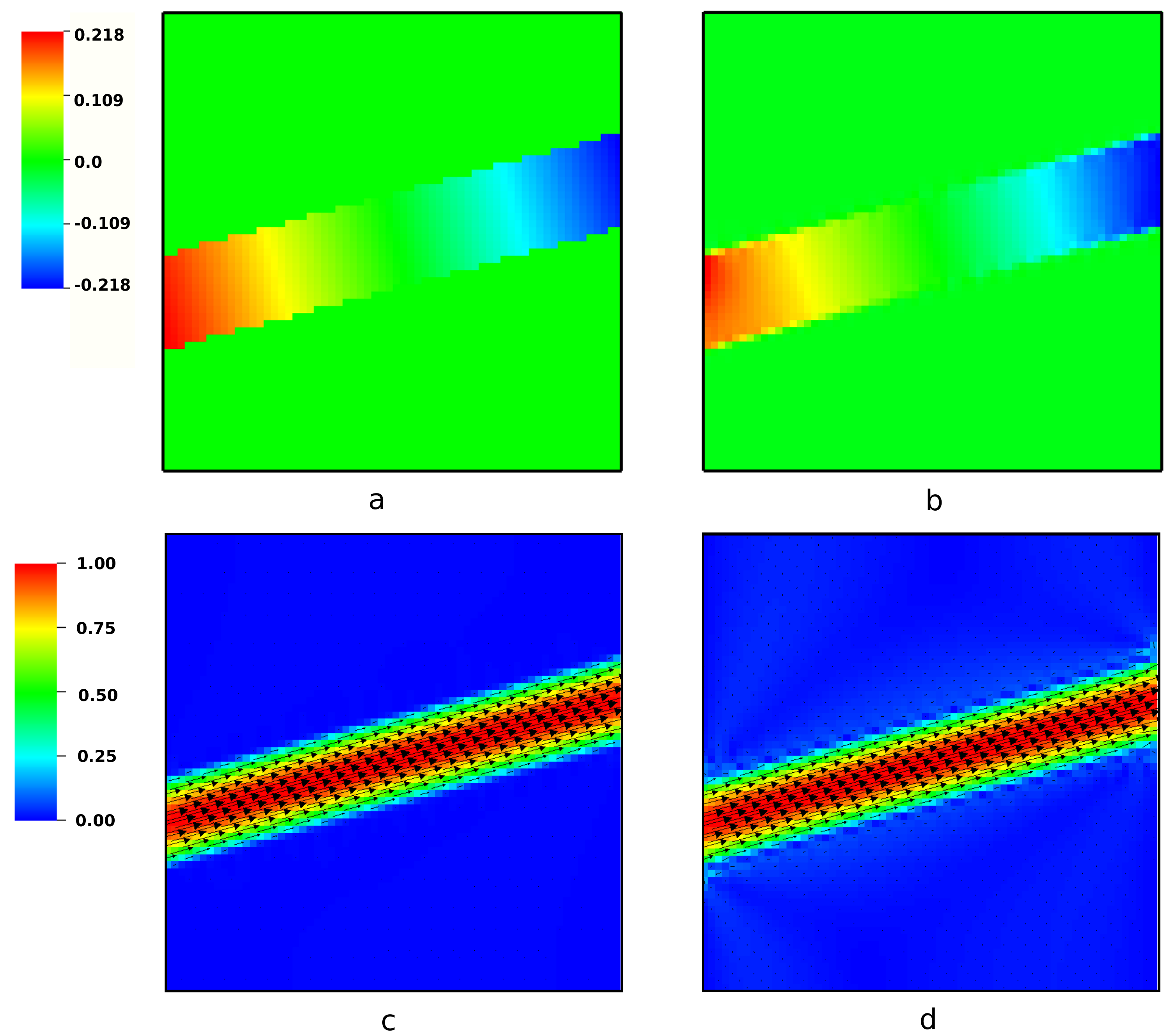}
		
		\caption{Comparison of the steady-state pressure (top panels) and velocity (bottom panels) distributions between the present IIM (left panels) and the standard IB method (right panels)
		     for the inclined channel with grid spacing $h=0.078$.  The results are qualitatively similar to those of the horizontal case shown in Fig.~\ref{fig:U_P_contours_straight_channel_2D_comparison}, but note that here we are using a finer computational grid, and so both the IIM and IB results are more accurate than those of Fig.~\ref{fig:U_P_contours_straight_channel_2D_comparison}.}	
\label{fig:U_P_contours_inclinde_channel_2D_comparison} 
\end{figure}

To provide a qualitative comparison of the present immersed interface method to the conventional IB method, Fig.~\ref{fig:U_P_contours_straight_channel_2D_comparison} shows the velocity and pressure in the horizontal channel for a relatively coarse grid spacing of $h=0.156$ using both numerical approaches.
Fig.~\ref{fig:U_P_contours_inclinde_channel_2D_comparison} is similar, but considers the skew configuration at a finer grid resolution of $h=0.078$.
As compared to the IB method, the present immersed interface scheme yields sharply resolved pressure and velocity fields.

\begin{figure}[t!!]
		\centering
			\includegraphics[width=0.9\textwidth]{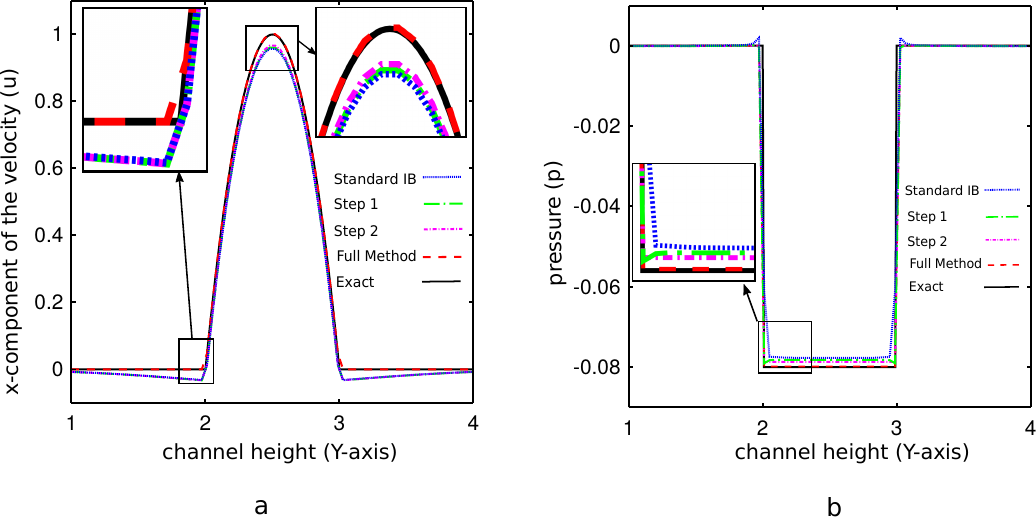}
		
		\caption{Comparisons of (a) velocity and (b) pressure profiles at $x=3.5$ for the plane-Poiseuille flow with grid spacing $h=0.078$.  Here, \textit{Standard IB} corresponds to a standard IB method with piecewise-linear regularized delta functions, \textit{Step 1} corresponds to a hybrid IB-IIM that only imposes jump conditions for the pressure, \textit{Step 2} corresponds to an IIM that imposes jump conditions for both the pressure and velocity gradient but uses uncorrected bilinear interpolation to determine the motion of the interface, and \textit{Full Method} corresponds to the IIM described in this paper.  Notice that the full method is substantially more accurate than the other methods considered here.}	
		\label{fig:compare_UP_profile_PPF}	
\end{figure}

To further analyze the improvement in accuracy offered by jump condition-based discretization approaches, Fig.~\ref{fig:compare_UP_profile_PPF} shows the incremental changes in the accuracy of the scheme as the jump conditions are systematically incorporated.
The horizontal channel with a grid spacing of $h=0.078$ is considered, and pressure and velocity profiles are plotted at $x=3.5$.
In the figure, \textit{Step 1} indicates results obtained by a simplified version of the present method in which only the normal component of the interfacial force is used to determine pressure jump conditions.
The tangential component of the force is transmitted to the background grid using a standard IB approach with a piecewise linear kernel function.
In addition, the velocity is also interpolated to the interface using the conventional IB approach with a piecewise linear kernel function.
Thus, in Step 1, only pressure discontinuities are sharply resolved, and discontinuities in the velocity gradient are regularized.
This scheme is similar to the IIM of Lee and LeVeque \cite{lee2003immersed}.  
The result demonstrates that adding the correction due to the pressure discontinuity yields a significant improvement in the accuracy of the pressure as compared to the conventional IB method, but the velocity profile at the immersed boundary still suffers from low accuracy.
\textit{Step 2} indicates results obtained by a version of the simplified method that is modified so that the tangential portion of the force is also used to impose jump conditions in the velocity gradient along the interface, but the velocity is still interpolated to the interface using the conventional IB approach with a piecewise linear kernel function.
Notice that there remains a sizable difference between the Step 2 velocity profile and the exact solution, which results from the use of uncorrected interpolation to determine the velocity of the interface.
\textit{Full Method} indicates results obtained by the full immersed interface method presented in this work.
It is clear that if all jump conditions are included both in applying the forces and in interpolating the velocities, the resulting method is in excellent agreement with the analytical solutions for both velocity and pressure.
The accuracy of the scheme is degraded by omitting any of these corrections.

\begin{figure}[t!!]
		\centering
			\includegraphics[width=0.5\textwidth]{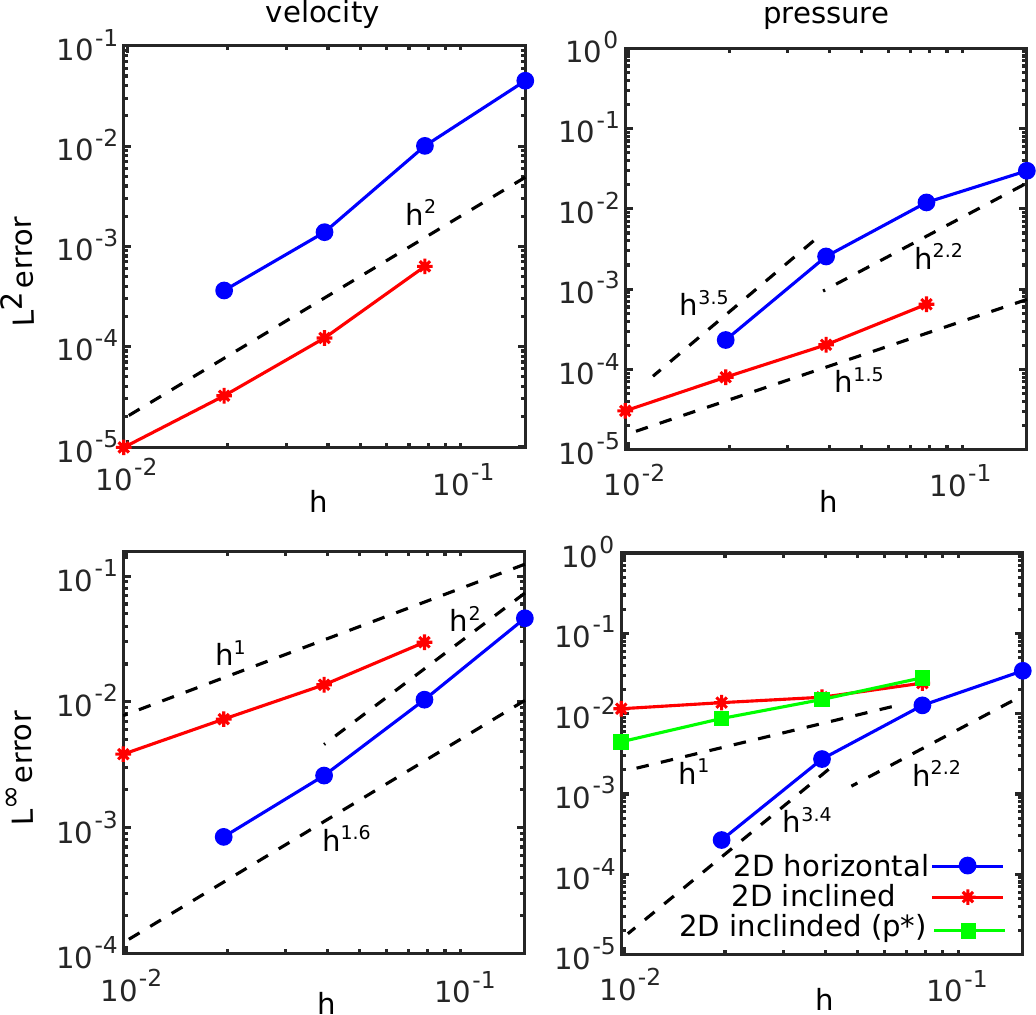}
		
		\caption{Convergence of the $\L2$ and $\Linf$ norms of the error of the Eulerian velocity and pressure for the horizontal and inclined plane-Poiseuille flows.  
		Simulation parameters include $\textrm{Re}=66.66$, $\Delta t=0.1h$, and $\mfac=2$.
		}
		\label{fig:convergence_ch2d_Eulerian} 
\end{figure}

\begin{figure}[t!!]
		\centering
			\includegraphics[width=0.95\textwidth]{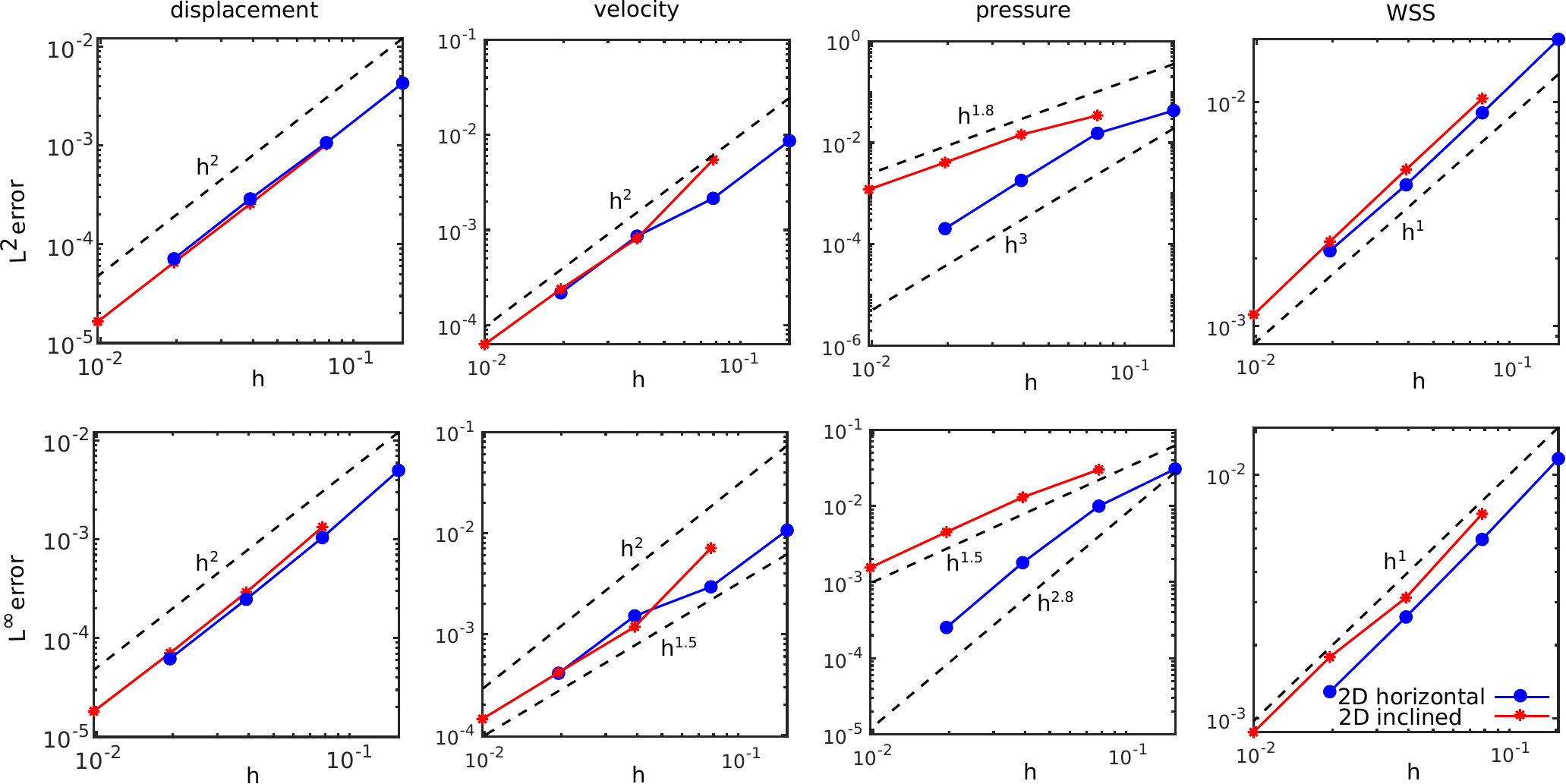}
		
		\caption{Convergence of the $\L2$ and $\Linf$ norms of the errors of the Lagrangian displacement, velocity, pressure, and wall shear stress (WSS) for the horizontal and inclined plane-Poiseuille flow, determined along the lower side of the channel.
		  Simulation parameters include  $\textrm{Re}=66.66$, $\Delta t=0.1h$, and $\mfac=2$.
		  }
		\label{fig:convergence_ch2d_Lagrangian} 
\end{figure}

Grid convergence studies are performed for both the horizontal and inclined configuration.
The errors in the computed Eulerian velocity and pressure in $\L2$ and $\Linf$ norms are summarized in Fig.~\ref{fig:convergence_ch2d_Eulerian}. 
These errors are also tabulated and reported in \ref{sec:appendix}.
Second-order convergence of the velocity is observed in the $\L2$ norm for both the horizontal and inclined cases, and the global order of accuracy in the $\L2$ norm of the pressure appears to be between first-~and second-order.
As already discussed, errors and convergence rates for discontinuous quantities such as the pressure are sensitive to the precise location of the discontinuity, and small perturbations in the position of the interface can result in large errors when compared to the analytic solution obtained with the exact interfacial geometry.
Consequently, we report $\Linf$ errors for both the full Eulerian pressure field $p$ and also for the field $p^*$ that is obtained by excluding grid cells within $2h$ of the interface.
Although the convergence of $p$ in the $\Linf$ norm is slow, $p^*$ shows first-order pointwise convergence rates.

Fig.~\ref{fig:convergence_ch2d_Lagrangian} reports the errors in the computed Lagrangian displacement, velocity, pressure, and wall shear stress, and tabulated errors are provided in \ref{sec:appendix}.
At least second-order convergence is observed for the displacement in the $L^2$ and $L^\infty$ norms and for the velocity in the $\L2$ norm.
The $\Linf$ norm of the velocity error in the inclined channel test appears to converge at slightly less than second order.
Also at least $1.5$-order accuracy is achieved in the pressure in the $\L2$ and $\Linf$ norms along with first-order convergence in the wall shear stress in the $\L2$ and $\Linf$ norms.

\subsubsection{Hagen-Poiseuille flow}
\label{subsubsec:Hagen-poiseuille}

\begin{figure}[t!!]
		\centering
			\includegraphics[width=0.9\textwidth]{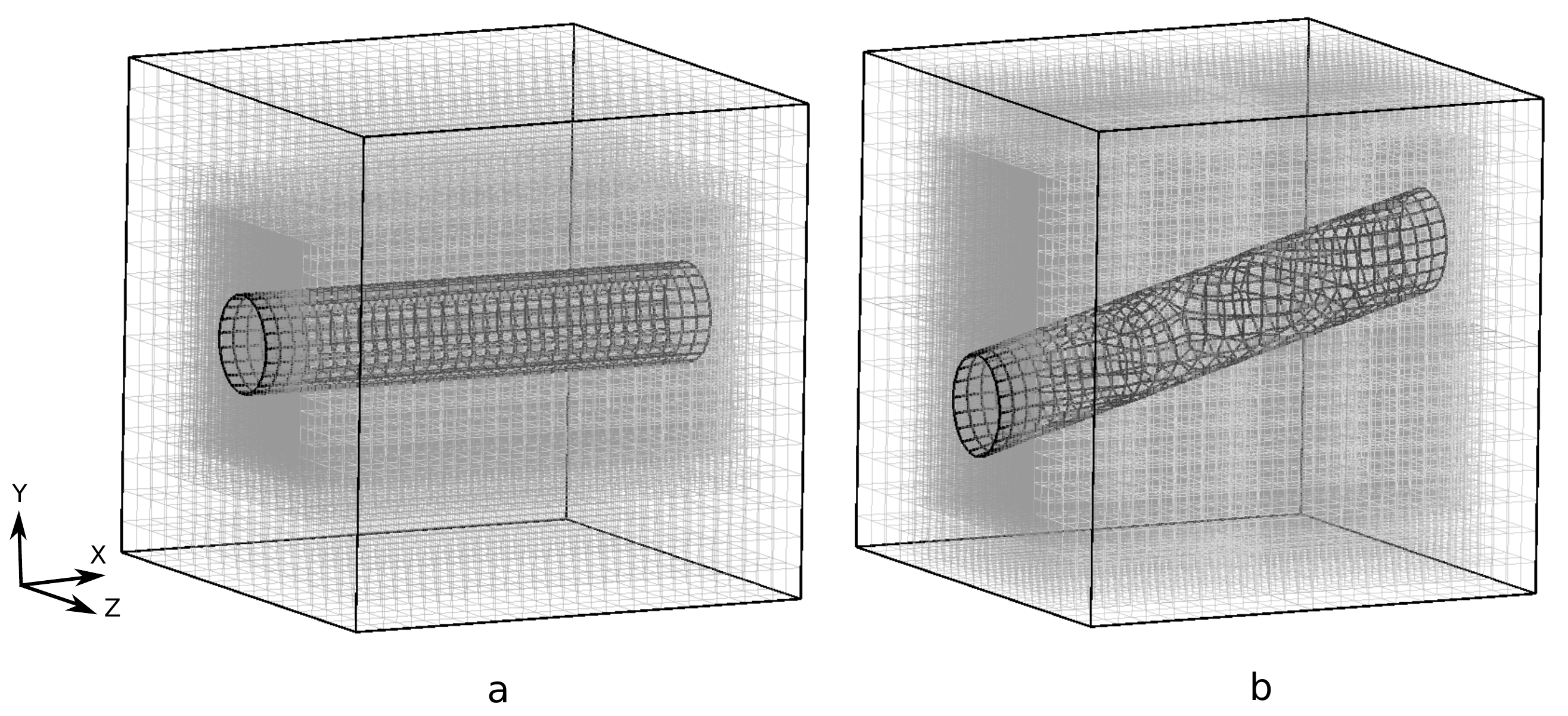}
		
		\caption{Computational meshes used to compute flow in (a) the horizontal pipe and (b) the inclined pipe.  A locally refined Cartesian grid is used in both cases. }
\label{fig:mesh-3D-channel} 
\end{figure}

The steady-state solution for a three-dimensional axisymmetric Hagen-Poiseuille flow in the $x$-direction is
\begin{align}
	\label{eq:PPF_u_3D} u(r)&=\frac{p_0}{2\mu L}\left(R^{2}-r^{2}\right),\\
	\label{eq:PPF_v_3D} v&=0, \\
	\label{eq:PPF_w_3D} w&=0, \\
	\label{eq:PPF_p_3D} p(x)&=p_{0}-2p_{0}x/L,
\end{align}
in which $r=\sqrt{(x-x_{0})^{2}+(y-y_{0})^{2}}$ and $R$ is the radius of the pipe.
The physical domain is $\Omega=[0,L]^3$.
A horizontal pipe of length $L$ extends across the middle of the domain and is described using quadrilateral surface elements with $\Mfac = 2$, see Fig.~\ref{fig:mesh-3D-channel}(a).
A pressure difference of $2p_0$ is applied across the pipe. 
As in the two-dimensional case, the boundary conditions at the inlet and outlet are imposed through applying a combination of normal-traction and zero-tangential-slip boundary conditions \cite{griffith2009}.
Solid-wall boundary conditions are imposed along the remaining parts of $\p\Omega$.
The computational domain is discretized using an adaptively refined Cartesian grid.
With $N$ grid levels and a refinement ratio of two between levels, the Cartesian grid spacing on the finest grid level is $h_\text{finest}=2^{-(N-1)} h_{\text{coarsest}}$, with $h_{\text{coarsest}}=\frac{L}{8}$ being the grid spacing on the coarsest level. 
The time step size is $\Delta t=0.05 h_{\text{finest}}$, giving a maximum advective CFL number of approximately $0.05$--$0.1$, and the spring stiffness is set to $\kappa=2.2\times 10^{-3}/\Delta t^2$.
Other material and simulation parameters are $L=5$, $R=0.5$, $\rho=1$, $\mu=0.01$, and $p_0=0.4$, resulting in a Reynolds number of $\Re=100$. 

As in the plane flow case, we also consider an inclined configuration.
In the skew case, the structure is placed in the same cubic Cartesian domain but is rotated around the $z$-axis with an angle of $\theta=\pi/12$; see Fig.~\ref{fig:mesh-3D-channel}(b). 
At the inlet and outlet, we specify normal traction boundary conditions and tangential velocity boundary conditions that are consistent with the exact solution, and solid wall boundary conditions are imposed along the remainder of $\p\Omega$.
Other simulation parameters are the same as the grid-aligned case.  

\begin{figure}[t!!]
		\centering
			\includegraphics[width=0.9\textwidth]{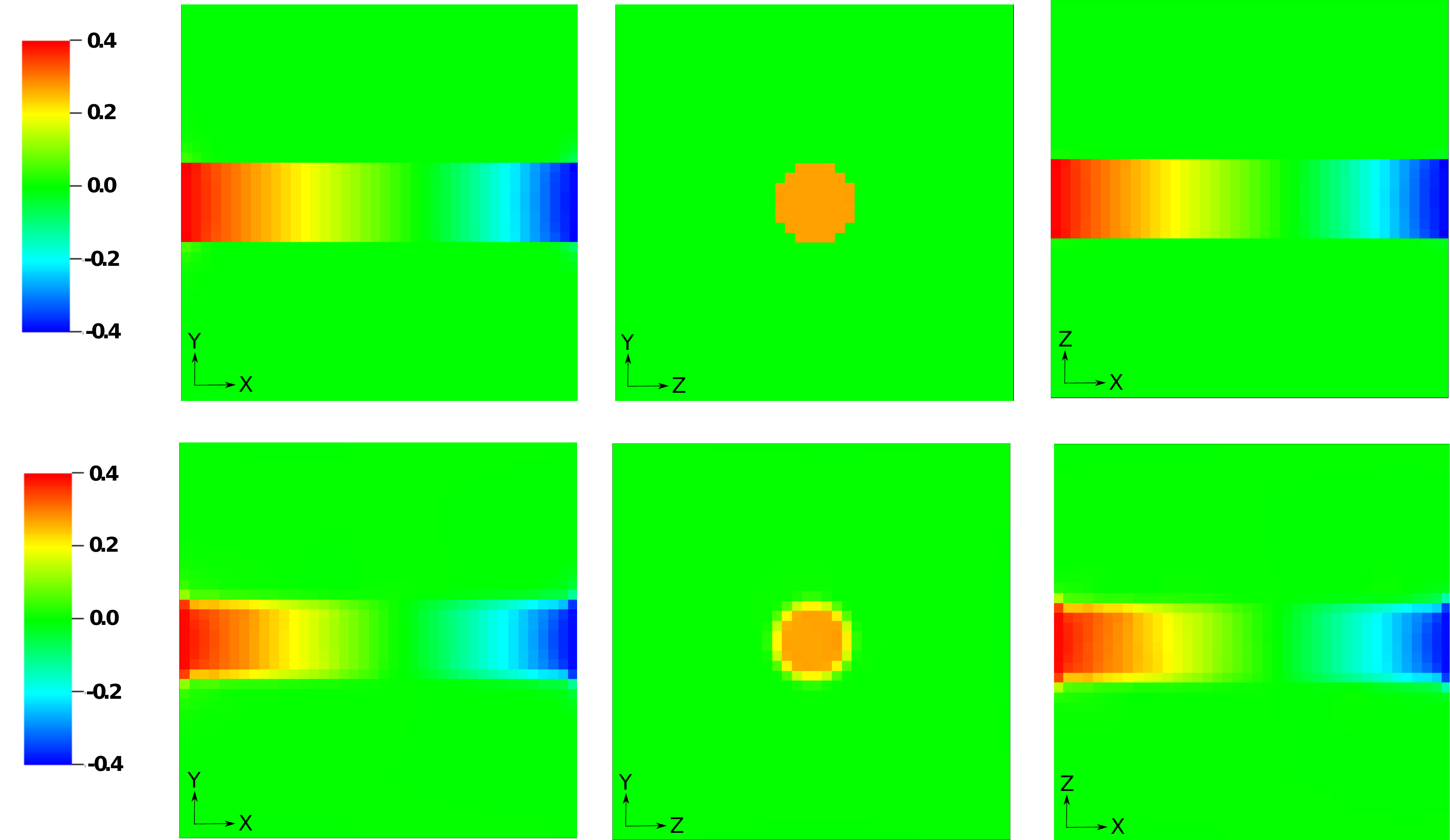}
		
		\caption{ Comparison of the steady-state pressure field along the planes $z=2.5$, $y=2.5$, and $x=1$ obtained using the present IIM (top panel) and the standard IB method (bottom panel)
					with a relatively coarse locally refined grid comprising $N=3$ levels. With the IB method,
		          spurious pressure oscillations are clearly observed, whereas the interface method yields a sharp pressure profile throughout the computational domain.}    
\label{fig:UP_HPF_channel_3D} 
\end{figure}

\begin{figure}[t!!]
		\centering
			\includegraphics[width=0.9\textwidth]{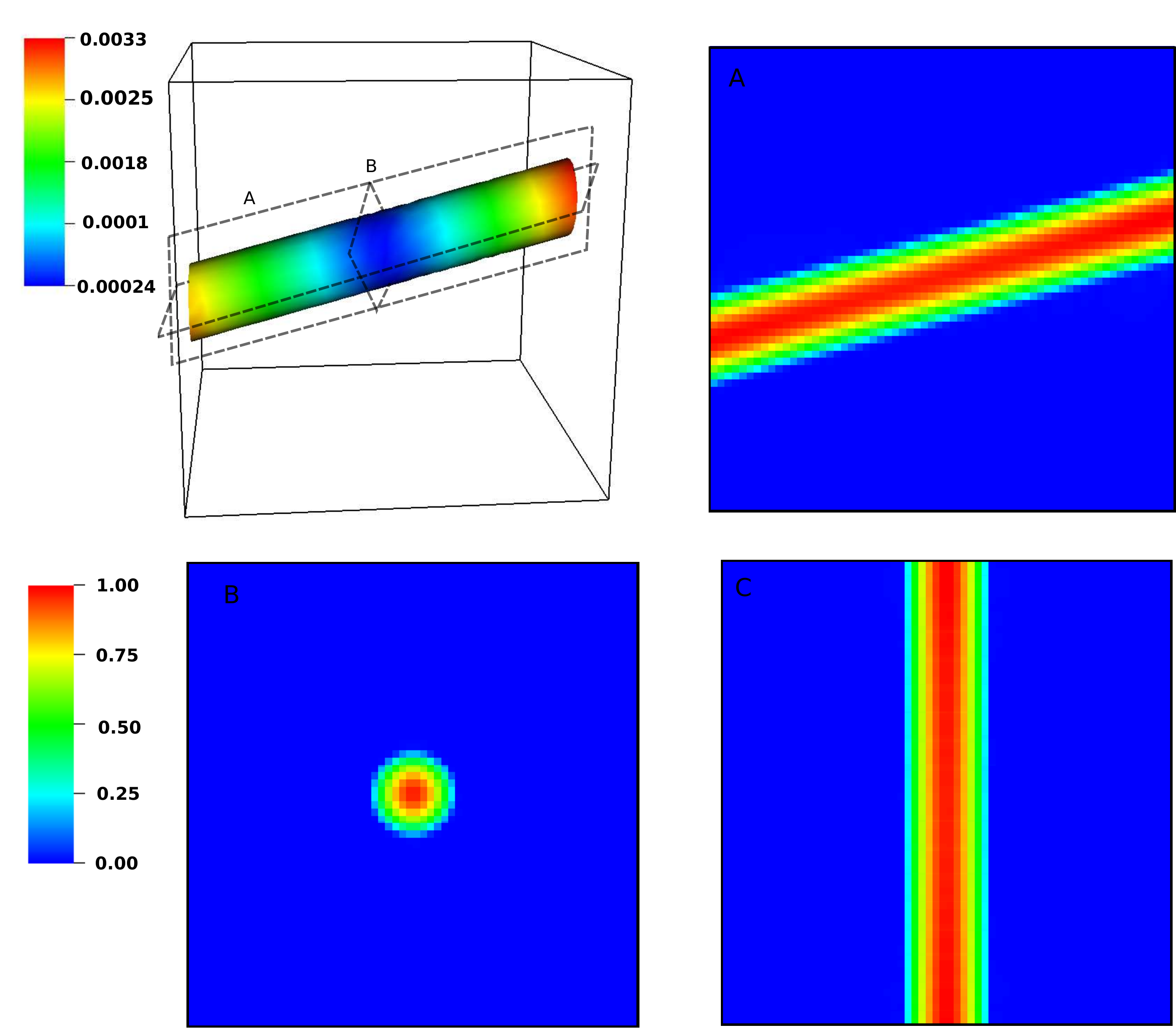}
		
		\caption{Maximum structural displacement magnitude for the three-dimensional slanted channel (top left panel). Section views of velocity magnitudes at different planes are shown in panels A, B, and C.
		       In this computation, the effective fine grid resolution is $h_\text{finest} = 0.156$, and the maximum structural displacement error is approximately $3.3\times 10^{-3}$, which is approximately a factor of 50 \textit{smaller} than $h_\text{finest}$.}	
\label{fig:disp-slanted-3d-vel} 
\end{figure}

To illustrate the improvement of the present IIM over the conventional IB method, Fig.~\ref{fig:UP_HPF_channel_3D} shows the pressure fields obtained using both methods on a relatively coarse locally refined grid with $N=3$ levels, which yields an effective grid spacing of $h = 0.156$. 
Cross sectional views of the planes $z=2.5$, $y=2.5$, and $x=1$ are shown.
The IIM clearly yields a more sharply resolved pressure field.
The solution accuracy of the velocity using the present immersed interface approach is qualitatively evaluated in Fig.~\ref{fig:disp-slanted-3d-vel}.
As in the two-dimensional case, the IIM sharply resolves the velocity field.
This figure also examines the maximum structural displacement magnitude on the surface of the cylinder.
The maximum value is on the order of $3.3\times 10^{-3}$, which is approximately $1/50$ of the Cartesian grid cell width $h_\text{finest}$.

\begin{figure}[t!!]
		\centering
			\includegraphics[width=0.5\textwidth]{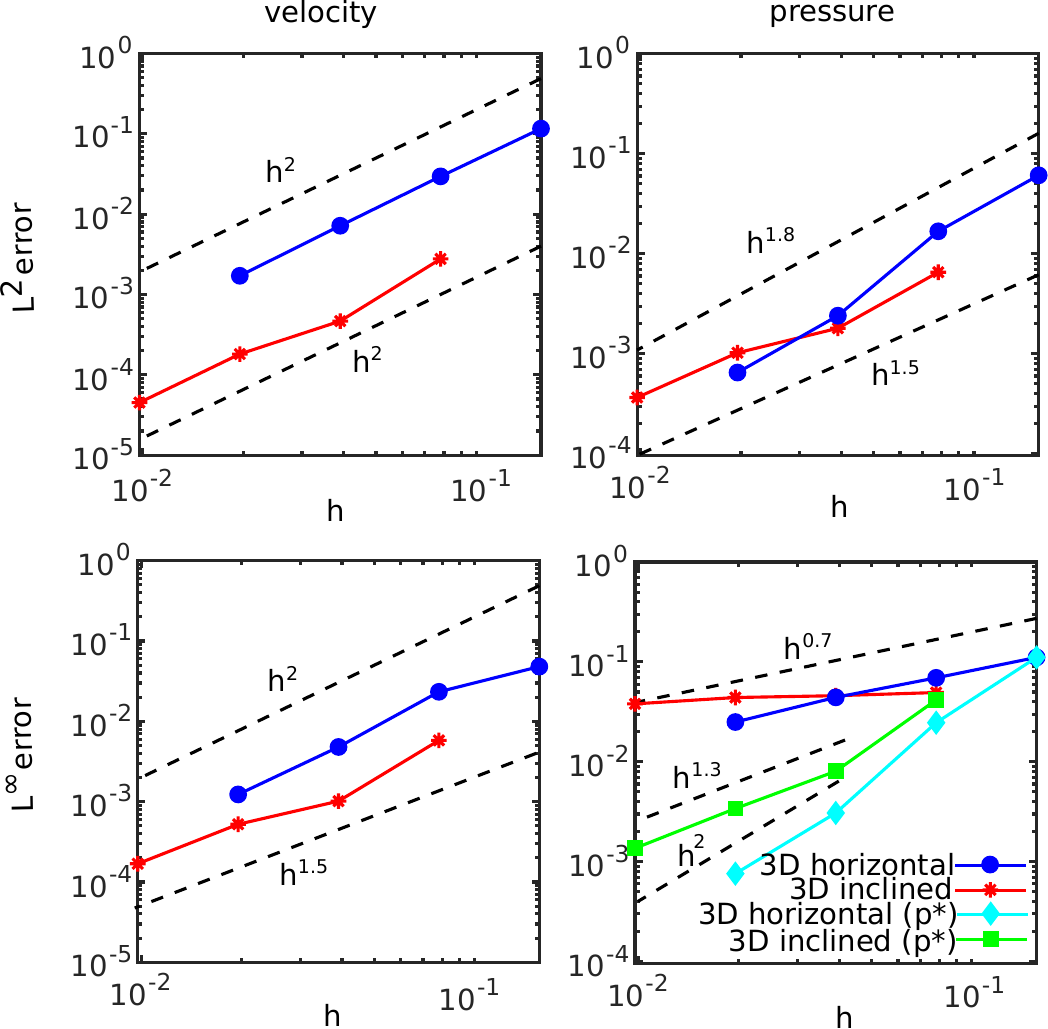}
		
		\caption{Convergence in the $\L2$ and $\Linf$ norms of the error in the Eulerian velocity and pressure for the horizontal and inclined 
		three-dimensional Hagen-Poiseuille flow. Simulation parameters include $\Re=100$, $\Delta t=0.1h_\text{finest}$, and $\mfac=2$.}
		\label{fig:convergence_ch3d_Eulerian} 
\end{figure}

\begin{figure}[t!!]
		\centering
			\includegraphics[width=0.95\textwidth]{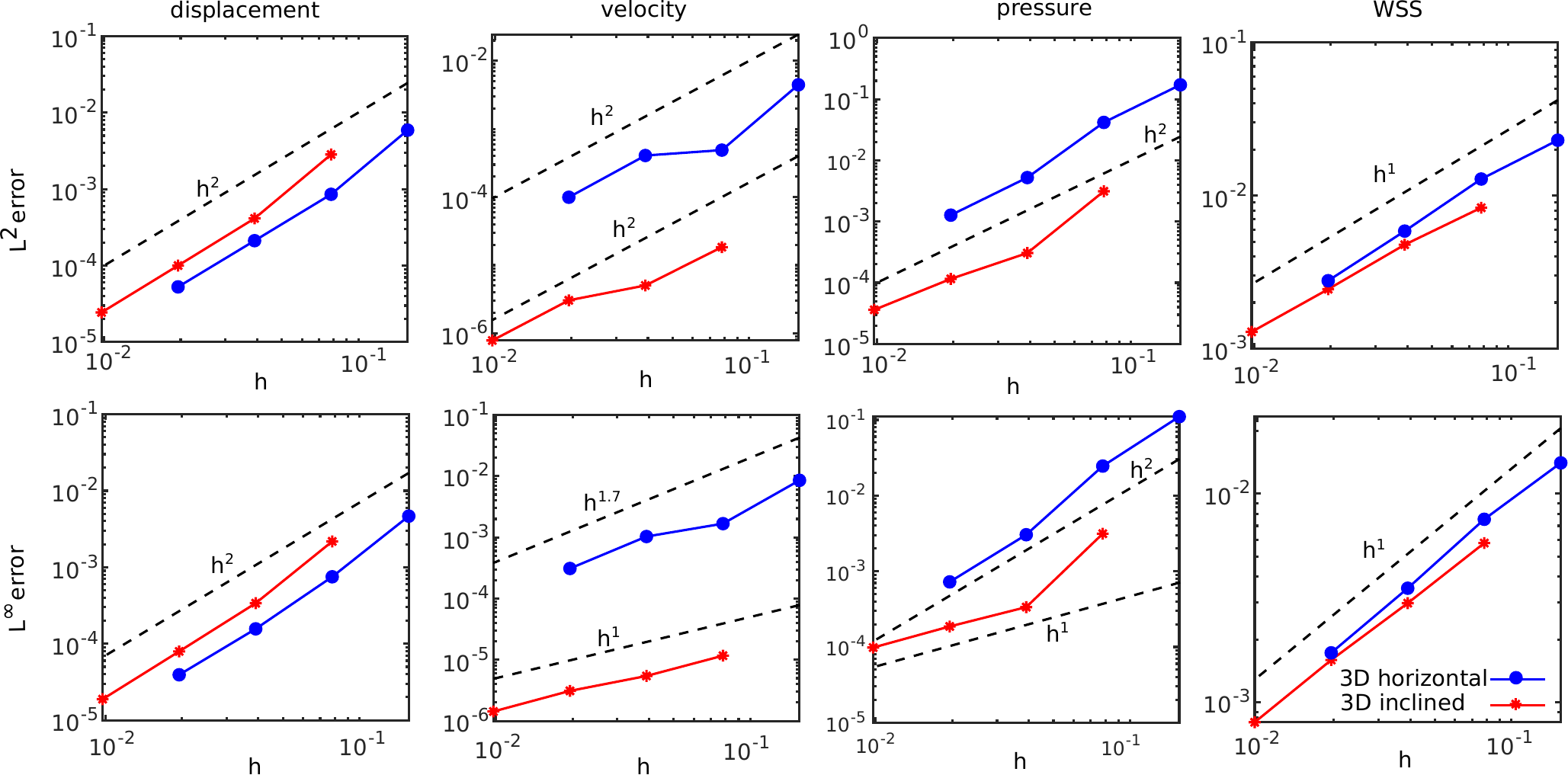}
		
		\caption{Convergence in the $\L2$ and $\Linf$ norms of the error in the Lagrangian displacement, velocity, pressure, and wall shear stress (WSS) for the horizontal and inclined three-dimensional Hagen-Poiseuille flow. 
		Simulation parameters include $\textrm{Re}=100$, $\Delta t=0.1h_\text{finest}$, and $\mfac=2$.}
		\label{fig:convergence_ch3d_Lagrangian} 
\end{figure}

Finally, a convergence study is performed by systematically increasing the number $N$ of AMR levels.
Fig.~\ref{fig:convergence_ch3d_Eulerian} reports the $\L2$ and $\Linf$ errors in Eulerian velocity and pressure. 
The method yields second-order convergence rates for the velocity in the $\L2$ norm for both grid-aligned and skew configurations. 
Slightly less than second-order accuracy is observed in the $\L2$ norm of the pressure, and although less than first-order convergence is observed for $p$ in the $\Linf$ norm over the entire domain, first-order pointwise convergence is recovered in the restricted domain $\Omega^*$.
Fig.~\ref{fig:convergence_ch3d_Lagrangian} summarizes the convergence results on the Lagrangian mesh.
Second-order convergence rates are achieved for the displacement in both $\L2$ and $\Linf$ error norms.
Second-order convergence is also observed for the Lagrangian displacement and velocity in the $\L2$ norm, whereas between first-~and second-order convergence is observed for the velocity in the $\Linf$ norm.
The Lagrangian pressure converges at approximately second order in the $\L2$ norm and between first-~and second-order in the $\Linf$ norm.
First-order convergence is observed for the wall shear stress in both the $\L2$ and $\Linf$ norms.
See \ref{sec:appendix} for tabulated errors of this example.

\subsection{Circular Couette flow in two and three dimensions}
\label{subsec:TCF}

\begin{figure}[t!!]
		\centering
			\includegraphics[width=0.5\textwidth]{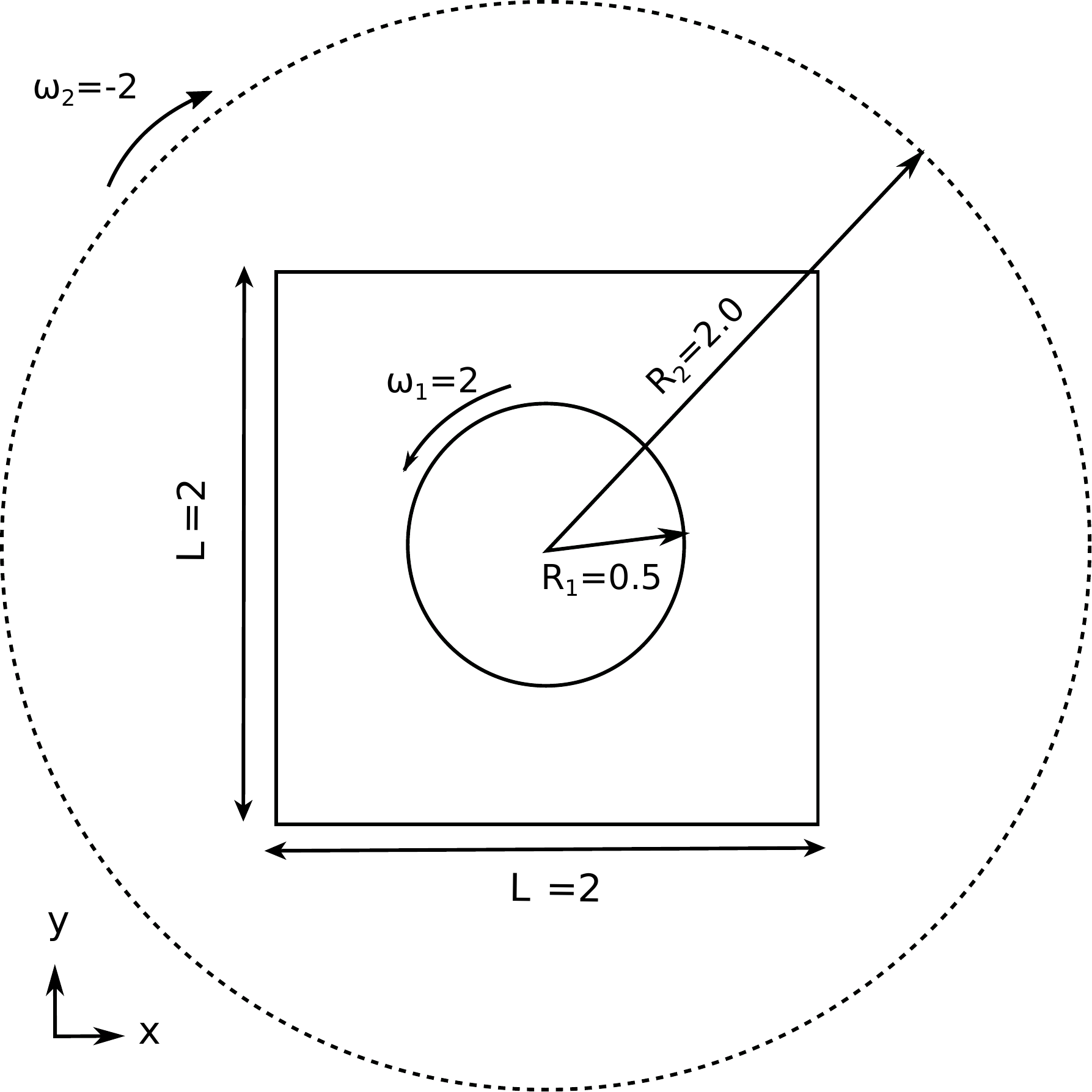}
		
\caption{Geometry and dimensions of the circular Couette flow.}
\label{fig:TCF_schematic} 
\end{figure}

This section considers the steady shear-driven Couette flow of a fluid confined between two rotating circular cylinders in two and three spatial dimensions.
This problem is an excellent example to test the tangential portion of the interfacial force on a curved geometry. 
Fig.~\ref{fig:TCF_schematic} provides a schematic of the geometry and dimensions of the circular Couette flow in two dimensions 
Analytical solutions for the steady-state velocity and the pressure for the two cylinders co-centered at $(x_{0},y_{0})$ are
\begin{equation}
u(x,y)=\begin{cases}
-\omega_{1}(y-y_0), & 0 \leq r\leq R_{1},\\
-(y-y_{0})(A+\frac{B}{r^{2}}), & R_{1}<r<R_{2},
\end{cases} 
\end{equation}
\begin{equation}
v(x,y)=\begin{cases}
\omega_{1}(x-x_{0}), & 0 \leq r\leq R_{1},\\
(x-x_{0})(A+\frac{B}{r^{2}}), & R_{1}<r<R_{2},
\end{cases}
\end{equation}
and
\begin{equation}
p(x,y)=\begin{cases}
\frac{\omega^2_{1}r^{2}}{2} + p_0, & 0 \leq r\leq R_{1},\\
\frac{A^{2}r^{2}}{2}-\frac{B^{2}}{2r^{2}}+AB\cdot\ln{r^{2}}, & R_{1}<r<R_{2},
\end{cases}
\end{equation}
in which $A=\frac{\omega_{2}R_{2}^{2}-\omega_{1}R_{1}^{2}}{R_{2}^{2}-R_{1}^{2}}$, $B=\frac{(\omega_{1}-\omega_{2})R_{1}^{2}R_{2}^{2}}{R_{2}^{2}-R_{1}^{2}}$,
$r=\sqrt{(x-x_{0})^{2}+(y-y_{0})^{2}}$, and $p_0$ is an arbitrary constant value.
We choose $R_1=0.5$, $R_2=2.0$, $\omega_{1}=2$, and $\omega_{1}=-2$.
The computational domain is $\Omega = [-1,1]^2$, and $(x_{0},y_{0})=(0,0)$ so that only the inner cylinder is embedded in the physical domain. 
Dirichlet conditions for the velocity are imposed along $\partial \Omega$ using the  analytic relation.
The time step size is set to $\Delta t=0.05h$, yielding a CFL number between $0.05$ and $0.1$. 
The penalty parameters are chosen to be $\kappa=7\times 10^{-3}/\Delta t^2$ and $\eta=0$. 
Other simulation parameters are $\Mfac = 2$, $\rho=1$, and $\mu=0.01$.
The Reynolds number, based on the velocity and diameter of the inner circle, is $\Re=\frac{2\rho R^2_{1}\omega_{1}}{\mu}=100$.

A three-dimensional extension of the  problem is also considered.  
The same simulation parameters are used as in the two-dimensional case, except that here, the domain is $\Omega = [-1,1]\times[-1,1]\times[-2,2]$.
Periodic boundary condition are applied in the $z$-direction. 
The computational domain is discretized using an $N$-level locally refined Cartesian grid with $h=h_{\text{finest}}=2^{-(N-1)}h_{\text{coarsest}}$ and $h_{\text{coarsest}}=\frac{L}{8}$ being the coarsest level.  
Because $\Re$ is relatively low, we do not observe Taylor vortices, but the additional dimension helps to further verify the implementation and consistency of results in two and three spatial dimensions.

\begin{figure}[t!!]
		\centering
			\includegraphics[width=0.95\textwidth]{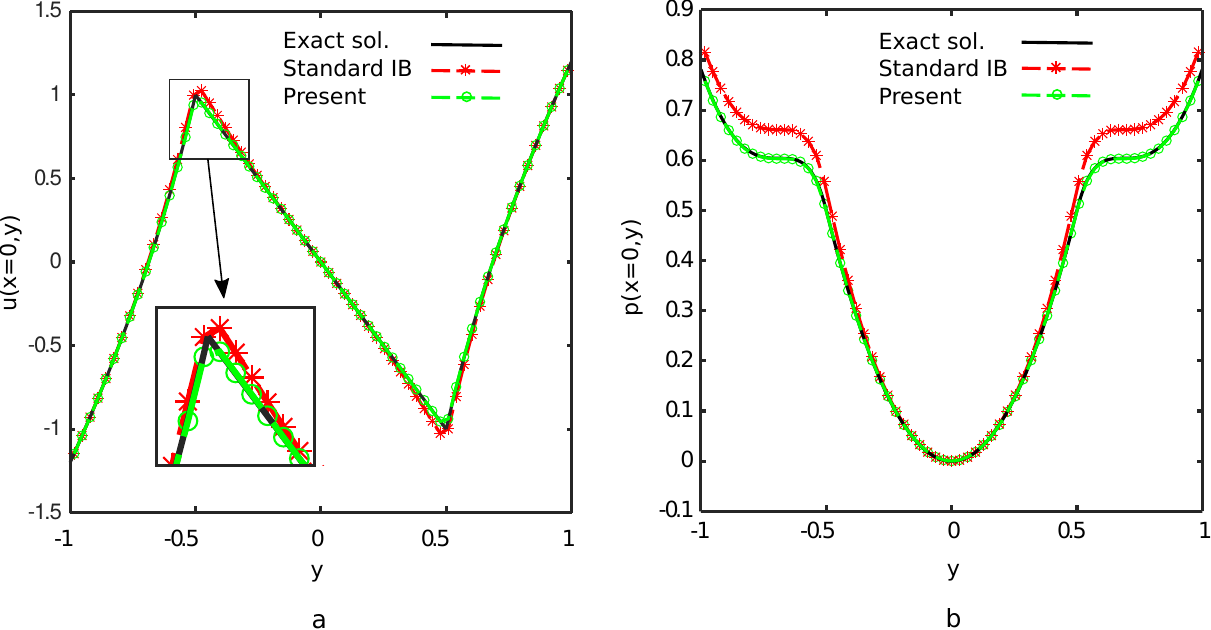}
		
\caption{(a) $u$-velocity and (b) pressure profiles for the IIM and IB solutions to the two-dimensional Couette flow at $\Re=100$ and $h=0.0625$. The present IIM 
is in outstanding agreement with the analytic solution.}
\label{fig:u-p-profile-TCF-2d} 
\end{figure}

Fig.~\ref{fig:u-p-profile-TCF-2d} compares the performance of the present IIM and the IB method for the two-dimensional steady-state velocity and pressure profiles of the present approach along $x=0$ using grid spacing $h=0.0625$.
The agreement between the analytic solution and the results generated by the IIM is outstanding.

\begin{figure}[t!!]
		\centering
			\includegraphics[width=0.5\textwidth]{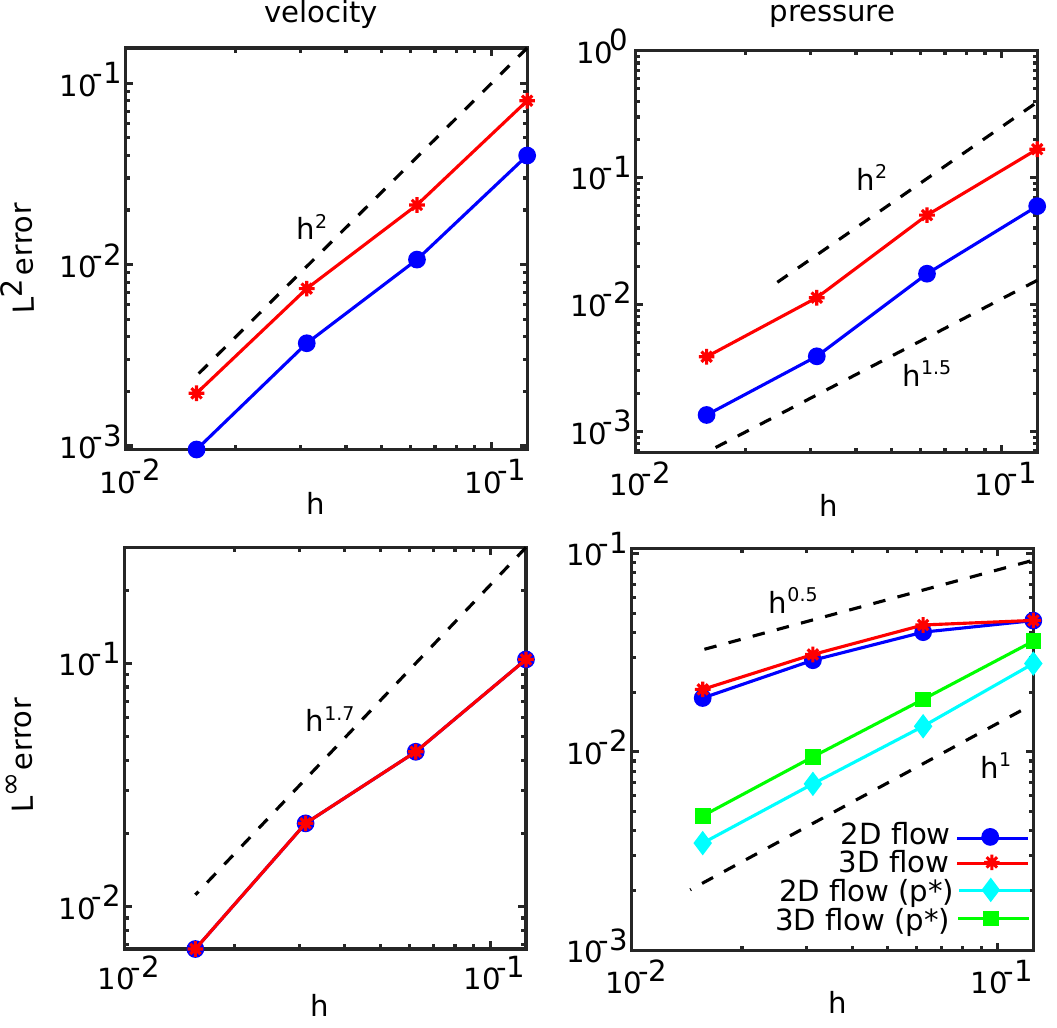}
		
			\caption{Convergence in the $\L2$ and $\Linf$  norms of the error in the Eulerian velocity and pressure for the two-~and three-dimensional circulat Couette flows.
	Simulation parameters include $\textrm{Re}=100$, $\Delta t=0.05h_\text{finest}$, and $\mfac=2$.}
		\label{fig:convergence_TCF_Eulerian} 
\end{figure}
\begin{figure}[t!!]
		\centering
			\includegraphics[width=0.95\textwidth]{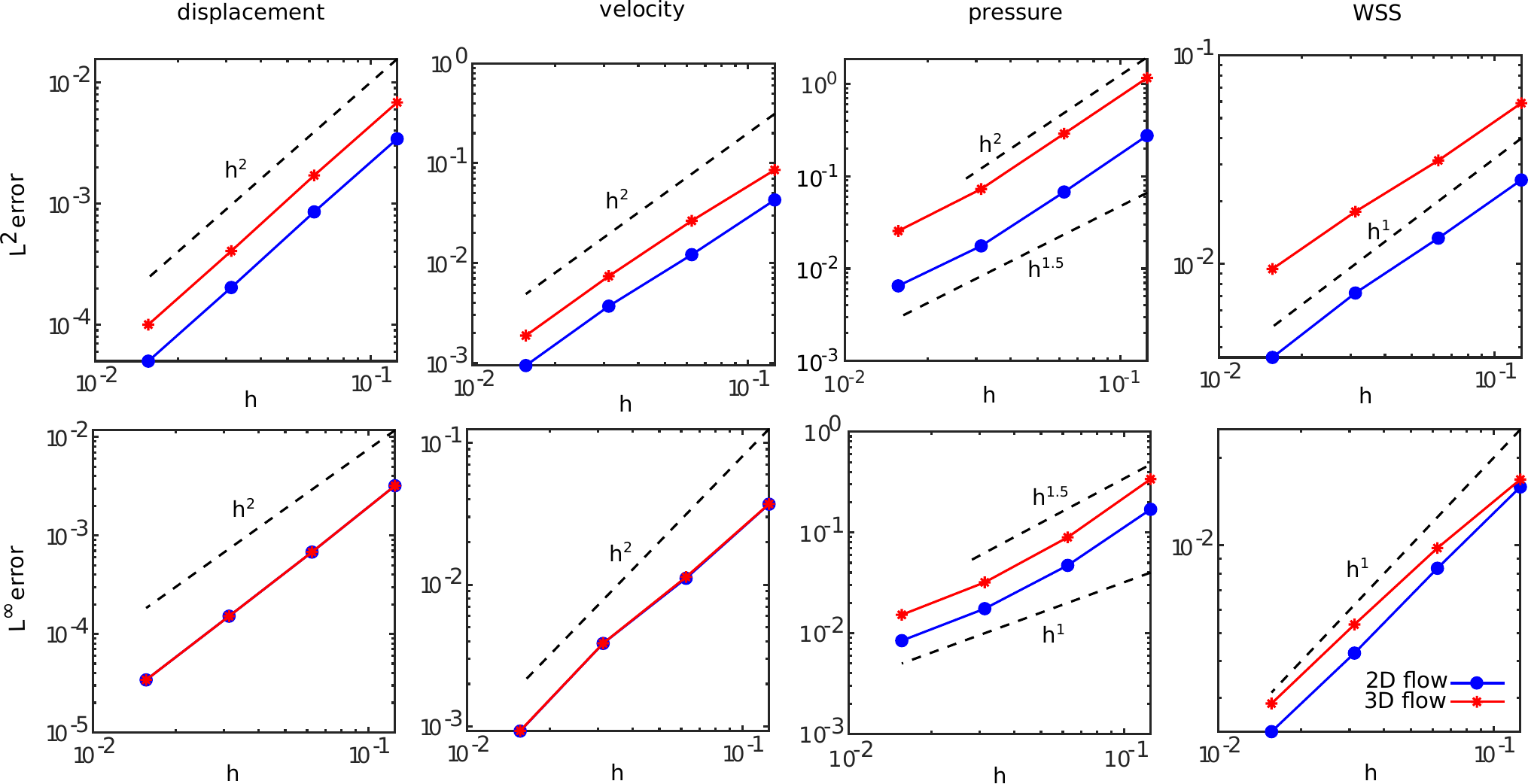}
		
		\caption{Convergence in the $\L2$ and $\Linf$ norms of the Lagrangian displacement, velocity, pressure, and wall shear stress (WSS) 
		for the two-~and three-dimensional circulat Couette flows. The Lagrangian pressure and wall shear stress 
		are computed from the side between the two cylinders.
		Simulation parameters include $\textrm{Re}=100$, $\Delta t=0.05h_\text{finest}$, and $\mfac=2$.}
		\label{fig:convergence_TCF_Lagrangian} 
\end{figure}

We also perform a convergence study for both the two-~and three-dimensional cases.
As seen in Fig.~\ref{fig:convergence_TCF_Eulerian} the method successfully achieves global second-order convergence rate of the Eulerian velocity in both two and three dimensions while yielding slightly less than second-order accuracy in the $\Linf$ norm.
At least $1.5$-order global convergence rates are obtained for the Eulerian pressure along with first-order convergence rates for the Eulerian pressure on $\Omega^*$.
The method also achieves pointwise second-order accuracy in the Lagrangian displacement and velocity along with between first-~and second-order pointwise convergence in the Lagrangian pressure and wall shear stress as seen in Fig.~\ref{fig:convergence_TCF_Lagrangian}. Note that the Lagrangian pressure and wall shear stress 
are computed from the side between the two cylinders.


\subsection{Flow within eccentric rotating cylinders}
\label{subsec:eccentric_cylinders}

\begin{figure}[t!!]
		\centering
			\includegraphics[width=0.5\textwidth]{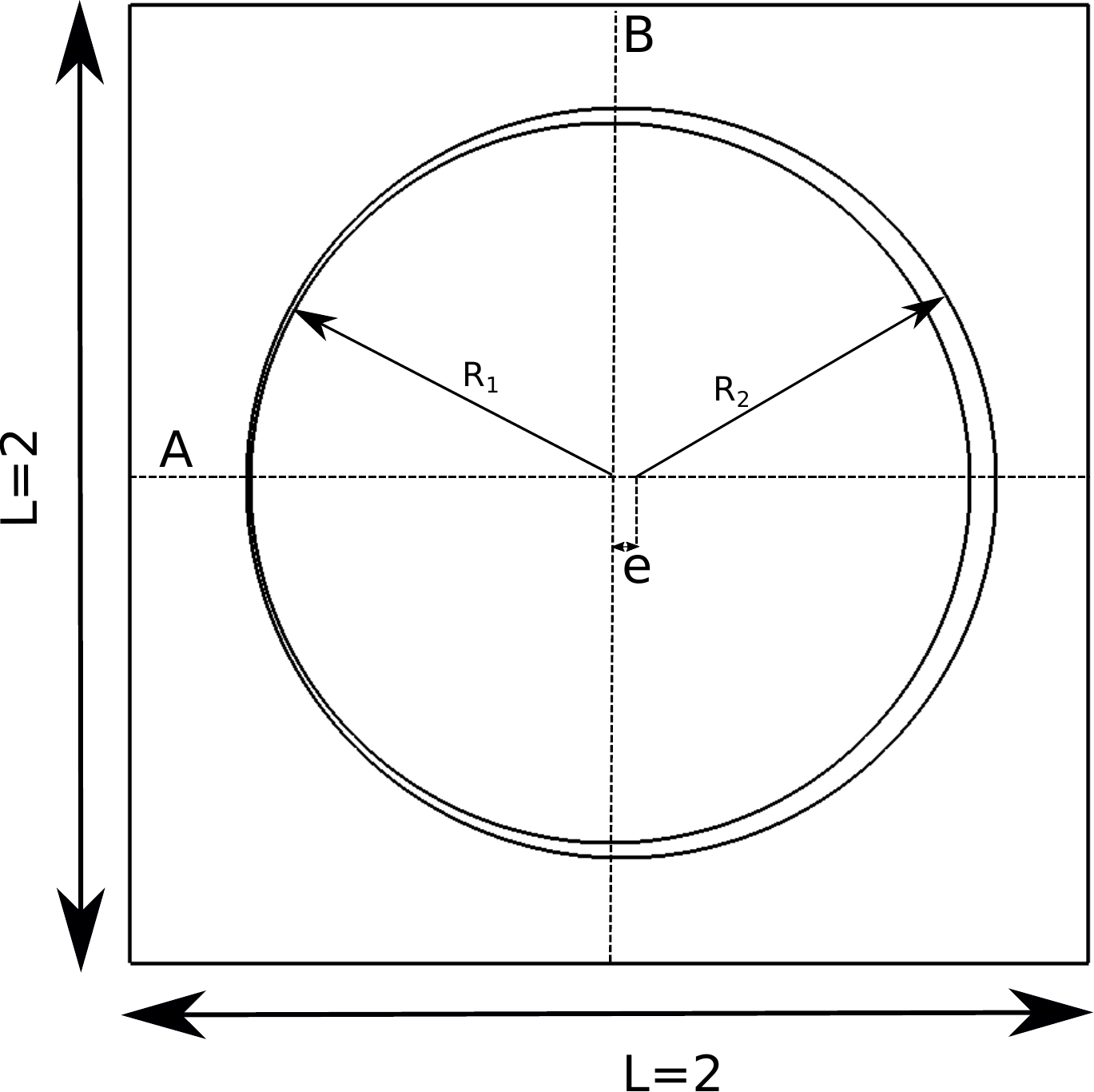}
		
\caption{Geometry and dimensions of the eccentric rotating cylinder test. The inner cylinder is centered at the origin, and the outer one
is centered at $(e, 0)=(3/128,0)$.
}
\label{fig:eccentric_cylinders_schematic} 
\end{figure}

This section considers the two-dimensional flow in the small gap between two eccentric rotating cylinders.
This example yields nontrivial fluid dynamics and includes immersed interfaces that are nearly in contact.
Fig.~\ref{fig:eccentric_cylinders_schematic} provides a schematic of the problem geometry and dimensions.
The parameter set including the dimensions of the problem are the same as Fai and Rycroft \cite{fai2018lubricated}.
The computational domain is taken to be $\Omega = [-1,1]^2$, and $R_1=3/4$ and $R_2=3/4\times(1+1/24)$ are the radii of inner and outer cylinders, respectively.
The inner cylinder is centered at the origin, and the outer cylinder is centered at $(e,0)=(3/128,0)$.
The outer cylinder is set to be stationary while the inner cylinder rotates angular velocity $\omega_1=8.33\times 10^{-4}$.
Asymptotic analytical solutions can be derived for the velocity and pressure of the flow within long eccentric cylinders using lubrication theory \cite{christensen1971hydrodynamic},
\begin{align}
u(x,y) &=\begin{cases}
-\omega_{1}y, & \sqrt{x^2 + y^2} < R_{1},\\
-\omega_{1}y\big(1 - \gamma - \frac{3\epsilon(\gamma-{\gamma}^2)(2x+3\epsilon\sqrt{x^2 + y^2}+\epsilon^2 x)}{(2+\epsilon^2)(1+\epsilon x)}\big), & \sqrt{x^2 + y^2} \ge R_{1} \, \text{and} \, \sqrt{(x-e)^2 + y^2} \leq R_{2},\\
0, & \sqrt{(x-e)^2 + y^2}> R_{2},
\end{cases} \\
v(x,y) &=\begin{cases}
\omega_{1}x, & \sqrt{x^2 + y^2} < R_{1},\\
\omega_{1}x\big( 1 - \gamma -  \frac{3 \epsilon (\gamma-{\gamma}^2)(2x+3\epsilon+\epsilon^2 x)}{(2+\epsilon^2)(\sqrt{x^2 + y^2} +\epsilon x)}\big), & \sqrt{x^2 + y^2} \ge R_{1} \, \text{and} \, \sqrt{(x-e)^2 + y^2} \leq R_{2},\\
0, & \sqrt{(x-e)^2 + y^2}> R_{2},
\end{cases}\\
\end{align}
and \\
 \begin{align}
p(x,y) &=\begin{cases}
\frac{\omega^2_{1}(x^2 + y^2)}{2} + p_0, & \sqrt{x^2 + y^2} < R_{1},\\
\big(\frac{6\epsilon\mu R_1}{c^2}\big)\big(\frac{2y(x^2 + y^2)+\epsilon xy}{(2+\epsilon^2)(x^2 + y^2 +\epsilon x)}\big) & \sqrt{x^2 + y^2} \ge R_{1} \, \text{and} \, \sqrt{(x-e)^2 + y^2} \leq R_{2},\\
0, & \sqrt{(x-e)^2 + y^2}> R_{2},\\ 
\end{cases}
\end{align}
in which the nondimensional thickness is $\epsilon =(R2-R1)/R1$, the radial clearance is $c =R_{2}-R_{1}$, and $p_0$ is an arbitrary constant.
The reduced radial coordinate parameter is $\gamma =(x^2 + y^2 - R_{1}\sqrt{x^2 + y^2})/(R_2 - R_1 + e x)$.
A combination of zero normal-traction and zero tangential velocity boundary conditions is imposed on $\p\Omega$.

\begin{figure}[t!!]
		\centering
			\includegraphics[width=0.95\textwidth]{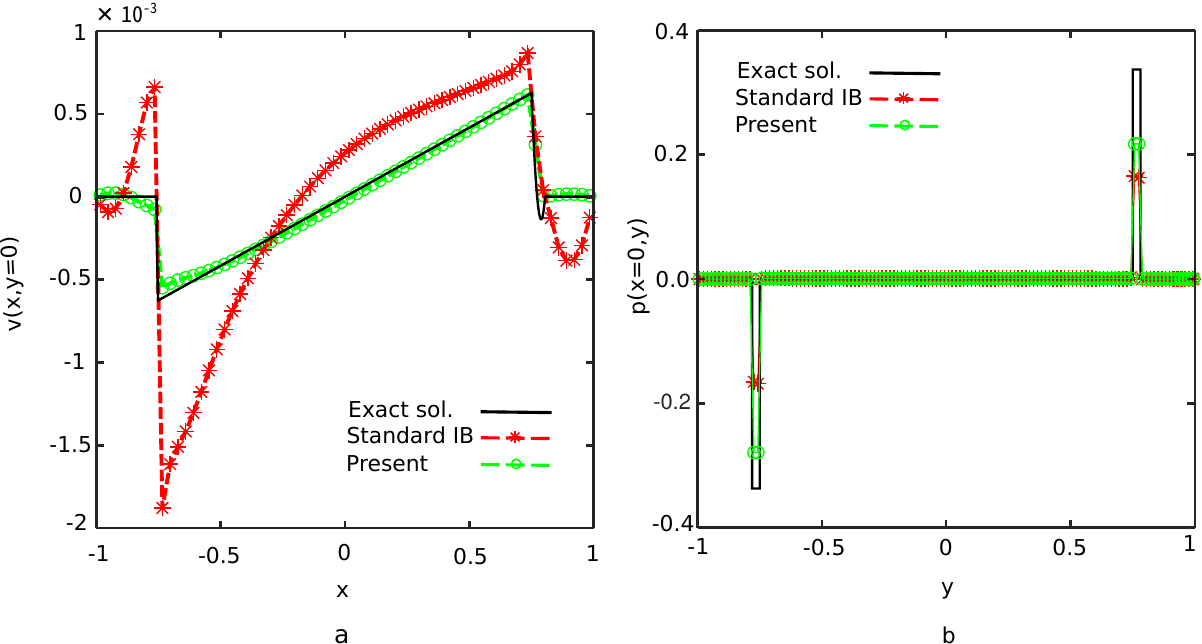}
		
\caption{(a) $v$-velocity and (b) pressure profiles for the IIM and IB solutions of the two-dimensional flow within the eccentric rotating cylinders on a relatively coarse $128\times 128$ Cartesian grid, for which the minimum distance between the cylinders is $h/2$.
The present IIM yields results that are in reasonable quantitative agreement with the asymptotic solution in this underresolved case.}
\label{fig:u-p-128-profile-eccentric-cylinders-2d} 
\end{figure}

\begin{figure}[t!!]
		\centering
			\includegraphics[width=0.95\textwidth]{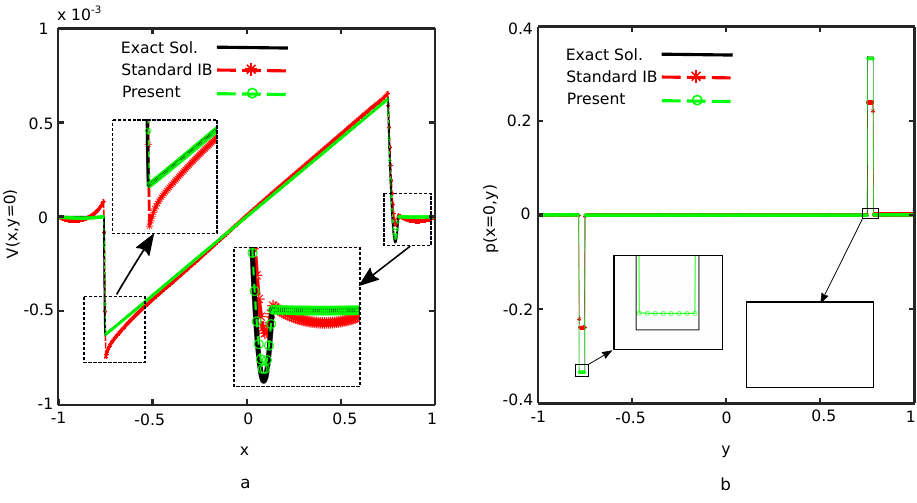}
		
\caption{(a) $v$-velocity and (b) pressure profiles for the IIM and IB solutions of the two-dimensional flow within eccentric rotating cylinders on a relatively fine $512\times 512$ Cartesian grid, for which the minimum distance between the cylinders is $2h$.
The present IIM yields results that are in excellent quantitative agreement with the asymptotic solution even in this marginally resolved case.}
\label{fig:u-p-512-profile-eccentric-cylinders-2d} 
\end{figure}

The steady state velocity and pressure solutions of the standard IB method and the present approach are obtained using two different mesh sizes.
Velocity and pressure profiles are respectively visualized along sections A and B shown in Fig.~\ref{fig:eccentric_cylinders_schematic}.
As in the previous example, a piecewise linear regularized delta function is used to spread the force in the standard IB method.
The penalty parameters are chosen to be $\kappa=2.0 \times 10^{-4}/\Delta t^2$ and $\eta = 0$.
In the first case, a relatively coarse $128 \times 128$ Cartesian grid is used, resulting in the two cylinders coming as close as $h/2$. 
Fig.~\ref{fig:u-p-128-profile-eccentric-cylinders-2d} compares the computed results to the asymptotic solution using this coarse Cartesian grid resoltion.
Although the gap is clearly underresolved in this case, the IIM gives much better accuracy in the velocity field than the standard IB method, which yields large overshoots near the boundaries.
We also use a relatively fine $512 \times 512$ Cartesian grid  is used, resulting in a minimal gap of size $2h$ between the cylinders.
Fig.~\ref{fig:u-p-512-profile-eccentric-cylinders-2d} shows that the IIM yields accurate approximations of both the velocity and pressure in this case.

\subsection{Flow past a stationary cylinder}
\label{subsec:cylinder}

Flow past a stationary cylinder is a widely used benchmark problem for testing numerical methods involving immersed boundaries.
Here, the immersed boundary is a disc centered at the origin with diameter $D=1$. 
The physical domain is $\Omega = [-15,45]\times[-30,30]$, a square of length $L=60$. 
This domain configuration corresponds to the problem setup in Taira and Colonius \cite{taira2007immersed}. 
A uniform inflow velocity $U=(1,0)$ is imposed on the left boundary ($x=-15$), 
and zero normal traction and tangential velocity is imposed at the right boundary ($x=45$) as an outflow condition.
Along the bottom ($x=-30$) and top ($x=30$) boundaries, the normal velocity and tangential traction are set to zero.

We set $\rho=1$ and use the inflow velocity $U$ as the characteristic velocity.
The Reynolds number is $\Re=\frac{\rho U D}{\mu}$.
Reynolds numbers between 20 to 200 are considered.
The computational domain is discretized using a locally refined grid with $N=6$ nested grid levels and a refinement ratio of two between levels.
The Cartesian grid spacing on the coarsest level is $h_{\text{coarsest}}=\frac{L}{64}$, and $h_{\text{finest}}=2^{-(N-1)}h_{\text{coarsest}}$ is the grid spacing on the finest grid level.
The time step size is $\Delta t=0.05h_{\text{finest}}$, yielding a maximum advective CFL number in the range 0.1--0.2.
Values of $\kappa$ and $\eta$ are tuned to ensure both rigidity of the structure and stability of the dynamics. 

\begin{table}[b!!]
	\centering	
	\caption{Comparison of computational and experimental values for the flow past a stationary cylinder at $\Re=20$ and $\Re=40$. Simulation parameters include $h_\text{finest}=0.0293$, $\Delta t=0.05h_\text{finest}$, and $\mfac=2$.}
	\label{table:Cylinder_drag_coeffs_20_40}
\begin{tabular}{l*{6}{c}r}
             & &$\Re = 20$& & & $\Re = 40$ & &\\
\cmidrule(lr){2-4}\cmidrule(lr){5-7}
 & $\Lwake$ & $\thetas$ & $\CD$ & $\Lwake$ & $\thetas$ & $\CD$  \\
\hline
Tritton (experimental) \cite{tritton1959}            & - & - & 2.22 & - &  - & 1.48 \\
Le et al. \cite{le2006immersed}          & 0.93 & 43.9 & 2.05 & 2.35 &  53.8 & 1.52 \\
Xu and Wang \cite{xu2006systematic}      & 0.92 & 44.2 & 2.23 & 2.21 &  53.5 & 1.66  \\
Calhoun \cite{calhoun2002}      & 0.91 & 45.5 & 2.19 & 2.18 & 54.2 & 1.62  \\
Present      & 0.93 & 44.4 & 2.10 & 2.31 &  54.1 & 1.58  \\
\end{tabular}
\end{table}

To quantitatively assess the computed dynamics, we compute nondimensional quantities including the drag coefficient $CD$ and
lift coefficient $CL$ as,
\begin{equation}
(\CD,\CL)=\frac{-\int_{\Gamma_0}\F(\s,t)\ \mathrm{d}A}{\frac{1}{2}\rho U^{2}D},
\end{equation}
and the Strouhal number,
\begin{equation}
\St=\frac{f_\text{s}D}{U},
\end{equation}
in which $F^{x}$ and $F^{y}$ are the $x$ and $y$ components of the penalty force, and $f_\text{s}$ is the vortex shedding frequency. 
Table~\ref{table:Cylinder_drag_coeffs_20_40} lists the drag coefficient ($\CD$), recirculation length ($L_{\text{wake}}$), and angle of separation (${\theta}_\text{s}$) for $\Re=20$ and $\Re=40$.
For this range of Reynolds numbers, the flow separates from the back of the cylinder, and a pair of vortices form that gradually approach a steady state.
Comparisons of the results obtained by our method with previous numerical and experimental results show excellent quantitative agreement in all flow characteristics. 

\begin{figure}[t!!]
		\centering
			\includegraphics[width=0.9\textwidth]{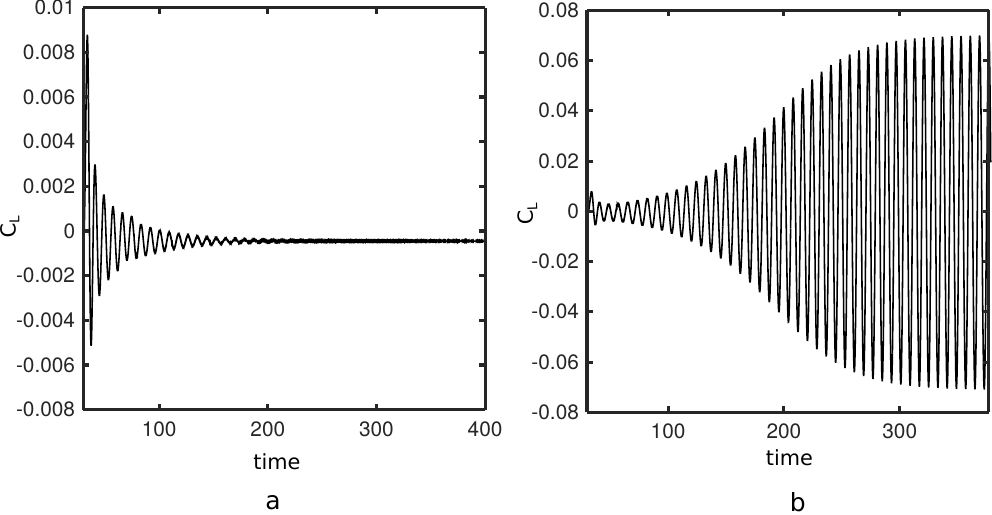}
		
		\caption{ The lift coefficient ($\CL$) of flow past a stationary cylinder over time for (a) $\Re=40$ and (b) $\Re=50$. The onset of the von Karman vortex street is detected at $\Re=50$. 
		Simulation parameters include $h_\text{finest}=0.0293$, $\Delta t=0.05h_\text{finest}$, and $\mfac=2$.}
		\label{fig:Cylinder_Re_40_50}		
\end{figure}

\begin{figure}[t!!]
		\centering
			\includegraphics[width=0.95\textwidth]{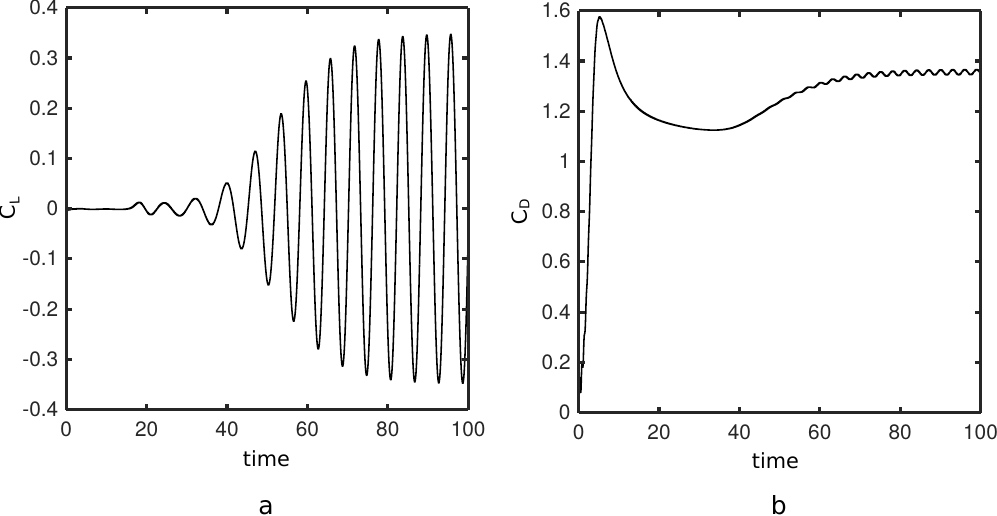}
		
		\caption{Lift ($\CL$) and drag ($\CD$) coefficients over time for flow past a cylinder at $\Re=100$.
		Simulation parameters include $h_\text{finest}=0.0293$, $\Delta t=0.05h_\text{finest}$, and $\mfac=2$.
		}	
		\label{fig:CD_CL_100}		
\end{figure}

A transition from steady flow to alternate vortex shedding occurs for $40<\Re<50$.
The present IIM recovers this transition, as demonstrated in Fig.~\ref{fig:Cylinder_Re_40_50}.
Although an initial instability decays over time at $\Re=40$, the same instability at $\Re=50$ eventually leads to the well-known von Karman vortex street. 
Vortex shedding continues for $\Re=100$ and $\Re=200$, and the flow becomes increasingly unsteady.
Fig.~\ref{fig:CD_CL_100} details the lift and drag coefficients at $\Re=100$.

\begin{figure}[t!!]
		\centering
			\includegraphics[width=0.95\textwidth]{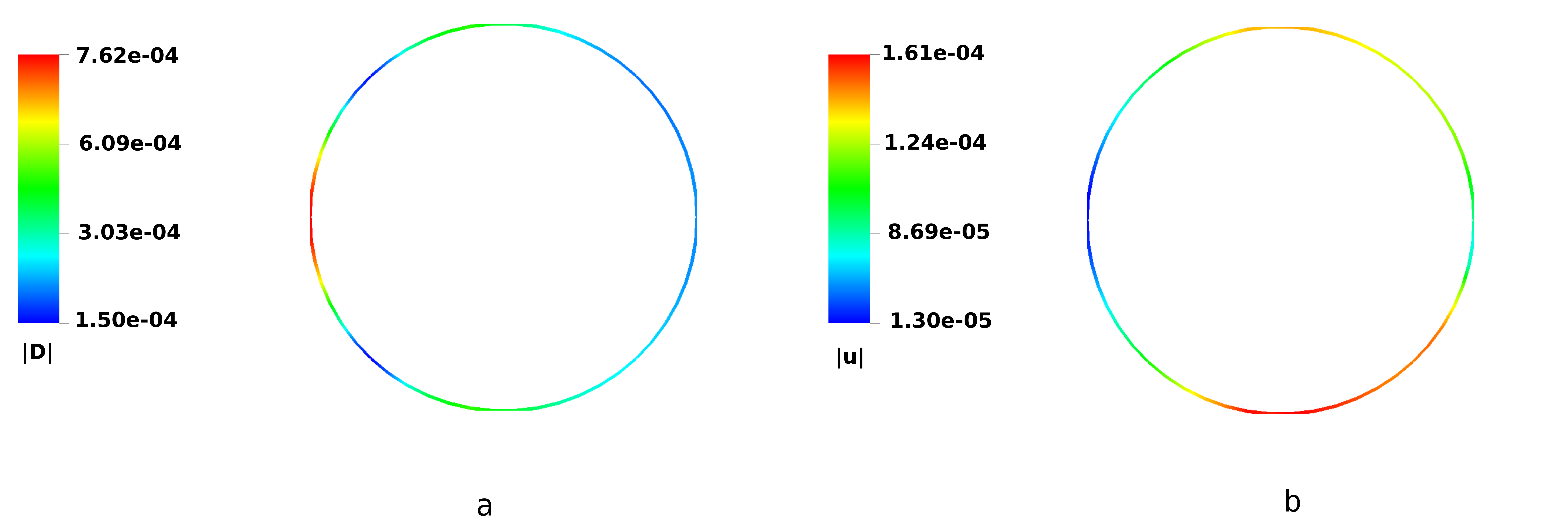}
		
		\caption{Interfacial (a) displacement and (b) velocity magnitudes for flow around a stationary cylinder at $\Re=200$ using the present IIM,
		          plotted at the same time as shown in Fig.~\ref{fig:FPC_zoom_IIM_IB_200_vorticity}.
		          Simulation parameters include $h_\text{finest}=0.0293$, $\Delta t=0.05h_\text{finest}$, and $\mfac=2$.
		          The maximum displacement error is approximately $2.6\%$ that of the background grid spacing $h_\text{finest}$, and the velocity error is approximately $0.016\%$ that of the free stream flow velocity.}
		\label{fig:FPC_IIM_200}
\end{figure}

Fig.~\ref{fig:FPC_IIM_200} shows representative interfacial displacement and velocity magnitudes at $\Re=200$.
The maximum displacement magnitude is approximately $2.6\%$ that of the background grid spacing ($h_\text{finest}=\frac{60}{64\times2^5}=0.0293$), and the velocity magnitude is approximately $0.016\%$ that of the free stream flow velocity ($U = 1$).

\begin{figure}[t!!]
		\centering
			\includegraphics[width=0.95\textwidth]{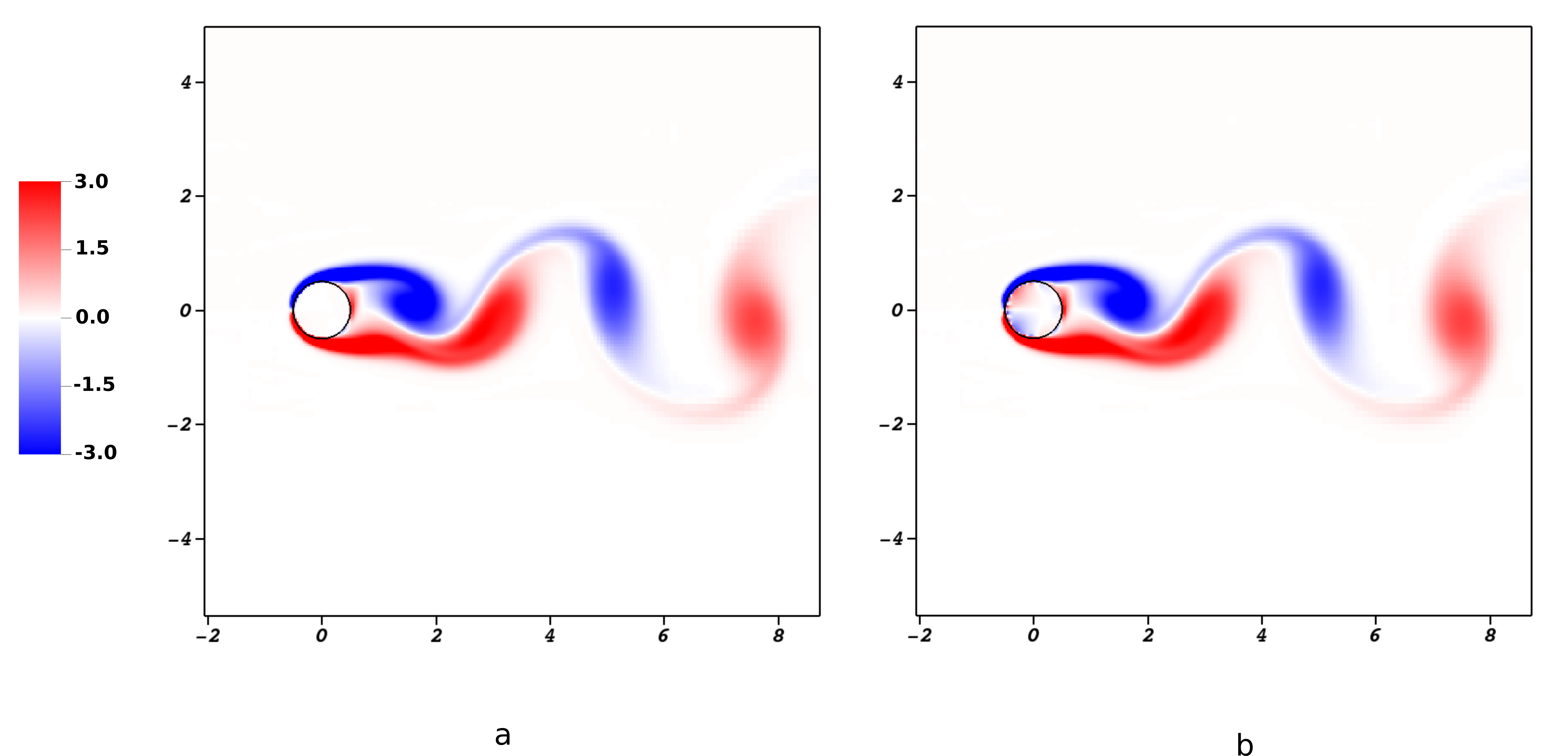}
		
		\caption{Vorticity fields for flow past a cylinder at $\Re = 200$ for (a) the present IIM and (b) the conventional IB method.
		Notice that the conventional IB method generates spurious interior flows that are eliminated with the present IIM.}	
\label{fig:FPC_zoom_IIM_IB_200_vorticity} 
\end{figure}

Fig.~\ref{fig:FPC_zoom_IIM_IB_200_vorticity} compares the vorticity fields generated by the immersed interface and IB methods at $\Re=200$.
It is clear that the IIM gives a very sharply resolved vorticity profile while allowing for essentially no flow in the interior of the structure.
This is in clear contrast to the conventional IB method, which generates spurious flow regions inside the cylinder, as seen in Fig.~\ref{fig:FPC_zoom_IIM_IB_200_vorticity}(b).

\begin{table}[t!!]
	\centering	
	\caption{Comparison of computational and experimental values of the drag coefficient ($\CD$), lift coefficient ($\CL$), and Strouhal numbers (St) for flow past a cylinder at $\Re=100$ and $\Re=200$.
		          Other simulation parameters include $h_\text{finest}=0.0293$, $\Delta t=0.05h_\text{finest}$, and $\mfac=2$. All values are computational unless otherwise noted.
	}
	\label{table:Cylinder_drag_coeffs_100_200}	
\begin{tabular}{l*{6}{c}r}
             & &$\Re = 100$& & & $\Re = 200$ & &\\
\cmidrule(lr){2-4}\cmidrule(lr){5-8}
 & $\CD$ & $\CL$ & $\text{St}$ & $\CD$ & $\CL$ & $\text{St}$  \\
\hline
Braza et al.~\cite{braza1986}			&  $1.360 \pm 0.015$ & $\pm 0.250$ & 0.160 & $1.400 \pm 0.050$ &  $\pm 0.75$ & 0.200   \\
Liu et al.~\cite{liu1998}     & $1.350 \pm 0.012 $ & $\pm 0.339$ & 0.164 & $1.310 \pm 0.049$ &  $\pm 0.69$ & 0.192   \\
Calhoun \cite{calhoun2002}          & $1.330 \pm 0.014$ & $\pm 0.298$ & 0.175 & $1.170 \pm 0.058$ &  $\pm 0.67$ & 0.202   \\
Le et al.~\cite{le2006immersed}        & $1.370 \pm 0.009$ &  $\pm 0.323$ & 0.160 & $1.340 \pm 0.030$ &  $\pm 0.43$ & 0.187   \\
Xu and Wang \cite{xu2006systematic}        & $1.423 \pm 0.013 $& $\pm 0.340$ & 0.171 & $1.420 \pm 0.040$  & $\pm 0.66$  & 0.202  \\
Lai and Peskin \cite{lai2000}           & - & - & 0.165 & - &  - & 0.190   \\
Roshko (experimental) \cite{roshko1961}        & - & - & 0.164 & - &  - & 0.190   \\
Williamson (experimental) \cite{williamson1996vortex}          & - & - & 0.166 & - &  - & 0.197   \\
Griffith and Luo \cite{BEGriffith17-ibfe}      & - & - & - & $1.360 \pm 0.046$ &  $\pm 0.70$ & 0.195  \\
Present      & $1.370 \pm 0.015$ & $\pm 0.351$ & $0.168$ & $1.390 \pm 0.060$  &  $\pm 0.75$ & 0.198  \\
\end{tabular}
\end{table}

Table \ref{table:Cylinder_drag_coeffs_100_200} compares $\CL$, $\CD$, and $\text{St}$ with previous experimental and computational studies at $\Re=100$ and $\Re=200$.
The force coefficients fall within the range of values reported in previous studies.

\begin{figure}[ht!]
		\centering
			\includegraphics[width=0.95\textwidth]{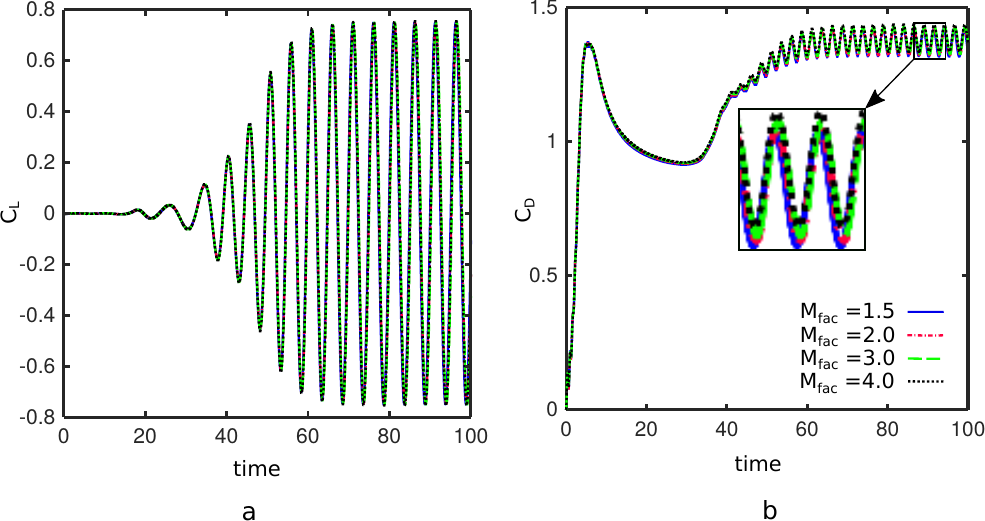}
		
		\caption{Lift ($\CL$) and drag ($\CD$) coefficients for flow past a cylinder at $\Re=200$. The method produces consistent results across a wide range of $\mfac$ values.}	
		\label{fig:CDL_P1_IIM_200}		
\end{figure}

We also investigate the effects of varying the relative mesh spacing parameter $\Mfac$.
A range of relative mesh spacing values are considered.
Fig.~\ref{fig:CDL_P1_IIM_200} shows the lift and drag coefficients at $\Re=200$ with $\Mfac=1.5$, $2$, $3$, and $4$.
The results are consistent across a relatively large range of $\Mfac$ values.
This mesh insensitivity is important when modeling complex geometries, for which large variations in the element sizes are likely.
We remark that it is not possible to use arbitrarily small values of $\Mfac$ with the present method.
To see why this is, recall that forces are transmitted from the interface to the background Cartesian grid through jump conditions that are determined by intersecting the finite difference stencils against the interface mesh.
Forces that are associated with elements that are not intersected by any finite difference stencil will have no physical effect on the flow dynamics.
This is physically unstable, and if it occurs in a simulation, the computation also will generally become unstable.
Consequently, we consider values of $\Mfac$ that are large enough to ensure that all surface elements are pierced by at least one finite difference stencil.

\subsection{Flow past a spinning cylinder}
\label{subsec:rotating_cylinder}

\begin{figure}[t!!]
		\centering
			\includegraphics[width=0.9\textwidth]{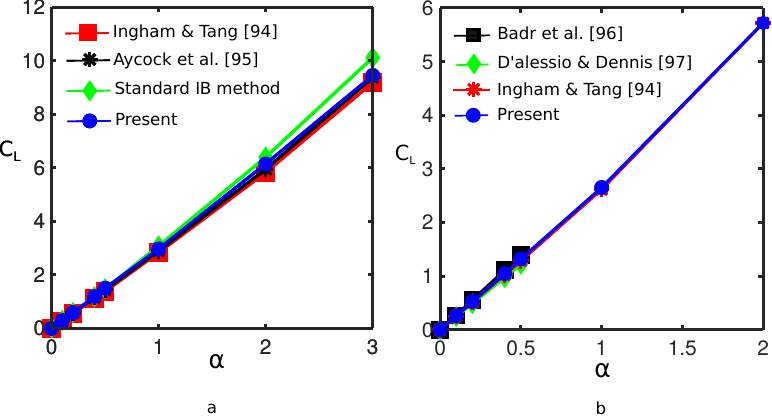}
		
		\caption{Lift coefficient ($\CL$) versus nondimensional angular velocity $\alpha$ for the flow past a rotating cylinder at (a) $\Re = 5$  and (b) $\Re = 20$.
		Both cases use $\mfac= 2$. Results produced by the present method are in excellent agreement with previous numerical studies.}	
		\label{fig:spinning_cylinder_comparison}	
\end{figure}

We also consider the performance of the method in determining the lift coefficient for flow over a spinning cylinder.
The computational domain and the Eulerian boundary conditions are the same as the case of the stationary cylinder in Sec.~\ref{subsec:cylinder}. 
The cylinder is prescribed to rotate about its central axis at nondimensional rotation rate $\alpha = \frac{\omega R}{U}$, in which $\omega$ is the angular velocity.
Fig.~\ref{fig:spinning_cylinder_comparison} compares the lift coefficient generated by the present method to previous work \cite{ingham1990numerical, aycock2017resolved,BEGriffith17-ibfe,badr1989steady, d1994vorticity} for $\alpha$ in the range 0--3 at $\Re = 5$ and $20$. 
The $\CL$ values generated by the present IIM are in excellent agreement with prior work for the entire range of $\alpha$ values considered.

\begin{figure}[t!!]
		\centering
			\includegraphics[width=0.9\textwidth]{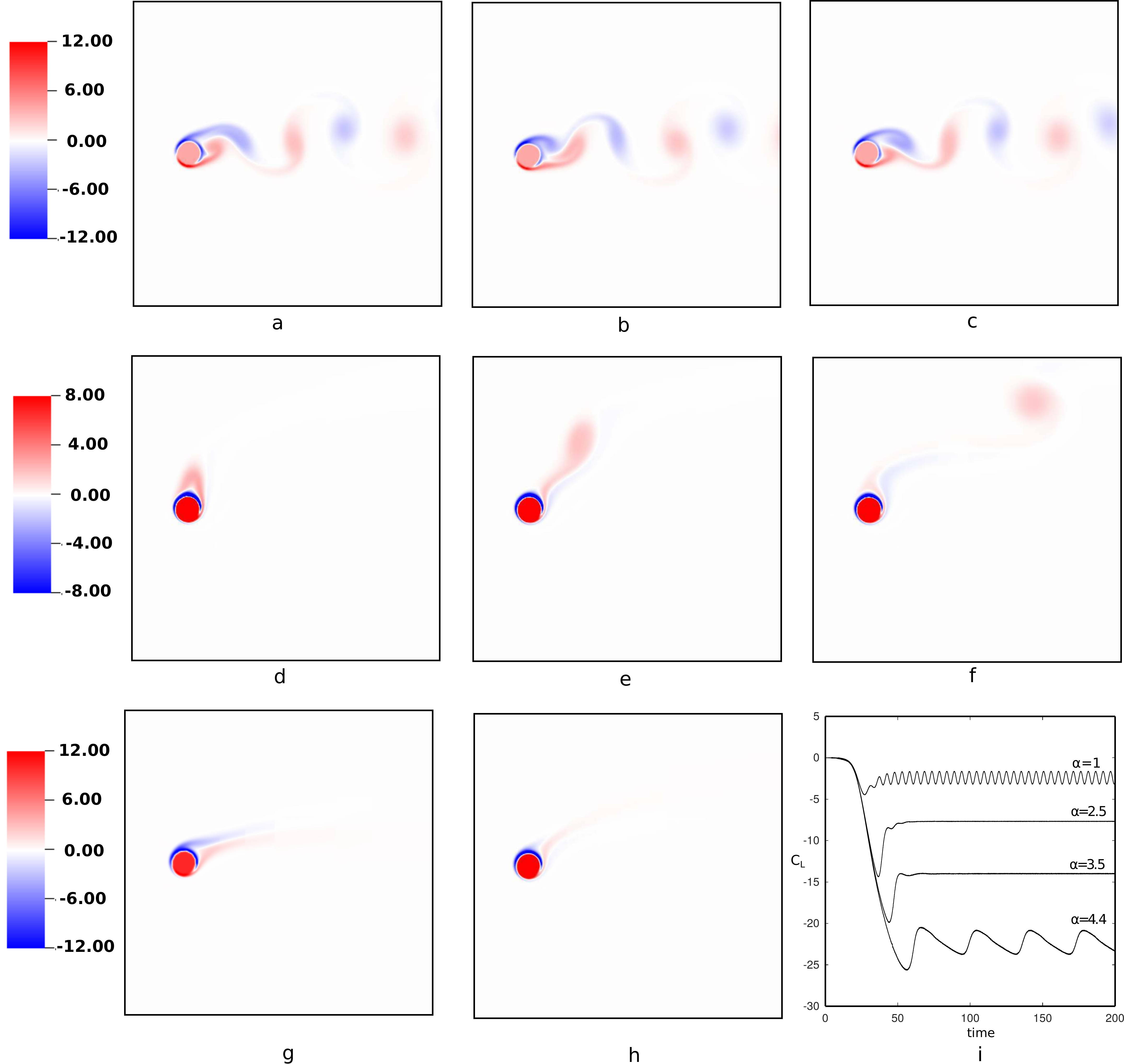}
		
		\caption{ Vorticity fields of flow past a rotating cylinder for four rotation rates at $\Re = 200$. Panels (a), (b), and (c) show
		vorticity fields corresponding to $\alpha=1.0$ at $t=128$, $130$, and $132$, respectively. 
		Vorticity contours of $\alpha=4.4$ are shown in panels (d), (e), and (f) at $t=132$, $138$, and $144$.
		Periodic flows with distinct behaviors are observed for $\alpha=1.0$ and $4.4$,
		  whereas the solutions at $\alpha=2.5$ shown in panel (g) and at $\alpha=3.5$ in panel (h) reach steady state.
		  Subfigure (i) shows time histories of the corresponding lift coefficient $\CL$ in the lower bottom panel.}    
\label{fig:flow_spinning_cylinder_vorticity} 
\end{figure}

Fig.~\ref{fig:flow_spinning_cylinder_vorticity} shows the vorticity fields for $\alpha=1$, $2$, $2.5$, $3.5$, and $4.4$ at $\Re = 200$.  
As previously demonstrated by Mittal and Kumar \cite{mittal2003flow}, the time history of $\CL$ reveals vortex shedding behavior for $\alpha < 1.9$.
At higher rotation rates, the flow is expected to achieve a steady state except for $4.34 < \alpha < 4.70$, for which the flow undergoes a different kind of instability, such that only one-sided vortex shedding occurs.
Fig.~\ref{fig:flow_spinning_cylinder_vorticity} demonstrates that these flow conditions are accurately predicted using the present IIM.

\subsection{Flow past a stationary sphere}
\label{subsec:sphere}

This section examines flow past a stationary sphere at low to moderate Reynolds numbers.
We consider a sphere with diameter of $D=1$ centered at the origin and placed in the computational domain $\Omega = [-15,45]\times[-30,30]\times[-30,30]$. 
A uniform inflow velocity $U=(1,0,0)$ is imposed at $x=-15$, and zero normal traction and tangential velocity is imposed at $x = 45$ as an outflow condition.
Along the other boundaries, the normal velocity and tangential traction are set to zero.


%

\begin{table}[t!]
	\centering	
	\caption{Drag coefficients for three-dimensional flow past a sphere at different Reynolds numbers determined by the present IIM and by several prior
	 computational \cite{fornberg1988steady,xu20083d,fadlun2000combined} and empirical \cite{turton1986short} studies.
	 	 Simulation parameters include $h_\text{finest}=0.0293$, $\Delta t=0.05h_\text{finest}$, and $\mfac=2$.
}
	\label{table:Sphere_drag_coeff}	
\begin{tabular}{l*{5}{c}r}
              & $\Re=20$ & $\Re=100$ & $\Re=200$ & $\Re=500$  \\
\hline
Fornberg \cite{fornberg1988steady} & - & 1.0852 & 0.7683 & 0.4818  \\
Turton and Levenspiel \cite{turton1986short} & 2.6866 & 1.0994 & 0.8025 & 0.5617   \\
Fadlun et al. \cite{fadlun2000combined}  & - & 1.0794 & 0.7567 & 0.4758   \\
Campregher et al. \cite{campregher2009} & - & 1.1781 & 0.8150 & 0.5200   \\ 
Xu and Wang \cite{xu20083d}           & 2.73 & 1.15 & 0.88 & -   \\
Present       & 2.6940 & 1.0920 & 0.7852 & 0.5348   \\
\end{tabular}
\end{table}

We set $\rho=1$ and use the inflow velocity $U$ as the characteristic velocity.
The Reynolds number is $\Re=\frac{\rho U D}{\mu}$, and we consider Reynolds numbers from 20 to 500.
The computational domain is discretized using an $N$-level locally refined grid with a refinement ratio of two between the grid levels.
The Cartesian grid spacing on the coarsest level is $h_{\text{coarsest}}=\frac{L}{64}$, and $h_{\text{finest}}=2^{-(N-1)}h_{\text{coarsest}}$ is the grid spacing on the finest grid level.
The time step size is $\Delta t=0.05h_{\text{finest}}$, yielding a maximum CFL number of approximately 0.1--0.2.
The surface of the sphere is described using bilinear quadrilateral surface elements with $\Mfac=2$.
The spring stiffnesses and damping parameter at $\Re=20$ are taken to be $\kappa=1000$ and $\eta=5$.
The values for $\Re=100$ are $\kappa=300$ and $\eta=0.5$.
At higher Reynolds numbers of $\Re=200$ and $\Re=500$, only spring forces are used, and the spring stiffness in these cases are chosen to $\kappa=250$ and $\kappa=200$, respectively.
We compute the force coefficients as
\begin{equation}
    (\CD,C_\text{L}^y,C_\text{L}^z) = \frac{-2\int_{\Gamma_0}\F \ \mathrm{d}A}{A_{\text{proj}}},
\end{equation}
in which $A_{\text{proj}}=\pi/4$ is the projected area of the sphere with $D=1$.
Table \ref{table:Sphere_drag_coeff} compares results generated by the present method to previous computational and experimental studies \cite{fornberg1988steady,turton1986short,fadlun2000combined,campregher2009,xu20083d} 
for $\Re$ between 20 to 500.
We observe excellent agreement between the values produced by the present method and these previous works.

\begin{figure}[t!!]
		\centering
			\includegraphics[width=0.8\textwidth]{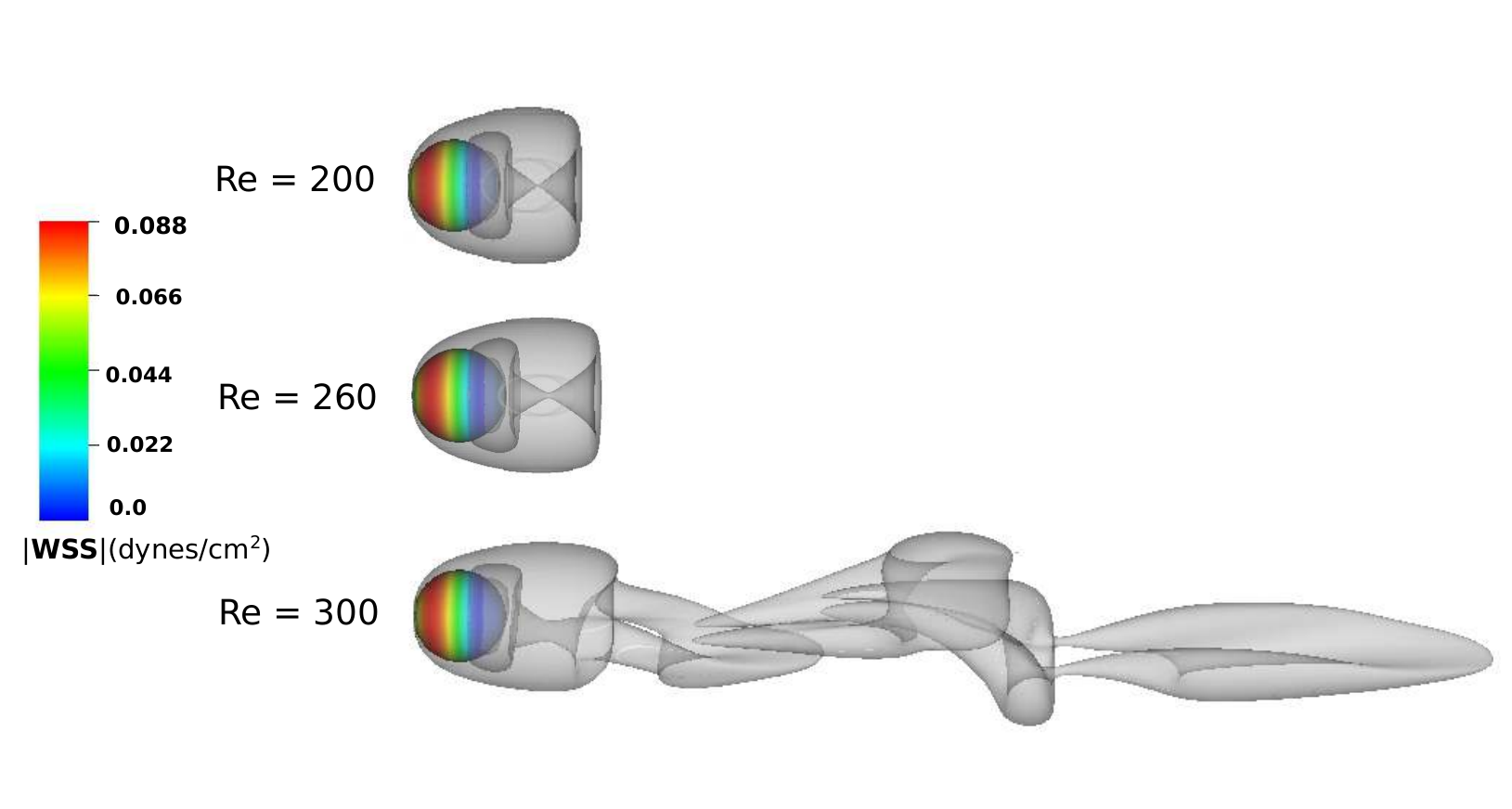}
		
\caption{Iso-surface of Q-criterion \cite{hunt1988eddies} showing the vortex dynamics for flow past a sphere at $\Re = 200$, $260$, and $300$.
The color contours on
the sphere show the distribution of the wall shear stress, with red and blue being regions of largest and lowest magnitudes, respectively.
 The onset of the unsteady vortex shedding regime occurs in the range  $290<\Re<400$. The unsteady behavior at $\Re=300$ predicted by the present method agrees with previous studies \cite{sakamoto1990study,johnson1999flow}.}
\label{fig:Sphere-Re-200-300} 
\end{figure}

To test the ability of the method in predicting the dynamics of the flow around the sphere in the transitional regime, simulations are performed at $\Re=200$, $260$, and $300$.
It is well established that the onset of the unsteady vortex shedding regime occurs for $290 < \Re < 400$ \cite{sakamoto1990study,johnson1999flow}.
Fig.~\ref{fig:Sphere-Re-200-300} shows that unsteady vortex shedding occurs using the present approach at $\Re=300$.
Fig.~\ref{fig:Sphere-Re-200-300} also shows the magnitude of the wall shear stress on the surface of the sphere.  

		


\subsection{Flow inside the inferior vena cava}
\label{subsec:ivc}

\begin{figure}[t!!]
		\centering
			\includegraphics[width=0.8\textwidth]{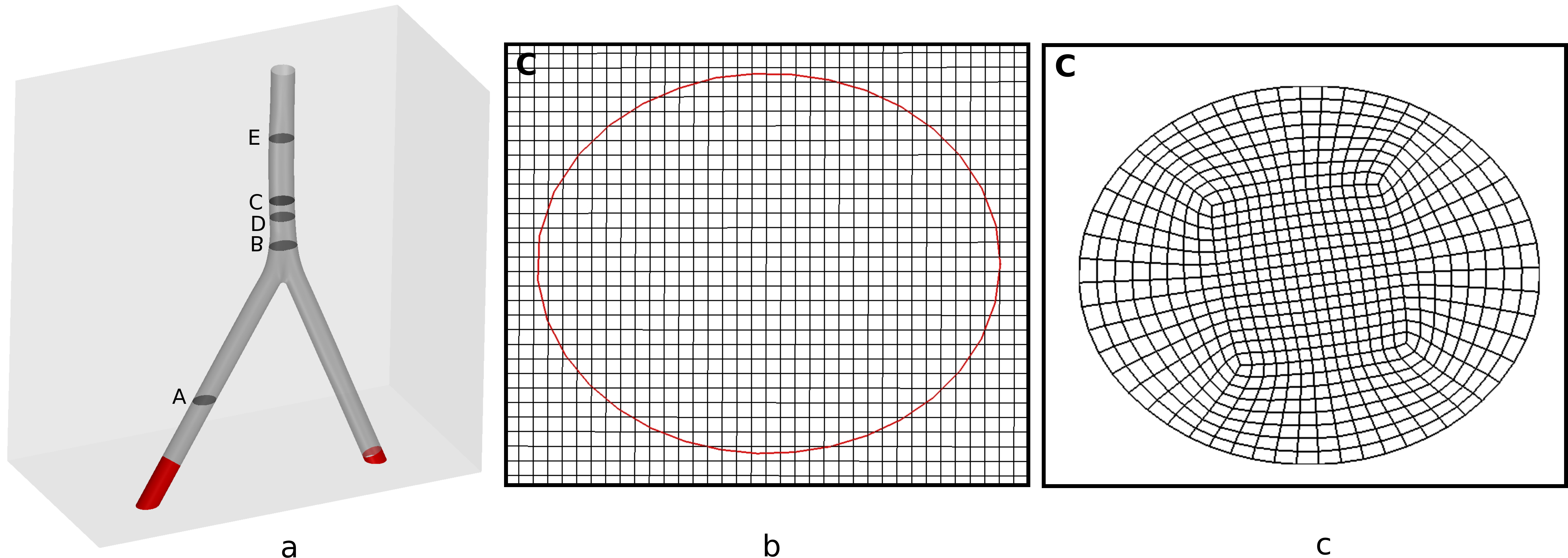}
		
		\caption{(a)~The patient-averaged inferior vena cava (IVC) model geometry used to compare the performance of the present IIM scheme to a body-fitted finite volume simulation performed with OpenFOAM.  Comparisons between the two approaches focus on the highlighted sections A--E.  (b) The projected Eulerian and interface geometry used in the IIM simulation along section C.  (c) The body fitted mesh used in the OpenFOAM simulation along section C.  The number of grid cells along the major and minor axes are roughly equal for both simulations.}
\label{fig:Eulerian-mesh-ivc} 
\end{figure}

\begin{figure}[t!!]
		\centering
			\includegraphics[width=0.6\textwidth]{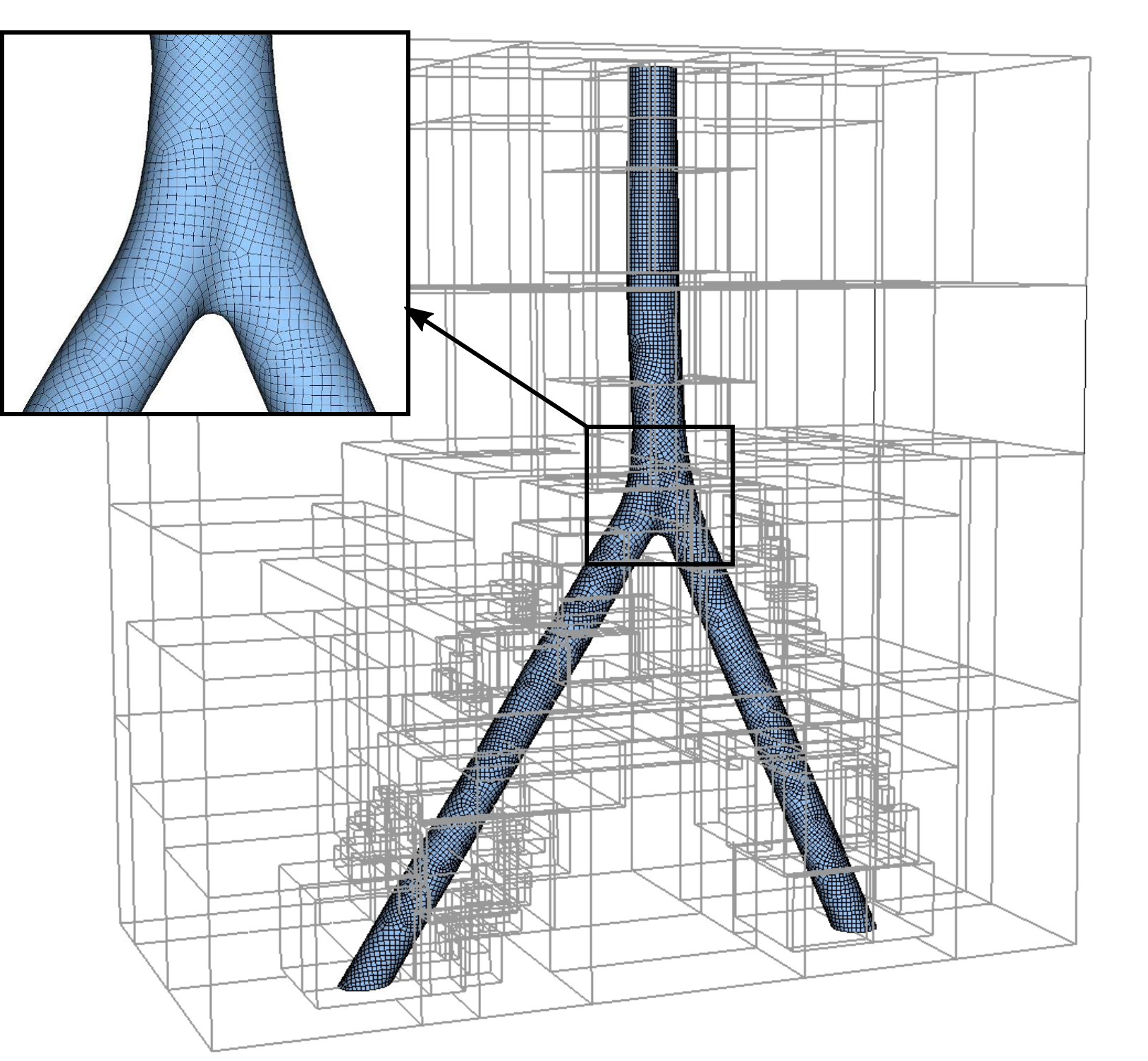}
		
		\caption{Computational mesh of the IVC including quadrilateral surface elements and a locally refined Cartesian grid.} 
\label{fig:Quad-mesh-ivc} 
\end{figure}

The inferior vena cava (IVC) is a large vein through which deoxygenated blood from the lower and middle body returns to the right atrium of the heart.
As a demonstration of the capability of the present method to handle complex geometries, we consider flow through a model of the human IVC.
This geometry used in this section is a modified version of the patient-averaged model by Rahbar et al.~\cite{rahbar2011three} that has been recently studied by Craven et al.~\cite{craven2018steady} and Gallagher et al.~\cite{gallagher2018steady}.
Fig.~\ref{fig:Eulerian-mesh-ivc}(a) shows the geometry of the IVC.
The morphological features of the IVC including the flow mixing at the junction of the iliac veins, infrarenal curvature, and non-circular vessel cross sections that all contribute to the complex hemodynamics and make this model a suitable demonstrative case to test the robustness of the present algorithm. 
The infrarenal IVC downstream of the iliac bifurcation has an average hydraulic diameter of $D_\text{h} =2.8 $ cm. 
A high flow rate of 100 cm$^{3}$/s is used, which corresponds to exercise flow conditions with a maximum Reynolds number of $\Re=1500$.
In accord with the experimental observations of Gallagher et al.~\cite{gallagher2018steady} in the same model and under the same flow conditions, here we assume that the flow is laminar.
The density of the fluid is set to $\rho=1.817 $ g/cm$^{3}$, and the viscosity is $\mu=5.487\times 10^{-2} $ g/cm s.
For comparison, an additional simulation is performed with the finite volume method using the SIMPLE (Semi-Implicit Method for Pressure-Linked Equations) solver in OpenFOAM (version 1812) with second-order accurate spatial discretization schemes. Importantly, the OpenFOAM simulation was performed using a body-fitted
 block-structured mesh that was generated to approximately match the spatial resolution of the IIM mesh (see Fig.~\ref{fig:Eulerian-mesh-ivc} panels b and c) to enable a consistent comparison.

Steady fully-developed parabolic velocity boundary conditions are imposed at the upstream inlets of the iliac veins.
The two inlets are circular with a diameter of $D=2.44 $ cm, which then transition to an elliptical shape a short distance downstream.
Because our current IIM implementation requires applying the Eulerian boundary conditions at the outer faces of the surrounding Cartesian domain, we have slightly extended the IVC geometry at the inlets to be able to apply the corresponding fully developed conditions at the resulting enclosed intersection ellipses such that the flow conditions are the same as in the OpenFOAM simulation; see Fig.~\ref{fig:Eulerian-mesh-ivc}(a).
The surface of the IVC is described using bilinear quadrilateral surface elements with $\Mfac \approx 2$, as shown in Fig.~\ref{fig:Quad-mesh-ivc}. 
The surface mesh of the IVC used in the IIM simulation was constructed using Trelis Mesh Generation Toolkit, which is based on the CUBIT software \cite{cubit}, and the mesh of the OpenFOAM simulation was generated using the Pointwise software.
The IVC is embedded in a rectangular computational domain of size $L_x\times L_y \times L_z = 50 \times 25 \times 50 \, \textrm{cm}^{3}$.
This Eulerian domain is discretized using a three-level locally refined grid with a refinement ratio of four between the grid levels, resulting in a grid spacing of $h_{\text{coarsest}}=\frac{25}{16}=1.5625$ cm on the coarsest level and $h_{\text{finest}}=\frac{25}{16\times 4^2} \approx 0.098$ cm on the finest grid level.
A fixed time step size of $\Delta t=10^{-4} \, \textrm{sec}$  is used.
At the outlet, the normal traction and tangential velocities are set to zero.
Solid-wall boundary conditions are imposed along the remainder of $\p\Omega$.
\begin{figure}[t!!]
		\centering
			\includegraphics[width=0.6\textwidth]{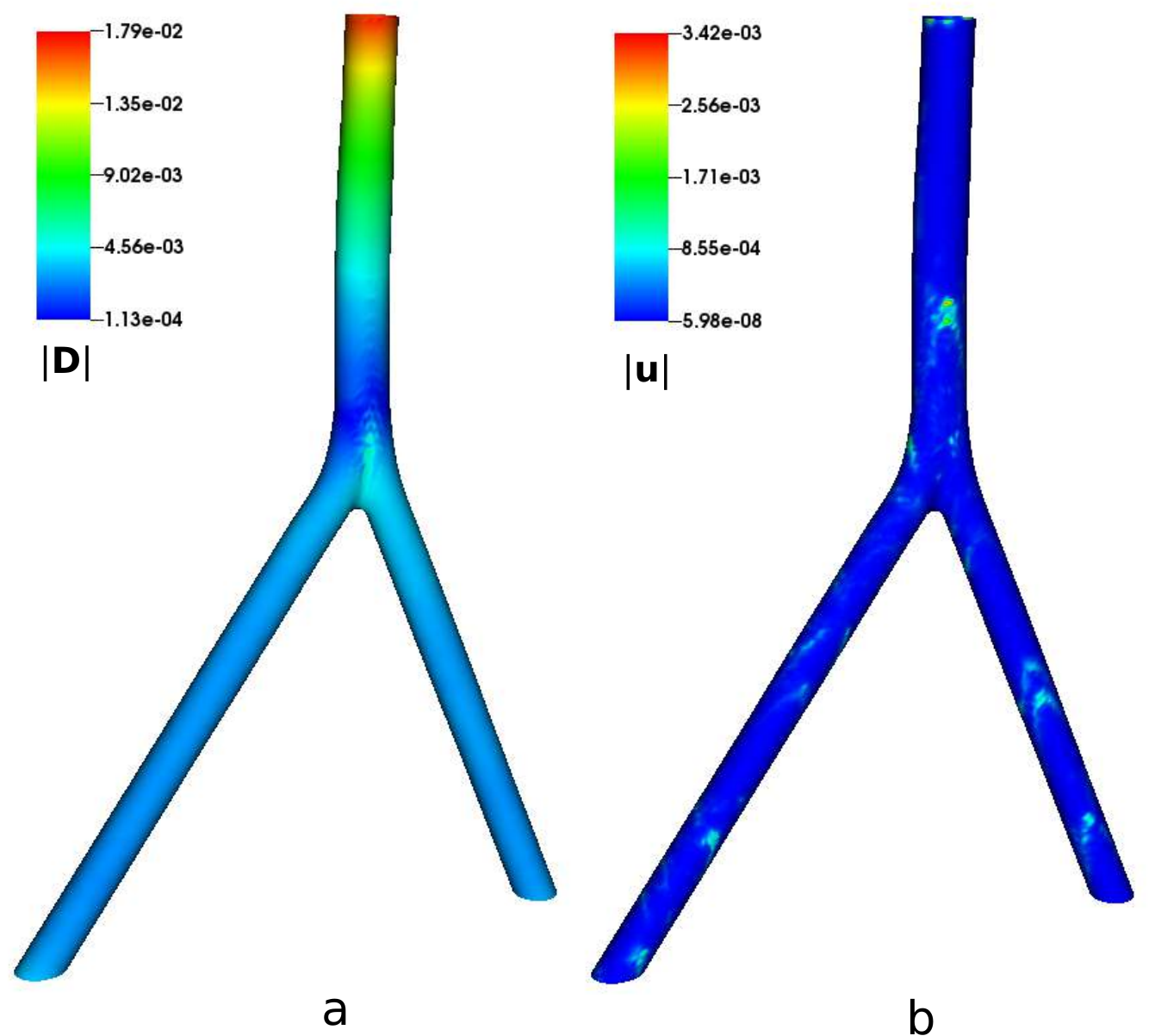}
		
		\caption{Lagrangian (a) displacement and (b) velocity magnitudes of the rigid IVC mesh using the present IIM. The maximum displacement is approximately $18\%$ that of the background 
		grid spacing, and the maximum velocity value is approximately $0.032\%$ that of the mean value velocity at the inlet ($\bar{U}\approx 10.7$ cm/s).}
\label{fig:disp-velocity-Lag-ivc} 
\end{figure}

\begin{figure}[t!!]
		\centering
			\includegraphics[width=0.9\textwidth]{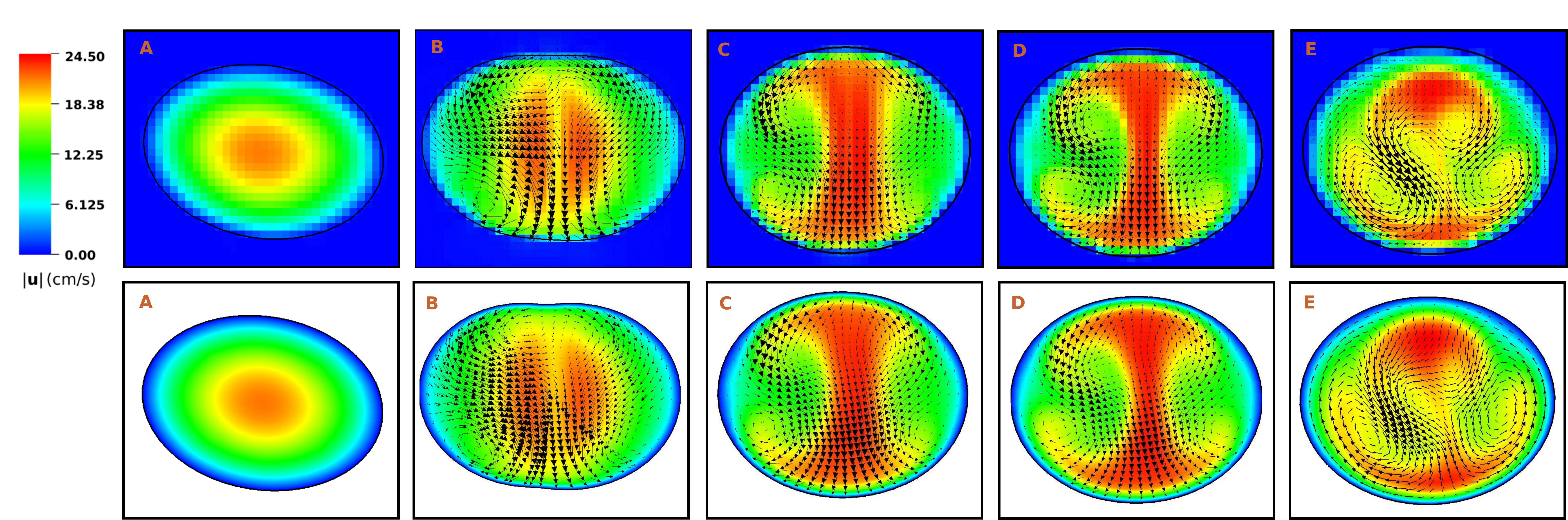}
		
\caption{Velocity magnitudes and vorticity vectors for flow in the inferior vena cava (IVC) at exercise flow conditions using the present IIM approach (top panels) and the OpenFOAM simulation (bottom panels).
 Good agreement between the two solutions is clearly obtained at all cross sections.
  Note that no additional post-processing smoothing or filtering is performed in the visualization of the results from the present method.}
 \label{fig:ivc_velocity} 
\end{figure}

Fig.~\ref{fig:disp-velocity-Lag-ivc} shows the Lagrangian displacement and velocity magnitudes of the stationary IVC using the present IIM. 
The maximum displacement of the mesh stays within $18\%$ of the Eulerian grid spacing, and the maximum velocity is approximately $0.032\%$ that of the mean value velocity at the inlet ($\bar{U}\approx 10.7 \, \textrm{cm/s}$).
Numerical experiments show that the fully-developed flows through the iliac veins remain almost completely undisturbed until they reach the confluence of the veins, where the two streams merge and form a high-velocity region in the center of the IVC lumen, with a pair of counter-rotating vortices on either side.
This is demonstrated in Fig.~\ref{fig:ivc_velocity}(a), where velocity contours of both the IIM computation and the OpenFOAM simulation at cross section A are shown.
Significant swirl and mixing occur after the confluence of the iliac veins.  
The velocity contours and the flow patterns are compared between the present method and the OpenFOAM simulation for four slices labeled B, C, D, and E past the confluence region.
Fig.~\ref{fig:ivc_velocity} shows there is a good agreement in both velocity magnitudes and the flow patterns between the two solutions. 

\begin{figure}[t!!]
		\centering
			\includegraphics[width=0.95\textwidth]{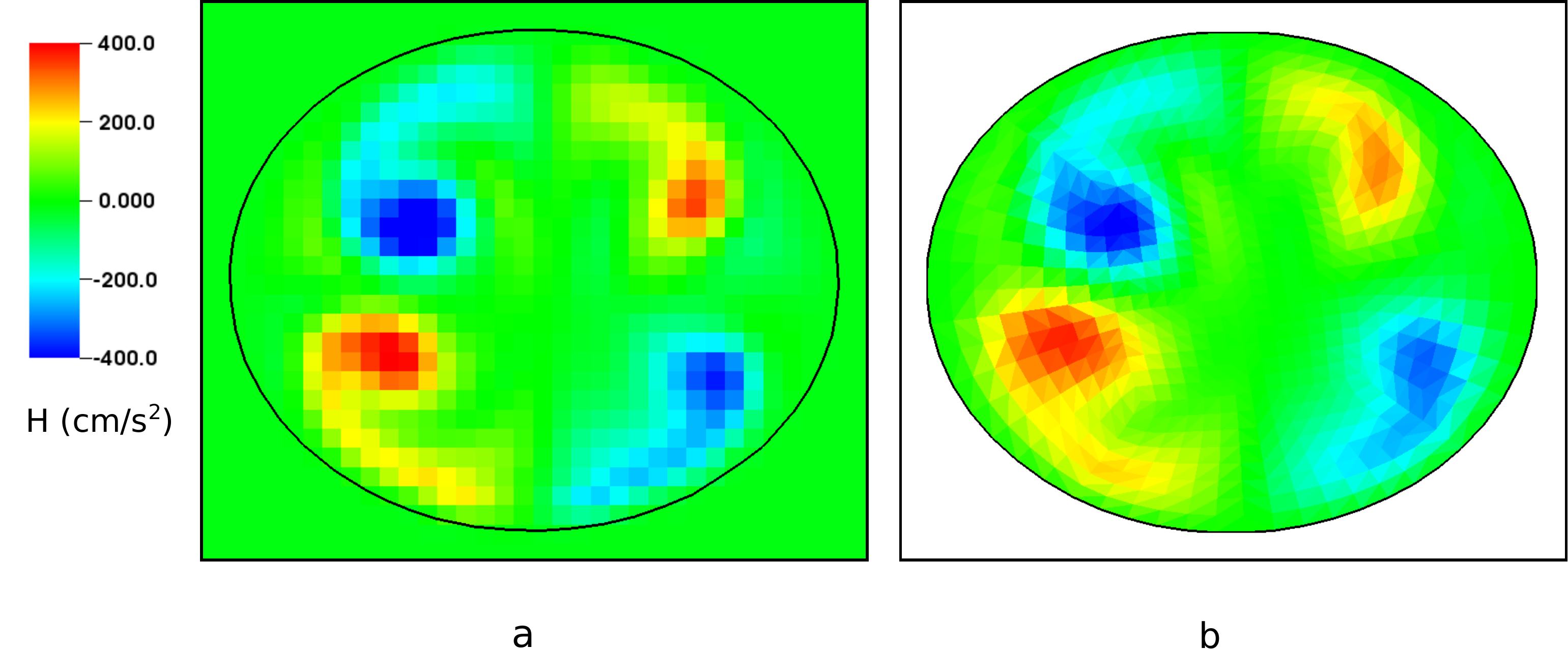}
		
		\caption{Local helicity density magnitudes in the IVC for (a) the present IIM and (b) the body-fitted OpenFOAM simulation.}
\label{fig:IIM-OpenFoam-Helicity} 
\end{figure}

To further quantify and compare the amount of swirl and mixing in the solution, we compute the helicity density $H=\u \cdot \vec{\omega}$.
This parameter is an indicator of local embolus particle transport in cardiovascular flows \cite{mukherjee2016numerical}.
Fig.~\ref{fig:IIM-OpenFoam-Helicity} illustrates the contours of helicity for the two solutions using the present IIM approach and OpenFOAM at the critical cross section E previously shown in Fig.~\ref{fig:Eulerian-mesh-ivc}.
Once again, overall there is good agreement between the two solutions.

\begin{figure}[t!!]
		\centering
			\includegraphics[width=0.8\textwidth]{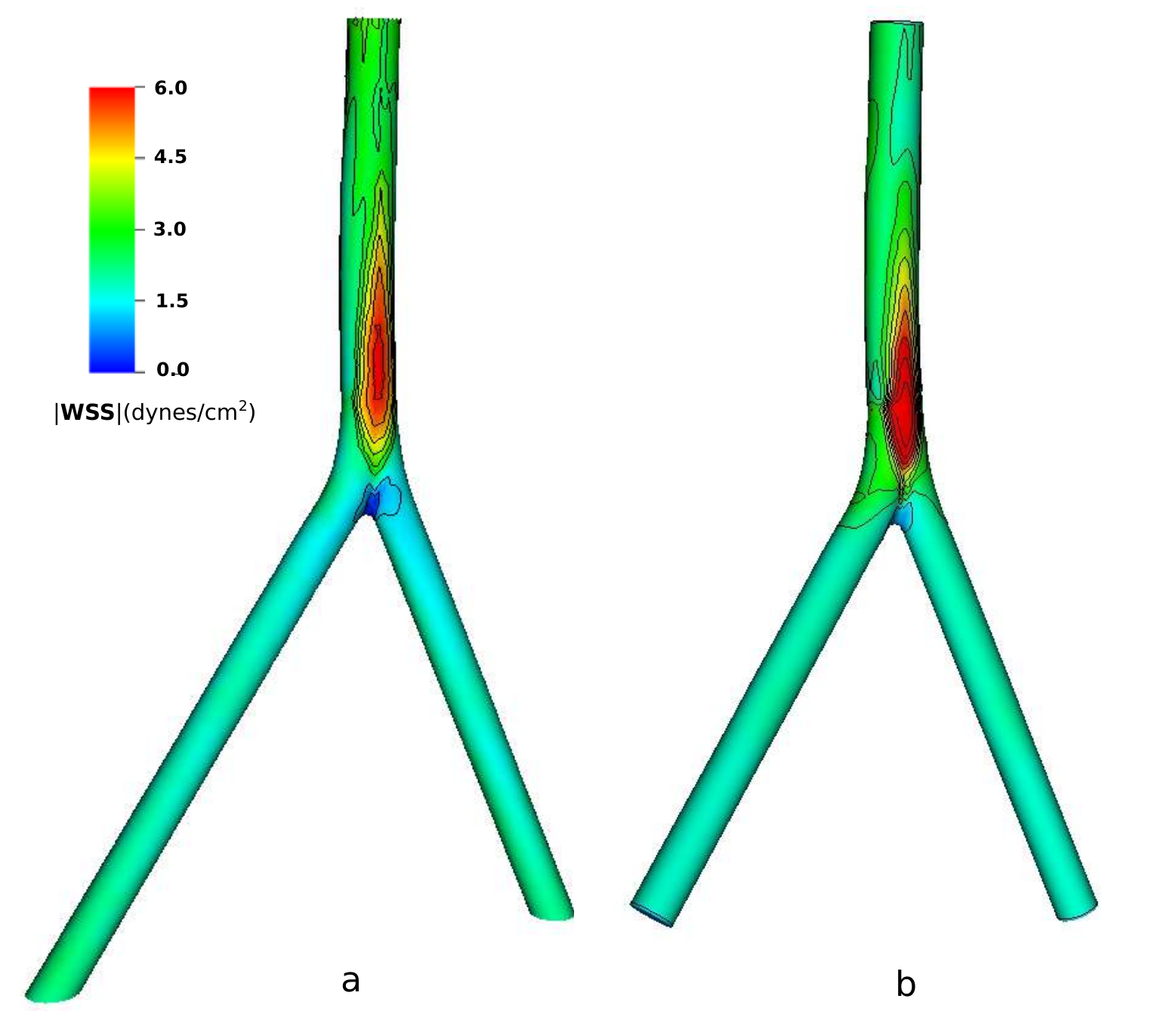}
		
\caption{Contours of wall shear stress magnitude for flow in the IVC at the exercise flow rate for (a) the present IIM and (b) the OpenFOAM simulation.}
 \label{fig:ivc_wss} 
\end{figure}

Finally, contours of wall shear stress magnitudes are plotted in Fig.~\ref{fig:ivc_wss} for both the present approach and the OpenFOAM simulation.
The overall distribution appears to be similar particularly in the region prior to the confluence of the veins and where the minimum wall shear stress occurs at the junction. However, slight differences are observed in the location where the maximum wall shear stress happens.
This is likely due to minor differences in the mesh resolution between the two cases at the iliac vein confluence due to the requirement for point-matching in the block-structured body-fitted OpenFOAM mesh, which yields a slightly coarser resolution than the IIM mesh in this region. 
Future work should more rigorously compare the two simulation methods by performing a full mesh refinement study, which is beyond the scope of the present work.

\section{Summary and conclusions}
\label{sec:conclusion}

This work describes an immersed interface method that sharply resolves fluid dynamics problems in complex geometries using a general finite element mesh description of the interface.
Like the IB method, the IIM does not require geometrically conforming spatial discretizations, which enables the use of fast Cartesian grid linear solvers. 
The primary contribution of the present work is to extend the IIM to treat general interface geometries that are described only in terms of $C^{0}$-continuous surface representations.
In particular, the present method does not require any analytic description of the interface geometry.
Methods for accurately modeling flows involving $C^{0}$-continuous immersed boundaries described by standard finite element meshes, without analytic information about the interfacial geometry, are needed to address applications involving geometries determined from experimental or clinical image data.

As in the original immersed interface methods, jumps in the pressure and velocity gradient are determined from interfacial forces imposed along the immersed boundary.
In this work, a continuous representation of these jump conditions is determined on the $C^0$ representation of the interface by projecting them onto continuous Lagrangian basis functions supported along the interface.
Although this requires the solution of a linear system of equations, this solve occurs only on the surface mesh and can be performed in an algorithmically optimal way using simple iterative methods.
Further, an accurate velocity interpolation scheme based on a modified bilinear (or, in three spatial dimensions, trilinear) scheme is developed that accounts for the imposed discontinuity in the velocity gradient.
Although the method uses only the lowest-order jump conditions, it is demonstrated to maintain global second-order accuracy in the Eulerian velocity and between first-~and second-order global accuracy in the Eulerian pressure along with second-order accuracy in the pointwise deformation and velocity of the interface, and at least first-order accuracy in the pointwise fluid stress along the interface. 
Benchmark examples in two and three spatial dimensions demonstrate that very high fidelity flow simulations are possible using this scheme and, further, that it is necessary to account for all of these lowest-order jump conditions to maximize the accuracy of the scheme.


Another novel aspect of the present method is that it uses only the lowest-order jump conditions.
It is common in immersed interface methods described in the literature to impose higher-order jumps conditions, including the jump in the first normal derivative of the pressure and the second normal derivative of the velocity.
These higher-order corrections require the evaluation of terms such as the curvature and surface divergence of the force.
Note that directly evaluating the curvature on the surface in a pointwise fashion requires $C^2$ continuity.
$C^{0}$-continuous surfaces do not provide enough regularity to obtain continuity of even the normal vector between elements.
To address this difficulty, specialized treatments have been developed that yield accurate approximations of higher-order geometrical quantities, including the surface normal and curvature, on $C^0$ continuous surface representations.
For triangulated surfaces, one approach is to adopt alternative basis functions, such as those used by subdivision surfaces \cite{cirak2000subdivision}, which provide enough regularity to yield accurate pointwise approximations to curvature and surface derivatives by directly differentiating the shape functions.
Another approach is to use methods of discrete differential geometry using averaging Voronoi cells \cite{meyer2003discrete} or finite element stabilization techniques \cite{cenanovic2017finite}.
In the immersogeometric FSI methodology  \cite{kamensky2015immersogeometric}, higher order spline-based surface representation can be used, following the original work on isogeometric analysis \cite{hughes2005isogeometric}.
Nonetheless, if global second-order accuracy in the velocity and first-order accuracy in the pressure and wall shear stress are sufficient, empirical tests reported in this paper demonstrate that only the lowest-order jump conditions are required.
Further, these results demonstrate that these jump conditions may be determined and imposed using simple, low continuity-order interface representations.
Work is underway to extend the present method to volumetric (codimension zero) problems of moving and flexible objects as well as thin flexible structures described using plate and shell models, and to apply this methodology to a range of applications in engineering, applied science, and medicine.

\section*{Acknowledgement}

We gratefully acknowledge research support through NIH Awards HL117063 and HL143336, NSF Awards DMS 1664645, CBET 175193, OAC 1450327, and OAC 1652541, and the U.S.~FDA Center for Devices and Radiological Health (CDRH) Critical Path program.
The research was supported in part by an appointment to the Research Participation Program at the U.S.~FDA administered by the Oak Ridge Institute for Science and Education through an interagency agreement between the U.S.~Department of Energy and FDA.
A.P.S.B.~also acknowledges support from San Diego State University.
Computations were performed using facilities provided by University of North Carolina at Chapel Hill through the Research Computing division of UNC Information Technology Services and the high-performance computing clusters at the U.S.~FDA.
The findings and conclusions in this article have not been formally disseminated by the FDA and should not be construed to represent any agency determination or policy.
The mention of commercial products, their sources, or their use in connection with material reported herein is not to be construed as either an actual or implied endorsement of such products by the Department of Health and Human Services.
We thank Simone Rossi and Jae Ho Lee for their constructive comments to improve the manuscript.
We also thank Kenneth Aycock for providing the patient-averaged IVC geometry.


\bibliography{AnIIMForDiscreteSurfces}


\appendix
\section{Convergence data}
\label{sec:appendix}

\begin{table}[ht]
    \scriptsize
	\centering	
	\caption{Convergence results of the Eulerian velocity and pressure for two-dimensional horizontal plane-Poiseuille considered in Sec.~\ref{subsubsec:plane-poiseuille}.
	 }\vspace{4pt}
	\label{table:error_2D_straight_pipe_Eulerian}	
\begin{tabular}{cc|c c | c c|l}
\cline{3-6}
& & \multicolumn{2}{ |c| }{ velocity} & \multicolumn{2}{ |c| }{ pressure }  \\  \cline{2-6}
&  \multicolumn{1}{ |c| }{h} & Error & Order & Error & Order & \\ \cline{1-6}
\multicolumn{1}{ |c  }{\multirow{4}{*}{$\L2$  } } &
\multicolumn{1}{ |c| }{$1.56\times 10^{-1}$} & $4.4886\times 10^{-2}$ & - & $2.9764\times 10^{-2}$ & - &    \\ \cline{2-6}
\multicolumn{1}{ |c  }{}                        &
\multicolumn{1}{ |c| }{$7.80\times 10^{-2}$} & $1.0044\times 10^{-2}$ & 2.16 & $1.2007\times 10^{-2}$ & 1.31 &    \\ \cline{2-6}
\multicolumn{1}{ |c  }{}                        &
\multicolumn{1}{ |c| }{$3.90\times 10^{-2}$} & $1.3822\times 10^{-3}$ & 2.86 & $2.5519\times 10^{-3}$ & 2.23 &    \\ \cline{2-6}
\multicolumn{1}{ |c  }{}                        &
\multicolumn{1}{ |c| }{$1.95\times 10^{-3}$} &  $3.6201\times 10^{-4}$ & 1.93 & $2.3269\times 10^{-4}$ & 3.46 &    \\ \hhline{======}
\multicolumn{1}{ |c  }{\multirow{4}{*}{$\Linf$ } } &
\multicolumn{1}{ |c| }{$1.56\times 10^{-1}$} &  $4.6171\times 10^{-2}$ & - & $3.4157\times 10^{-2}$ & - & \\ \cline{2-6}
\multicolumn{1}{ |c  }{}                        &
\multicolumn{1}{ |c| }{$7.80\times 10^{-2}$} &  $1.0396\times 10^{-2}$ & 2.15 & $1.2673\times 10^{-2}$ & 1.43 &\\ \cline{2-6}
\multicolumn{1}{ |c  }{}                        &
\multicolumn{1}{ |c| }{$3.90\times 10^{-2}$} &  $2.5911\times 10^{-3}$ & 2.00 & $2.7193\times 10^{-3}$ & 2.22 &    \\ \cline{2-6}
\multicolumn{1}{ |c  }{}                        &
\multicolumn{1}{ |c| }{$1.95\times 10^{-3}$} &  $8.4275\times 10^{-4}$ & 1.62 & $2.6682\times 10^{-4}$ & 3.35 &    \\ \cline{1-6}
\end{tabular}
\end{table} 
\begin{table}[ht]
    \scriptsize
	\centering	
	\caption{Convergence results of the displacement, velocity, pressure and WSS on the surface mesh of the two-dimensional horizontal plane-Poiseuille flow considered in Sec.~\ref{subsubsec:plane-poiseuille}.
	 }\vspace{4pt}
	\label{table:error_2D_straight_pipe_Lagrangian}	
\begin{tabular}{cc|c c |c c | c c| c c|l}
\cline{3-10}
& & \multicolumn{2}{ |c| }{displacement} & \multicolumn{2}{ |c| }{ velocity} & \multicolumn{2}{ |c| }{ pressure } & \multicolumn{2}{ |c| }{ WSS } \\  \cline{2-10}
&  \multicolumn{1}{ |c| }{h} & Error & Order & Error & Order & Error & Order & Error & Order & \\ \cline{1-10}
\multicolumn{1}{ |c  }{\multirow{4}{*}{$\L2$  } } &
\multicolumn{1}{ |c| }{$1.56\times 10^{-1}$} & $4.2943\times 10^{-3}$ & - & $8.6485\times 10^{-3}$ & - & $4.2904\times 10^{-2}$ & - & $1.9004\times 10^{-2}$ & - &    \\ \cline{2-10}
\multicolumn{1}{ |c  }{}                        &
\multicolumn{1}{ |c| }{$7.80\times 10^{-2}$} & $1.0702\times 10^{-3}$ & 2.00 & $2.1488\times 10^{-3}$ & 2.01 & $1.5327\times 10^{-2}$ & 1.49 & $8.9271\times 10^{-3}$ & 1.09 &    \\ \cline{2-10}
\multicolumn{1}{ |c  }{}                        &
\multicolumn{1}{ |c| }{$3.90\times 10^{-2}$} & $2.8847\times 10^{-4}$ & 1.89 & $8.6136\times 10^{-4}$ & 1.32 & $1.8084\times 10^{-3}$ & 2.55 & $4.2621\times 10^{-3}$ & 1.07 &    \\ \cline{2-10}
\multicolumn{1}{ |c  }{}                        &
\multicolumn{1}{ |c| }{$1.95\times 10^{-3}$} & $7.1420\times 10^{-5}$ & 2.01 & $2.1852\times 10^{-4}$ & 1.98 & $2.0321\times 10^{-4}$ & 3.15 & $2.1631\times 10^{-3}$ & 0.98 &    \\ \hhline{==========}
\multicolumn{1}{ |c  }{\multirow{4}{*}{$\Linf$ } } &
\multicolumn{1}{ |c| }{$1.56\times 10^{-1}$} &  $5.0088\times 10^{-3}$ & - & $1.0645\times 10^{-2}$ & -   & $1.0489\times 10^{-2}$ & - & $1.1273\times 10^{-2}$ & - & \\ \cline{2-10}
\multicolumn{1}{ |c  }{}                        &
\multicolumn{1}{ |c| }{$7.80\times 10^{-2}$} &  $1.0452\times 10^{-3}$ & 2.26 &  $2.9337\times 10^{-3}$ & 1.86 & $4.3378\times 10^{-3}$ & 1.27 & $5.3873\times 10^{-3}$ & 1.07 &\\ \cline{2-10}
\multicolumn{1}{ |c  }{}                        &
\multicolumn{1}{ |c| }{$3.90\times 10^{-2}$} &  $2.4787\times 10^{-4}$ & 2.08 & $1.5121\times 10^{-3}$ & 0.96 & $6.8581\times 10^{-4}$ & 2.66 & $2.5205\times 10^{-3}$ & 1.09 &    \\ \cline{2-10}
\multicolumn{1}{ |c  }{}                        &
\multicolumn{1}{ |c| }{$1.95\times 10^{-3}$} &  $6.1995\times 10^{-5}$ & 2.00 & $4.0972\times 10^{-4}$ & 1.88 & $1.0069\times 10^{-4}$ & 2.78 & $1.2724\times 10^{-3}$ & 0.99 &    \\ \cline{1-10}
\end{tabular}
\end{table} 
\begin{table}[ht]
    \scriptsize
	\centering	
	\caption{Convergence results of the Eulerian velocity and pressure for two-dimensional inclined plane-Poiseuille flow considered in Sec.~\ref{subsubsec:plane-poiseuille}.
	}\vspace{4pt}
	\label{table:error_2D_inclined_pipe_Eulerian}	
\begin{tabular}{cc|c c | c c| c c|l}
\cline{3-8}
& & \multicolumn{2}{ |c| }{ velocity} & \multicolumn{2}{ |c| }{ pressure  } & \multicolumn{2}{ |c| }{ pressure on $\Omega^*$ } \\  \cline{2-8}
&  \multicolumn{1}{ |c| }{h} & Error & Order & Error & Order & Error & Order & \\ \cline{1-8}
\multicolumn{1}{ |c  }{\multirow{4}{*}{$\L2$  } } &
\multicolumn{1}{ |c| }{$7.80\times 10^{-2}$} & $1.3072\times 10^{-3}$ & - & $7.4462\times 10^{-4}$ & - & - & - &    \\ \cline{2-8}
\multicolumn{1}{ |c  }{}                        &
\multicolumn{1}{ |c| }{$3.90\times 10^{-2}$} & $3.1417\times 10^{-4}$ & 2.06 & $2.2671\times 10^{-4}$ & 1.72 & - & - &    \\ \cline{2-8}
\multicolumn{1}{ |c  }{}                        &
\multicolumn{1}{ |c| }{$1.95\times 10^{-3}$} & $7.6507\times 10^{-5}$ & 2.04 & $9.4735\times 10^{-5}$ & 1.26 & - & - &    \\ \cline{2-8}
\multicolumn{1}{ |c  }{}                        &
\multicolumn{1}{ |c| }{$9.75\times 10^{-4}$} & $1.9933\times 10^{-5}$ & 1.94 & $4.0041\times 10^{-5}$ & 1.24 & - & - &    \\ \hhline{========}
\multicolumn{1}{ |c  }{\multirow{4}{*}{$\Linf$ } } &
\multicolumn{1}{ |c| }{$7.80\times 10^{-2}$} &  $5.8251\times 10^{-2}$ & - & $4.5021\times 10^{-2}$ & - & $5.2985\times 10^{-2}$ & - & \\ \cline{2-8}
\multicolumn{1}{ |c  }{}                        &
\multicolumn{1}{ |c| }{$3.90\times 10^{-2}$} &  $3.0388\times 10^{-2}$ & 0.94 & $2.1777\times 10^{-2}$ & 1.48 & $1.9759\times 10^{-2}$ & 1.42 &\\ \cline{2-8}
\multicolumn{1}{ |c  }{}                        &
\multicolumn{1}{ |c| }{$1.95\times 10^{-3}$} &  $1.5299\times 10^{-2}$ & 0.99 & $1.5221\times 10^{-2}$ & 0.52 & $1.0058 \times 10^{-2}$ & 0.97 &    \\ \cline{2-8}
\multicolumn{1}{ |c  }{}                        &
\multicolumn{1}{ |c| }{$9.75\times 10^{-4}$} &  $7.6174\times 10^{-3}$ & 1.00 & $1.4029\times 10^{-2}$ & 0.12 & $5.0311\times 10^{-3}$ & 1.00 &    \\ \cline{1-8}
\end{tabular}
\end{table} 
\begin{table}[ht]
    \scriptsize
	\centering	
	\caption{Convergence results of the displacement, velocity, pressure and WSS on the surface mesh of the two-dimensional inclined plane-Poiseuille flow considered in Sec.~\ref{subsubsec:plane-poiseuille}.
	 }\vspace{4pt}
	\label{table:error_2D_inclined_pipe_Lagrangian}	
\begin{tabular}{cc|c c |c c | c c| c c|l}
\cline{3-10}
& & \multicolumn{2}{ |c| }{displacement} & \multicolumn{2}{ |c| }{ velocity} & \multicolumn{2}{ |c| }{ Pressure } & \multicolumn{2}{ |c| }{ WSS } \\  \cline{2-10}
&  \multicolumn{1}{ |c| }{h} & Error & Order & Error & Order & Error & Order & Error & Order & \\ \cline{1-10}
\multicolumn{1}{ |c  }{\multirow{4}{*}{$\L2$  } } &
\multicolumn{1}{ |c| }{$7.80\times 10^{-2}$} & $1.0256\times 10^{-3}$ & - & $5.4813\times 10^{-3}$ & - & $3.4355\times 10^{-2}$ & - & $1.0341\times 10^{-2}$ & - &    \\ \cline{2-10}
\multicolumn{1}{ |c  }{}                        &
\multicolumn{1}{ |c| }{$3.90\times 10^{-2}$} & $2.5688\times 10^{-4}$ & 2.00 & $8.1599\times 10^{-4}$ & 2.75 & $1.4422\times 10^{-2}$ & 1.25 & $4.9935\times 10^{-2}$ & 1.05 &    \\ \cline{2-10}
\multicolumn{1}{ |c  }{}                        &
\multicolumn{1}{ |c| }{$1.95\times 10^{-3}$} & $6.5622\times 10^{-5}$ & 1.97 & $2.3934\times 10^{-4}$ & 1.77 & $4.0782\times 10^{-3}$ & 1.82 & $2.3804\times 10^{-2}$ & 1.07 &    \\ \cline{2-10}
\multicolumn{1}{ |c  }{}                        &
\multicolumn{1}{ |c| }{$9.75\times 10^{-4}$} & $1.6489\times 10^{-5}$ & 1.99 & $6.2818\times 10^{-5}$ & 1.93 & $1.1990\times 10^{-3}$ & 1.77 & $1.1239\times 10^{-3}$ & 1.08 &    \\ \hhline{==========}
\multicolumn{1}{ |c  }{\multirow{4}{*}{$\Linf$ } } &
\multicolumn{1}{ |c| }{$7.80\times 10^{-2}$} &  $1.3363\times 10^{-3}$ & - & $7.0881\times 10^{-3}$ & -   & $2.9897\times 10^{-2}$ & - & $6.9130\times 10^{-3}$ & - & \\ \cline{2-10}
\multicolumn{1}{ |c  }{}                        &
\multicolumn{1}{ |c| }{$3.90\times 10^{-2}$} &  $2.9123\times 10^{-4}$ & 2.20 &  $1.1809\times 10^{-3}$ & 2.59 & $1.2979\times 10^{-2}$ & 1.20 & $3.1308\times 10^{-3}$ & 1.14 &\\ \cline{2-10}
\multicolumn{1}{ |c  }{}                        &
\multicolumn{1}{ |c| }{$1.95\times 10^{-3}$} &  $7.0620\times 10^{-5}$ & 2.04 & $4.1602\times 10^{-4}$ & 1.51 & $4.5329\times 10^{-3}$ & 1.52 & $1.7871\times 10^{-3}$ & 0.81 &    \\ \cline{2-10}
\multicolumn{1}{ |c  }{}                        &
\multicolumn{1}{ |c| }{$9.75\times 10^{-4}$} &  $1.8191\times 10^{-5}$ & 1.96 & $1.4445\times 10^{-4}$ & 1.53 & $1.5480\times 10^{-3}$ & 1.55 & $8.7958\times 10^{-4}$ & 1.02 &    \\ \cline{1-10}
\end{tabular}
\end{table} 
\begin{table}[ht]
    \scriptsize
	\centering	
	\caption{Convergence results of the Eulerian velocity and pressure for three-dimensional horizontal Hagen-Poiseuille flow considered in Sec.~\ref{subsubsec:Hagen-poiseuille}.
	}\vspace{4pt}
	\label{table:error_3D_straight_pipe_Eulerian}	
\begin{tabular}{cc|c c | c c| c c|l}
\cline{3-8}
& & \multicolumn{2}{ |c| }{ velocity} & \multicolumn{2}{ |c| }{ pressure } & \multicolumn{2}{ |c| }{ pressure on $\Omega^*$ } \\  \cline{2-8}
&  \multicolumn{1}{ |c| }{h} & Error & Order & Error & Order & Error & Order & \\ \cline{1-8}
\multicolumn{1}{ |c  }{\multirow{4}{*}{$\L2$  } } &
\multicolumn{1}{ |c| }{$1.56\times 10^{-1}$} & $1.1629\times 10^{-1}$ & - & $6.0942\times 10^{-2}$ & - & - & - &    \\ \cline{2-8}
\multicolumn{1}{ |c  }{}                        &
\multicolumn{1}{ |c| }{$7.80\times 10^{-2}$} & $2.9558\times 10^{-2}$ & 1.98 & $1.6789\times 10^{-2}$ & 1.86 & - & - &    \\ \cline{2-8}
\multicolumn{1}{ |c  }{}                        &
\multicolumn{1}{ |c| }{$3.90\times 10^{-2}$} & $7.6507\times 10^{-3}$ & 2.04 & $2.4049\times 10^{-3}$ & 2.80 & - & - &    \\ \cline{2-8}
\multicolumn{1}{ |c  }{}                        &
\multicolumn{1}{ |c| }{$1.95\times 10^{-3}$} & $1.7100\times 10^{-3}$ & 2.07 & $6.4996\times 10^{-4}$ & 1.89 & - & - &    \\ \hhline{========}
\multicolumn{1}{ |c  }{\multirow{4}{*}{$\Linf$ } } &
\multicolumn{1}{ |c| }{$1.56\times 10^{-1}$} &  $4.8008\times 10^{-2}$ & - & $1.1099\times 10^{-1}$ & - & $1.1099\times 10^{-1}$ & - & \\ \cline{2-8}
\multicolumn{1}{ |c  }{}                        &
\multicolumn{1}{ |c| }{$7.80\times 10^{-2}$} &  $2.3307\times 10^{-2}$ & 1.04 & $6.9275\times 10^{-2}$ & 0.68 & $2.4625\times 10^{-2}$ & 2.17 &\\ \cline{2-8}
\multicolumn{1}{ |c  }{}                        &
\multicolumn{1}{ |c| }{$3.90\times 10^{-2}$} &  $4.7934\times 10^{-3}$ & 2.28 & $4.4148\times 10^{-2}$ & 0.65 & $3.0722 \times 10^{-3}$ & 3.00 &    \\ \cline{2-8}
\multicolumn{1}{ |c  }{}                        &
\multicolumn{1}{ |c| }{$1.95\times 10^{-3}$} &  $1.2291\times 10^{-3}$ & 1.96 & $2.5007\times 10^{-2}$ & 0.82 & $7.5856\times 10^{-4}$ & 2.02 &    \\ \cline{1-8}
\end{tabular}
\end{table} 
\begin{table}[ht]
    \scriptsize
	\centering	
	\caption{Convergence results of the displacement, velocity, pressure and WSS on the surface mesh of the three-dimensional horizontal Hagen-Poiseuille flow considered in Sec.~\ref{subsubsec:Hagen-poiseuille}.
	 }\vspace{4pt}
	\label{table:error_3D_straight_pipe_Lagrangian}	
\begin{tabular}{cc|c c |c c | c c| c c|l}
\cline{3-10}
& & \multicolumn{2}{ |c| }{displacement} & \multicolumn{2}{ |c| }{ velocity} & \multicolumn{2}{ |c| }{ pressure } & \multicolumn{2}{ |c| }{ WSS } \\  \cline{2-10}
&  \multicolumn{1}{ |c| }{h} & Error & Order & Error & Order & Error & Order & Error & Order & \\ \cline{1-10}
\multicolumn{1}{ |c  }{\multirow{4}{*}{$\L2$  } } &
\multicolumn{1}{ |c| }{$1.56\times 10^{-1}$} & $5.8931\times 10^{-3}$ & - & $4.4541\times 10^{-3}$ & - & $1.7140\times 10^{-1}$ & - & $2.3026\times 10^{-2}$ & - &    \\ \cline{2-10}
\multicolumn{1}{ |c  }{}                        &
\multicolumn{1}{ |c| }{$7.80\times 10^{-2}$} & $8.5408\times 10^{-4}$ & 2.79 & $4.8937\times 10^{-4}$ & 3.19 & $4.1685\times 10^{-2}$ & 2.04 & $1.2835\times 10^{-2}$ & 1.09 &    \\ \cline{2-10}
\multicolumn{1}{ |c  }{}                        &
\multicolumn{1}{ |c| }{$3.90\times 10^{-2}$} & $2.1148\times 10^{-4}$ & 2.01 & $4.0652\times 10^{-4}$ & 0.26 & $5.2332\times 10^{-3}$ & 2.99 & $5.8642\times 10^{-3}$ & 1.07 &    \\ \cline{2-10}
\multicolumn{1}{ |c  }{}                        &
\multicolumn{1}{ |c| }{$1.95\times 10^{-3}$} & $5.2629\times 10^{-5}$ & 2.00 & $9.9410\times 10^{-5}$ & 2.03 & $1.2725\times 10^{-3}$ & 2.04 & $2.7800\times 10^{-3}$ & 0.98 &    \\ \hhline{==========}
\multicolumn{1}{ |c  }{\multirow{4}{*}{$\Linf$ } } &
\multicolumn{1}{ |c| }{$1.56\times 10^{-1}$} &  $4.6561\times 10^{-3}$ & - & $8.5493\times 10^{-3}$ & -   & $1.1037\times 10^{-1}$ & - & $1.4015\times 10^{-2}$ & - & \\ \cline{2-10}
\multicolumn{1}{ |c  }{}                        &
\multicolumn{1}{ |c| }{$7.80\times 10^{-2}$} &  $7.4997\times 10^{-4}$ & 2.63&  $1.6748\times 10^{-3}$ & 2.35 & $2.4535\times 10^{-2}$ & 1.27 & $7.5266\times 10^{-3}$ & 0.90 &\\ \cline{2-10}
\multicolumn{1}{ |c  }{}                        &
\multicolumn{1}{ |c| }{$3.90\times 10^{-2}$} &  $1.5650\times 10^{-4}$ & 2.26 & $1.0312\times 10^{-3}$ & 0.70 & $3.0398\times 10^{-3}$ & 2.66 & $3.5189\times 10^{-3}$ & 1.10 &    \\ \cline{2-10}
\multicolumn{1}{ |c  }{}                        &
\multicolumn{1}{ |c| }{$1.95\times 10^{-3}$} &  $3.9278\times 10^{-5}$ & 1.99 & $3.0996\times 10^{-4}$ & 1.73 & $7.2375\times 10^{-4}$ & 2.78 & $1.7256\times 10^{-3}$ & 1.03 &    \\ \cline{1-10}
\end{tabular}
\end{table} 
\begin{table}[ht]
    \scriptsize
	\centering	
	\caption{Convergence results of the Eulerian velocity and pressure for three-dimensional inclined Hagen-Poiseuille flow considered in Sec.~\ref{subsubsec:Hagen-poiseuille}.
	 }\vspace{4pt}
	\label{table:error_3D_inclined_pipe_Eulerian}	
\begin{tabular}{cc|c c | c c| c c|l}
\cline{3-8}
& & \multicolumn{2}{ |c| }{ velocity} & \multicolumn{2}{ |c| }{ pressure  } & \multicolumn{2}{ |c| }{ pressure on $\Omega^*$ } \\  \cline{2-8}
&  \multicolumn{1}{ |c| }{h} & Error & Order & Error & Order & Error & Order & \\ \cline{1-8}
\multicolumn{1}{ |c  }{\multirow{4}{*}{$\L2$  } } &
\multicolumn{1}{ |c| }{$7.80\times 10^{-2}$} & $2.7733\times 10^{-3}$ & - & $6.5465\times 10^{-2}$ & - & - & - &    \\ \cline{2-8}
\multicolumn{1}{ |c  }{}                        &
\multicolumn{1}{ |c| }{$3.90\times 10^{-2}$} & $4.6659\times 10^{-4}$ & 2.57 & $1.8020\times 10^{-3}$ & 1.86 & - & - &    \\ \cline{2-8}
\multicolumn{1}{ |c  }{}                        &
\multicolumn{1}{ |c| }{$1.95\times 10^{-3}$} & $1.8128\times 10^{-4}$ & 1.36 & $1.0238\times 10^{-3}$ & 0.82 & - & - &    \\ \cline{2-8}
\multicolumn{1}{ |c  }{}                        &
\multicolumn{1}{ |c| }{$9.75\times 10^{-4}$} & $4.5095\times 10^{-5}$ & 2.01 & $3.6563\times 10^{-4}$ & 1.49 & - & - &    \\ \hhline{========}
\multicolumn{1}{ |c  }{\multirow{4}{*}{$\Linf$ } } &
\multicolumn{1}{ |c| }{$7.80\times 10^{-2}$} &  $8.4694\times 10^{-3}$ & - & $4.5885\times 10^{-2}$ & - & $3.9710\times 10^{-2}$ & - & \\ \cline{2-8}
\multicolumn{1}{ |c  }{}                        &
\multicolumn{1}{ |c| }{$3.90\times 10^{-2}$} &  $2.2634\times 10^{-3}$ & 1.90 & $4.0355\times 10^{-2}$ & 0.19 & $5.7205\times 10^{-3}$ & 2.80 &\\ \cline{2-8}
\multicolumn{1}{ |c  }{}                        &
\multicolumn{1}{ |c| }{$1.95\times 10^{-3}$} &  $7.9415\times 10^{-4}$ & 1.51 & $2.1088\times 10^{-3}$ & 0.94 & $2.1088\times 10^{-3}$ & 1.44 &    \\ \cline{2-8}
\multicolumn{1}{ |c  }{}                        &
\multicolumn{1}{ |c| }{$9.75\times 10^{-4}$} &  $2.2101\times 10^{-4}$ & 1.85 & $1.4058\times 10^{-3}$ & 0.59 & $1.0088\times 10^{-3}$ & 1.06 &    \\ \cline{1-8}
\end{tabular}
\end{table} 
\begin{table}[ht]
    \scriptsize
	\centering	
	\caption{Convergence results of the displacement, velocity, pressure and WSS on the surface mesh of the three-dimensional inclined Hagen-Poiseuille flow considered in Sec.~\ref{subsubsec:Hagen-poiseuille}.
	}\vspace{4pt}
	\label{table:error_3D_inclined_pipe_Lagrangian}	
\begin{tabular}{cc|c c |c c | c c| c c|l}
\cline{3-10}
& & \multicolumn{2}{ |c| }{displacement} & \multicolumn{2}{ |c| }{ velocity} & \multicolumn{2}{ |c| }{ pressure } & \multicolumn{2}{ |c| }{ WSS } \\  \cline{2-10}
&  \multicolumn{1}{ |c| }{h} & Error & Order & Error & Order & Error & Order & Error & Order & \\ \cline{1-10}
\multicolumn{1}{ |c  }{\multirow{4}{*}{$\L2$  } } &
\multicolumn{1}{ |c| }{$7.80\times 10^{-2}$} & $2.8315\times 10^{-3}$ & - & $1.8397\times 10^{-5}$ & - & $3.1258\times 10^{-3}$ & - & $8.3374\times 10^{-3}$ & - &    \\ \cline{2-10}
\multicolumn{1}{ |c  }{}                        &
\multicolumn{1}{ |c| }{$3.90\times 10^{-2}$} & $4.1330\times 10^{-4}$ & 2.78 & $5.0375\times 10^{-6}$ & 1.87 & $3.0608\times 10^{-4}$ & 3.35 & $4.7784\times 10^{-3}$ & 0.80 &    \\ \cline{2-10}
\multicolumn{1}{ |c  }{}                        &
\multicolumn{1}{ |c| }{$1.95\times 10^{-3}$} & $1.0007\times 10^{-4}$ & 2.05 & $3.0499\times 10^{-6}$ & 0.72 & $1.1570\times 10^{-4}$ & 1.40 & $2.4453\times 10^{-3}$ & 0.97 &    \\ \cline{2-10}
\multicolumn{1}{ |c  }{}                        &
\multicolumn{1}{ |c| }{$9.75\times 10^{-4}$} & $2.4407\times 10^{-5}$ & 2.04 & $7.8403\times 10^{-7}$ & 1.96 & $3.6156\times 10^{-5}$ & 1.68 & $1.2870\times 10^{-3}$ & 0.93 &    \\ \hhline{==========}
\multicolumn{1}{ |c  }{\multirow{4}{*}{$\Linf$ } } &
\multicolumn{1}{ |c| }{$7.80\times 10^{-2}$} &  $2.1752\times 10^{-3}$ & - & $1.1652\times 10^{-5}$ & -   & $3.1373\times 10^{-3}$ & - & $5.7894\times 10^{-2}$ & - & \\ \cline{2-10}
\multicolumn{1}{ |c  }{}                        &
\multicolumn{1}{ |c| }{$3.90\times 10^{-2}$} &  $3.3874\times 10^{-4}$ & 2.68 &  $5.4306\times 10^{-6}$ & 1.10 & $3.3588\times 10^{-4}$ & 3.22 & $2.9846\times 10^{-3}$ & 0.96 &\\ \cline{2-10}
\multicolumn{1}{ |c  }{}                        &
\multicolumn{1}{ |c| }{$1.95\times 10^{-3}$} &  $7.9366\times 10^{-5}$ & 2.09 & $3.0997\times 10^{-6}$ & 0.81 & $1.8660\times 10^{-4}$ & 0.85 & $1.5948\times 10^{-3}$ & 0.90 &    \\ \cline{2-10}
\multicolumn{1}{ |c  }{}                        &
\multicolumn{1}{ |c| }{$9.75\times 10^{-4}$} &  $1.8852\times 10^{-5}$ & 2.07 & $1.4285\times 10^{-6}$ & 1.12 & $9.8211\times 10^{-5}$ & 0.93 & $8.0262\times 10^{-3}$ & 0.99 &    \\ \cline{1-10}
\end{tabular}
\end{table} 
\begin{table}[ht]
    \scriptsize
	\centering	
	\caption{Convergence results of the Eulerian velocity and pressure for two-dimensional circular Couette flow considered in Sec.~\ref{subsec:TCF}.
	 }\vspace{4pt}
	\label{table:error_2D_TCF_Eulerian}	
\begin{tabular}{cc|c c | c c| c c|l}
\cline{3-8}
& & \multicolumn{2}{ |c| }{ velocity} & \multicolumn{2}{ |c| }{ pressure } & \multicolumn{2}{ |c| }{ pressure on $\Omega^*$ } \\  \cline{2-8}
&  \multicolumn{1}{ |c| }{h} & Error & Order & Error & Order & Error & Order & \\ \cline{1-8}
\multicolumn{1}{ |c  }{\multirow{4}{*}{$\L2$  } } &
\multicolumn{1}{ |c| }{$1.2500\times 10^{-1}$} & $3.9912\times 10^{-2}$ & - & $5.7864\times 10^{-2}$ & - & - & - &    \\ \cline{2-8}
\multicolumn{1}{ |c  }{}                        &
\multicolumn{1}{ |c| }{$6.2500\times 10^{-2}$} & $1.0652\times 10^{-2}$ & 1.91 & $1.6002\times 10^{-2}$ & 1.85 & - & - &    \\ \cline{2-8}
\multicolumn{1}{ |c  }{}                        &
\multicolumn{1}{ |c| }{$3.1250\times 10^{-2}$} & $3.6850\times 10^{-3}$ & 1.53 & $4.4463\times 10^{-3}$ & 1.85 & - & - &    \\ \cline{2-8}
\multicolumn{1}{ |c  }{}                        &
\multicolumn{1}{ |c| }{$1.5625\times 10^{-2}$} & $9.5675\times 10^{-4}$ & 1.95 & $1.6014\times 10^{-3}$ & 1.74 & - & - &    \\ \hhline{========}
\multicolumn{1}{ |c  }{\multirow{4}{*}{$\Linf$ } } &
\multicolumn{1}{ |c| }{$1.2500\times 10^{-1}$} &  $1.2378\times 10^{-1}$ & - & $4.5941\times 10^{-2}$ & - & $2.7895\times 10^{-2}$ & - & \\ \cline{2-8}
\multicolumn{1}{ |c  }{}                        &
\multicolumn{1}{ |c| }{$6.2500\times 10^{-2}$} &  $4.0347\times 10^{-2}$ & 1.62 & $4.0176\times 10^{-2}$ & 0.19 & $1.3461\times 10^{-2}$ & 1.05 &\\ \cline{2-8}
\multicolumn{1}{ |c  }{}                        &
\multicolumn{1}{ |c| }{$3.1250\times 10^{-2}$} &  $2.2022\times 10^{-2}$ & 0.87 & $2.8997\times 10^{-2}$ & 0.47 & $6.8975\times 10^{-3}$ & 0.96 &    \\ \cline{2-8}
\multicolumn{1}{ |c  }{}                        &
\multicolumn{1}{ |c| }{$1.5625\times 10^{-2}$} &  $5.7101\times 10^{-3}$ & 1.95 & $1.8695\times 10^{-2}$ & 0.63 & $3.4660\times 10^{-3}$ & 0.99 &    \\ \cline{1-8}
\end{tabular}
\end{table} 
\begin{table}[ht]
    \scriptsize
	\centering	
	\caption{Convergence results of the displacement, velocity, pressure and WSS on the surface mesh of the two-dimensional circular Couette flow considered in Sec.~\ref{subsec:TCF}.
	 }\vspace{4pt}
	\label{table:error_TCF_Lagrangian}	
\begin{tabular}{cc|c c |c c | c c| c c|l}
\cline{3-10}
& & \multicolumn{2}{ |c| }{displacement} & \multicolumn{2}{ |c| }{ velocity} & \multicolumn{2}{ |c| }{ pressure } & \multicolumn{2}{ |c| }{ WSS } \\  \cline{2-10}
&  \multicolumn{1}{ |c| }{h} & Error & Order & Error & Order & Error & Order & Error & Order & \\ \cline{1-10}
\multicolumn{1}{ |c  }{\multirow{4}{*}{$\L2$  } } &
\multicolumn{1}{ |c| }{$1.2500\times 10^{-1}$} & $3.4268\times 10^{-3}$ & - & $4.2827\times 10^{-2}$ & - & $2.7572\times 10^{-1}$ & - & $2.5266\times 10^{-2}$ & - &    \\ \cline{2-10}
\multicolumn{1}{ |c  }{}                        &
\multicolumn{1}{ |c| }{$6.2500\times 10^{-2}$} & $8.5475\times 10^{-4}$ & 2.00 & $1.2115\times 10^{-2}$ & 1.82 & $6.7620\times 10^{-1}$ & 2.03 & $1.4307\times 10^{-2}$ & 0.82 &    \\ \cline{2-10}
\multicolumn{1}{ |c  }{}                        &
\multicolumn{1}{ |c| }{$3.1250\times 10^{-2}$} & $2.0304\times 10^{-4}$ & 2.07 & $3.6922\times 10^{-3}$ & 1.71 & $1.7609\times 10^{-2}$ & 1.94 & $8.2638\times 10^{-3}$ & 0.79 &    \\ \cline{2-10}
\multicolumn{1}{ |c  }{}                        &
\multicolumn{1}{ |c| }{$1.5625\times 10^{-2}$} & $5.0057\times 10^{-5}$ & 2.02 & $9.3743\times 10^{-4}$ & 1.98 & $6.4982\times 10^{-3}$ & 1.44 & $4.3695\times 10^{-3}$ & 0.92 &    \\ \hhline{==========}
\multicolumn{1}{ |c  }{\multirow{4}{*}{$\Linf$ } } &
\multicolumn{1}{ |c| }{$1.2500\times 10^{-1}$} &  $3.2031\times 10^{-3}$ & - & $3.6900\times 10^{-2}$ & -   & $1.6974\times 10^{-1}$ & - & $1.5876\times 10^{-2}$ & - & \\ \cline{2-10}
\multicolumn{1}{ |c  }{}                        &
\multicolumn{1}{ |c| }{$6.2500\times 10^{-2}$} &  $6.7864\times 10^{-4}$ & 2.24 &  $1.1092\times 10^{-2}$ & 1.73 & $4.7125\times 10^{-2}$ & 1.85 & $8.3447\times 10^{-3}$ & 0.93 &\\ \cline{2-10}
\multicolumn{1}{ |c  }{}                        &
\multicolumn{1}{ |c| }{$3.1250\times 10^{-2}$} &  $1.5153\times 10^{-4}$ & 2.16 & $3.8568\times 10^{-3}$ & 1.52 & $1.7521\times 10^{-2}$ & 1.43 & $4.0669\times 10^{-3}$ & 1.04 &    \\ \cline{2-10}
\multicolumn{1}{ |c  }{}                        &
\multicolumn{1}{ |c| }{$1.5625\times 10^{-2}$} &  $3.4132\times 10^{-5}$ & 2.15 & $9.3019\times 10^{-4}$ & 2.05 & $8.4316\times 10^{-3}$ & 1.06 & $2.2926\times 10^{-3}$ & 0.83 &    \\ \cline{1-10}
\end{tabular}
\end{table} 
\begin{table}[ht]
    \scriptsize
	\centering	
	\caption{Convergence results of the Eulerian velocity and pressure for three-dimensional circular Couette flow considered in Sec.~\ref{subsec:TCF}.
    }\vspace{4pt}
	\label{table:error_3D_TCF_Eulerian}	
\begin{tabular}{cc|c c | c c| c c|l}
\cline{3-8}
& & \multicolumn{2}{ |c| }{ velocity} & \multicolumn{2}{ |c| }{ pressure } & \multicolumn{2}{ |c| }{ pressure on $\Omega^*$ } \\  \cline{2-8}
&  \multicolumn{1}{ |c| }{h} & Error & Order & Error & Order & Error & Order & \\ \cline{1-8}
\multicolumn{1}{ |c  }{\multirow{4}{*}{$\L2$  } } &
\multicolumn{1}{ |c| }{$1.2500\times 10^{-1}$} & $7.9945\times 10^{-2}$ & - & $1.6719\times 10^{-1}$ & - & - & - &    \\ \cline{2-8}
\multicolumn{1}{ |c  }{}                        &
\multicolumn{1}{ |c| }{$6.2500\times 10^{-2}$} & $2.1306\times 10^{-2}$ & 1.91 & $5.0511\times 10^{-2}$ & 1.73 & - & - &    \\ \cline{2-8}
\multicolumn{1}{ |c  }{}                        &
\multicolumn{1}{ |c| }{$3.1250\times 10^{-2}$} & $7.3705\times 10^{-3}$ & 1.53 & $1.1295\times 10^{-2}$ & 2.16 & - & - &    \\ \cline{2-8}
\multicolumn{1}{ |c  }{}                        &
\multicolumn{1}{ |c| }{$1.5625\times 10^{-2}$} & $1.9527\times 10^{-3}$ & 1.92 & $3.8657\times 10^{-3}$ & 1.55 & - & - &    \\ \hhline{========}
\multicolumn{1}{ |c  }{\multirow{4}{*}{$\Linf$ } } &
\multicolumn{1}{ |c| }{$1.2500\times 10^{-1}$} &  $1.0378\times 10^{-1}$ & - & $4.6053\times 10^{-2}$ & - & $3.6172\times 10^{-2}$ & - & \\ \cline{2-8}
\multicolumn{1}{ |c  }{}                        &
\multicolumn{1}{ |c| }{$6.2500\times 10^{-2}$} &  $4.3354\times 10^{-2}$ & 1.26 & $4.3662\times 10^{-2}$ & 0.077 & $1.8408\times 10^{-2}$ & 0.97 &\\ \cline{2-8}
\multicolumn{1}{ |c  }{}                        &
\multicolumn{1}{ |c| }{$3.1250\times 10^{-2}$} &  $2.0026\times 10^{-2}$ & 1.11 & $3.1003\times 10^{-2}$ & 0.49 & $9.4395\times 10^{-3}$ & 0.96 &    \\ \cline{2-8}
\multicolumn{1}{ |c  }{}                        &
\multicolumn{1}{ |c| }{$1.5625\times 10^{-2}$} &  $6.7114\times 10^{-3}$ & 1.58 & $2.0668\times 10^{-2}$ & 0.58 & $4.7551\times 10^{-3}$ & 0.99 &    \\ \cline{1-8}
\end{tabular}
\end{table}
\begin{table}[ht]
    \scriptsize
	\centering	
	\caption{Convergence results of the displacement, velocity, pressure and WSS on the surface mesh of the three-dimensional circular Couette flow considered in Sec.~\ref{subsec:TCF}.
	  }\vspace{4pt}
	\label{table:error_3D_TCF_Lagrangian}	
\begin{tabular}{cc|c c |c c | c c| c c|l}
\cline{3-10}
& & \multicolumn{2}{ |c| }{displacement} & \multicolumn{2}{ |c| }{ velocity} & \multicolumn{2}{ |c| }{ pressure } & \multicolumn{2}{ |c| }{ WSS } \\  \cline{2-10}
&  \multicolumn{1}{ |c| }{h} & Error & Order & Error & Order & Error & Order & Error & Order & \\ \cline{1-10}
\multicolumn{1}{ |c  }{\multirow{4}{*}{$\L2$  } } &
\multicolumn{1}{ |c| }{$1.2500\times 10^{-1}$} & $6.8526\times 10^{-3}$ & - & $8.5642\times 10^{-2}$ & - & $1.16276\times 10^{0}$ & - & $5.8980\times 10^{-2}$ & - &    \\ \cline{2-10}
\multicolumn{1}{ |c  }{}                        &
\multicolumn{1}{ |c| }{$6.2500\times 10^{-2}$} & $1.7098\times 10^{-3}$ & 2.00 & $2.6461\times 10^{-2}$ & 1.69 & $2.8960\times 10^{-1}$ & 2.01 & $3.1341\times 10^{-2}$ & 0.91 &    \\ \cline{2-10}
\multicolumn{1}{ |c  }{}                        &
\multicolumn{1}{ |c| }{$3.1250\times 10^{-2}$} & $4.0606\times 10^{-4}$ & 2.07 & $7.3760\times 10^{-3}$ & 1.84 & $7.2821\times 10^{-2}$ & 1.99 & $1.7820\times 10^{-2}$ & 0.81 &    \\ \cline{2-10}
\multicolumn{1}{ |c  }{}                        &
\multicolumn{1}{ |c| }{$1.5625\times 10^{-2}$} & $1.0023\times 10^{-4}$ & 2.02 & $1.8712\times 10^{-3}$ & 1.98 & $2.5551\times 10^{-2}$ & 1.51 & $9.4697\times 10^{-3}$ & 0.91 &    \\ \hhline{==========}
\multicolumn{1}{ |c  }{\multirow{4}{*}{$\Linf$ } } &
\multicolumn{1}{ |c| }{$1.2500\times 10^{-1}$} &  $3.2012\times 10^{-3}$ & - & $3.7025\times 10^{-2}$ & -   & $3.3877\times 10^{-1}$ & - & $1.6843\times 10^{-2}$ & - & \\ \cline{2-10}
\multicolumn{1}{ |c  }{}                        &
\multicolumn{1}{ |c| }{$6.2500\times 10^{-2}$} &  $6.7975\times 10^{-4}$ & 2.24 &  $1.1325\times 10^{-2}$ & 1.71 & $8.9525\times 10^{-2}$ & 1.92 & $9.7928\times 10^{-3}$ & 0.78 &\\ \cline{2-10}
\multicolumn{1}{ |c  }{}                        &
\multicolumn{1}{ |c| }{$3.1250\times 10^{-2}$} &  $1.5146\times 10^{-4}$ & 2.17 & $3.8613\times 10^{-3}$ & 1.55 & $3.1973\times 10^{-2}$ & 1.49 & $5.3490\times 10^{-3}$ & 0.87 &    \\ \cline{2-10}
\multicolumn{1}{ |c  }{}                        &
\multicolumn{1}{ |c| }{$1.5625\times 10^{-2}$} &  $3.4150\times 10^{-5}$ & 2.15 & $9.3525\times 10^{-4}$ & 2.04 & $1.5225\times 10^{-2}$ & 1.07 & $2.8636\times 10^{-3}$ & 0.90 &    \\ \cline{1-10}
\end{tabular}
\end{table}

\end{document}